\documentclass[12pt, a4paper, USenglish]{scrartcl}
\def\mastercall{}

\ifx\preambleloaded\undefined
   \ifx\mastercall\undefined
\documentclass[12pt, a4paper, USenglish]{scrartcl}
\fi

\usepackage[USenglish]{babel}

\usepackage{geometry}
\usepackage[applemac]{inputenc}

\usepackage{ifthen}

\usepackage{dsfont}
\usepackage[T1]{fontenc}
\usepackage{verbatim}
\usepackage{cite}

\usepackage{graphicx}
\graphicspath{ {images/} }

\usepackage{amsmath, amscd, amsthm}
\usepackage{thmtools}
\usepackage{amssymb}
\usepackage{mathtools}

\usepackage{stmaryrd}

\usepackage{tikz-cd}
\usepackage{tikz}

\usepackage{setspace}
\usepackage{enumitem}
\usepackage{fancyhdr}
\usepackage{footnote}

\usepackage{tocloft}
\usepackage[titletoc]{appendix}				%
\usepackage[nottoc,notlot,notlof]{tocbibind} 	
\addto\captionsenglish{
  \renewcommand{\contentsname}%
    {Table of contents}%
}
\usepackage{abstract} 
\newtheoremstyle{newtheorem}{}{}{}{}{\bfseries}

\newcommand{\theoremend}{\ensuremath{\heartsuit}}
\newcommand{\definitionend}{\ensuremath{\clubsuit}}
\newcommand{\exampleend}{\ensuremath{\diamondsuit}}
\newcommand{\questionend}{\ensuremath{\spadesuit}}
\newcommand{\proofend}{\ensuremath{\square}}

\theoremstyle{definition}

\newtheorem{Masterthm}{Masterthm}[subsection]
\declaretheorem[name=Theorem,style=definition,qed={\theoremend},sibling=Masterthm]{Thm}
\declaretheorem[name=Goal,style=definition,qed={\questionend},sibling=Masterthm]{Goal}

\declaretheorem[name=Objection,style=definition,qed={\questionend},sibling=Masterthm]{Obj}
\declaretheorem[name=Question,style=definition,qed={\questionend},sibling=Masterthm]{Qstn}
\declaretheorem[name=Proof,style=proof,qed={\proofend}, numbered=no]{Prf}
\declaretheorem[name=Proposition,style=definition,qed={\theoremend},sibling=Masterthm]{Prop}
\declaretheorem[name=Claim,style=definition,qed={\theoremend},sibling=Masterthm]{Claim}
\declaretheorem[name=Lemma,style=definition,qed={\theoremend},sibling=Masterthm]{Lem}
\declaretheorem[name=Corollary,style=definition,qed={\proofend},sibling=Masterthm]{dCor}
\declaretheorem[name=Definition,style=definition,qed={\definitionend},sibling=Masterthm]{Def}
\declaretheorem[name=Construction ,style=definition,qed={\definitionend},sibling=Masterthm]{Cstr}
\declaretheorem[name=Theorem \& Definition,style=definition,qed={\theoremend\definitionend},sibling=Masterthm]{ThmDef}
\declaretheorem[name=Proposition \& Definition,style=definition,qed={\theoremend\definitionend},sibling=Masterthm]{PropDef}

\theoremstyle{remark}
\declaretheorem[name=Remark ,style=remark,qed={\exampleend},sibling=Masterthm]{Rem}
\declaretheorem[name=Warning ,style=remark,qed={\exampleend},sibling=Masterthm]{War}
\declaretheorem[name=Example ,style=remark,qed={\exampleend},sibling=Masterthm]{Expl}

\makeatletter				%
\@addtoreset{Thm}{section}	%
\@addtoreset{Rem}{section}	
\makeatother				%

\numberwithin{equation}{section}

\newcommand{\Proofof}[1]{of #1}

\usepackage{hyperref}

\newcommand{\qquote}[1]{``#1''}
\newcommand{\roughly}[1]{``#1''}

\newcommand{\ifisemptythenelse}[3]{
  \if\relax\detokenize{#1}\relax
    #2%
  \else
    #3%
  \fi
}
\newcommand{\maybebrackets}[1]{
				\ifisemptythenelse{#1}{}{{(#1)}}
				}
				
\newcommand{\introduce}{\textbf}
\newcommand{\buzzword}{\emph}
\newcommand{\axiomlabelstyle}[1]{(\bfseries{#1}\arabic*)} 
\newcommand{\axiomlabel}[2]{(\bfseries{#1}#2)} 
\newcommand{\axiomrefstyle}[1]{(#1\arabic*)} 
\newcommand{\axiomref}[2]{({#1}#2)} 

\newcommand{\mathblackboardboldnew}{\mathds}

\newcommand{\BC}{\mathblackboardboldnew C}

\newcommand{\BF}{\mathblackboardboldnew F}

\newcommand{\BN}{\mathblackboardboldnew N}

\newcommand{\BP}{\mathblackboardboldnew P}
\newcommand{\BQ}{\mathblackboardboldnew Q}

\newcommand{\BZ}{\mathblackboardboldnew Z}
\newcommand{\calA}{\mathcal{A}}
\newcommand{\calB}{\mathcal{B}}
\newcommand{\calC}{\mathcal{C}}
\newcommand{\calD}{\mathcal{D}}
\newcommand{\calE}{\mathcal{E}}
\newcommand{\calF}{\mathcal{F}}
\newcommand{\calG}{\mathcal{G}}
\newcommand{\calH}{\mathcal{H}}

\newcommand{\calL}{\mathcal{L}}
\newcommand{\calM}{\mathcal{M}}

\newcommand{\calO}{\mathcal{O}}
\newcommand{\calP}{\mathcal{P}}
\newcommand{\calQ}{\mathcal{Q}}

\newcommand{\calU}{\mathcal{U}}
\newcommand{\calV}{\mathcal{V}}
\newcommand{\calW}{\mathcal{W}}
\newcommand{\calX}{\mathcal{X}}
\newcommand{\calY}{\mathcal{Y}}

\newcommand{\frakA}{\mathfrak{A}}

\newcommand{\frakF}{\mathfrak{F}}

\newcommand{\frakH}{\mathfrak{H}}

\newcommand{\frakg}{\mathfrak{g}}

\DeclareMathOperator{\colim}{colim}
\DeclareMathOperator{\Hom}{Hom}
\DeclareMathOperator{\Ext}{Ext}

\DeclareMathOperator{\End}{End}
\DeclareMathOperator{\Mat}{Mat}
\DeclareMathOperator{\Aut}{Aut}
\DeclareMathOperator{\GL}{GL}

\DeclareMathOperator{\Sym}{Sym}

\DeclareMathOperator{\Rep}{Rep}
\DeclareMathOperator{\Res}{Res}
\DeclareMathOperator{\Ind}{Ind}

\DeclareMathOperator{\Ker}{Ker}
\DeclareMathOperator{\Cok}{Coker}
\DeclareMathOperator{\Stab}{Stab}
\DeclareMathOperator{\Dom}{Dom}

\DeclareMathOperator{\Iso}{Iso}
\renewcommand{\Im}{\operatorname{Im}}

\DeclareMathOperator{\Quot}{Quot}

\DeclareMathOperator{\rk}{rk}

\DeclareMathOperator{\B}{B}
\DeclareMathOperator{\K}{K}

\renewcommand{\hom}[1]{{\mathop{\kern0pt #1}}}
\newcommand{\Id}{\mathrm{id}}
\newcommand{\One}{\mathds{1}} 
\newcommand{\ev}{\mathrm{ev}}
\newcommand{\pr}{\mathrm{pr}}

\newcommand{\tr}{\mathrm{tr}}
\newcommand{\op}{\mathrm{op}}

\newcommand{\lax}{\mathrm{lax}}
\newcommand{\ab}{\mathrm{ab}}

\newcommand{\ch}{\mathrm{ch}}
\newcommand{\suml}{\sum\limits}
\newcommand{\prodl}{\prod\limits}

\renewcommand{\(}{\left(}
\renewcommand{\)}{\right)}
\newcommand{\ol}{\overline}

\renewcommand{\cal}{\mathcal}

\newcommand{\LHS}{\mathrm{LHS}}
\newcommand{\RHS}{\mathrm{RHS}}

\newcommand{\Fone}{{\BF_1}}
\newcommand{\Fq}{{\BF_q}}

\newcommand{\SG}[1]{{S_{#1}}}
\newcommand{\Sn}{\SG n}
\newcommand{\Sm}{\SG m}

\newcommand{\lreg}[2]{{_{#1}{#2}[{#1}]}}
\newcommand{\dual}{{\vee}}
\newcommand{\unit}{{1}}
\newcommand{\triv}{{\{\unit\}}}
\newcommand{\inv}{{-1}}
\newcommand{\units}{\times}
\newcommand{\fact}[1]{#1!}
\newcommand{\shriek}{!}
\newcommand{\finite}{\mathrm{finite}}
\newcommand{\actson}{\curvearrowright}
\newcommand{\actiondot}{.}
\newcommand{\lcaction}{\triangleright}
\newcommand{\rcaction}{\triangleleft}
\newcommand{\blank}{-}
\newcommand{\?}{{?}}
\newcommand{\multiset}[2]{\left(\binom{#1}{#2}\right)}
\newcommand{\Bigmid}{\,\,\Big|\,\,}

\newcommand{\multipl}[2]{{[#1:#2]}}
\newcommand{\terminal}{\{\star\}}
\newcommand{\initial}{\varnothing}
\newcommand{\htimes}{\overset h\times}
\newcommand{\conv}{\mathbin\ast }
\newcommand{\dotcup}{\mathbin{\dot\cup}}
\newcommand{\pair}[2]{\(\ifisemptythenelse{#1}\blank{#1};\ifisemptythenelse{#2}\blank{#2}\)}
\newcommand{\apair}[2]{\langle\ifisemptythenelse{#1}\blank{#1};\ifisemptythenelse{#2}\blank{#2}\rangle}
\newcommand{\bpair}[3]{\(\ifisemptythenelse{#1}\blank{#1};\ifisemptythenelse{#2}\blank{#2};\ifisemptythenelse{#3}\blank{#3}\)}
\newcommand{\overcat}[2]{(#1\downarrow#2)}
\newcommand{\spanby}[1]{\langle{#1}\rangle}

 \newcommand\rquot[2]{
        \mathchoice
            {
                \text{\raise1ex\hbox{${#1}$}}\Big/\lower1ex\hbox{${#2}$}%
            }
            {
                #1\,/\,#2
            }
            {
                #1/#2
            }
            {
                #1/#2
            }
    }
\newcommand{\doubleslash}{\mathbin{
  \mathchoice{\Big/\mkern-10mu\Big/}
    {/\mkern-6mu/}
    {/\mkern-5mu/}
    {/\mkern-5mu/}}}
\newcommand{\qquot}[2]{
\mathchoice
	    {
                \text{\raise1ex\hbox{${#1}$}}\doubleslash\lower1ex\hbox{${#2}$}%
            }
	{{#1}\doubleslash{#2}}
	{{#1}\doubleslash{#2}}
	{{#1}\doubleslash{#2}}
} 
    
\newcommand\lquot[2]{
        \mathchoice
            {
                \lower1ex\hbox{${#2}$}\Big\backslash\text{\raise1ex\hbox{${#1}$}}%
            }
            {
                #2\setminus#1
            }
            {
                #2\setminus#1
            }
            {
                #1\setminus#2
            }
    }
    
 \newcommand\lrquot[3]{
        \mathchoice
            {
                \lower1ex\hbox{${#2}$}\Big\backslash\text{\raise1ex\hbox{${#1}$}}\Big/\lower1ex\hbox{${#3}$}%
            }
            {
                #2\backslash#1/#3
            }
            {
                #2\backslash#1/#2
            }
            {
                #1\backslash#2/#2
            }
    }

\newcommand{\category}[1]{\mathrm{\mathbf{#1}}}
\newcommand{\lMod}[1]{#1\cat{-Mod}}
\newcommand{\lmod}[1]{#1\cat{-mod}}

\newcommand{\lfree}[1]{#1\cat{-free}}
\newcommand{\lrep}[2]{{#1}\cat{-rep}_{#2}}
\newcommand{\lRep}[2]{{#1}\cat{-Rep}_{#2}}

\newcommand{\cat}[1]{\mathrm{\mathbf{#1}}}
\newcommand{\Set}{\cat{Set}}

\newcommand{\Grpd}{\cat{Grpd}}
\newcommand{\Cat}{\cat{Cat}}
\newcommand{\CAT}{\cat{CAT}}
\newcommand{\Vect}[1]{\cat{Vect}_{#1}}
\newcommand{\vect}[1]{\cat{vect}_{#1}}

\newcommand{\vectC}{\vect\BC}
\newcommand{\Mor}[3]{#1\left(#2,#3\right)}
\newcommand{\mors}{\to}

\newcommand\restr[3][]{{
  \left.\kern-\nulldelimiterspace 
  {#2} 
  \vphantom{\big|} 
  \right|_{#3}^{#1} 
  }}

\newlength{\longarrow}
\newlength{\arrow}
\newlength{\blablaone}
\newlength{\blablatwo}

\newcommand{\rsqa}{\rightsquigarrow}
\newcommand{\lra}{\longrightarrow}
\newcommand{\lla}{\longleftarrow}
\newcommand{\la}{\leftarrow}
\newcommand{\ra}{\rightarrow}
\newcommand{\La}{\Leftarrow}
\newcommand{\Ra}{\Rightarrow}
\newcommand{\LRa}{\Leftrightarrow}

\newcommand{\lRa}{\Longrightarrow}
\newcommand{\lLRa}{\Longleftrightarrow}
\newcommand{\hla}{\hookleftarrow}
\newcommand{\hra}{\hookrightarrow}
\newcommand{\thla}{\twoheadleftarrow}
\newcommand{\thra}{\twoheadrightarrow}

\newcommand{\xRa}{\xRightarrow}
\newcommand{\xLa}{\xLeftarrow}

\newcommand{\xhra}[2][]{
	\settowidth{\blablaone}{\scriptsize$#1$}
	\settowidth{\blablatwo}{\scriptsize$#2$}
	\pgfmathsetlength{\arrow}{max(\blablaone,\blablatwo,\longarrow)}
	\xhookrightarrow[{\mathmakebox[\arrow]{#1}}]{\mathmakebox[\arrow]{#2}}%
}
\newcommand{\xra}[2][]{
	\settowidth{\blablaone}{\scriptsize$#1$}
	\settowidth{\blablatwo}{\scriptsize$#2$}
	\pgfmathsetlength{\arrow}{max(\blablaone,\blablatwo,\longarrow)}
	\xrightarrow[{\mathmakebox[\arrow]{#1}}]{\mathmakebox[\arrow]{#2}}%
}
\newcommand{\xla}[2][]{
	\settowidth{\blablaone}{\scriptsize$#1$}
	\settowidth{\blablatwo}{\scriptsize$#2$}
	\pgfmathsetlength{\arrow}{max(\blablaone,\blablatwo,\longarrow)}
	\xleftarrow[{\mathmakebox[\arrow]{#1}}]{\mathmakebox[\arrow]{#2}}%
}
\newcommand{\xlra}[2][]{
	\settowidth{\blablaone}{\scriptsize$#1$}
	\settowidth{\blablatwo}{\scriptsize$#2$}
	\pgfmathsetlength{\arrow}{max(\blablaone,\blablatwo,\longarrow)}
	\xleftrightarrow[{\mathmakebox[\arrow]{#1}}]{\mathmakebox[\arrow]{#2}}%
}

\newcommand{\xthra}[2][]{\xra[#1]{#2}\!\!\!\!\!\ra}

\let\amsamp=& 

\author{Tashi Walde}
\title{{Hall} monoidal categories and categorical modules}
\date{}

\newcommand{\itemsecondlayer}{\ensuremath{\ast }}

\newcommand{\hasadj}{\dashv}

\newcommand{\Groth}[2]{{#1\rtimes#2}}
\newcommand{\mon}[2][]{{\mathcal{#2}^\boxtimes_{#1}}}

\newcommand{\monF}[1][]{\mon[#1]F}
\newcommand{\isoA}[1][]{{\calA_{#1}^{\cong}}}

\newcommand{\finitary}{\mathrm{f}}

\newcommand{\Hecke}{\frakH} 
\newcommand{\hall}[2][]{\mathrm{hall}_{#1}\maybebrackets{#2}}		
\newcommand{\fHall}[2][]{\mathrm{Hall}_{#1}\maybebrackets{#2}}		
\newcommand{\Hall}[1]{\mathrm{Hall}^\boxtimes\maybebrackets{#1}}			
\newcommand{\HallV}[1]{\mathrm{Hall}_\calV^\boxtimes\maybebrackets{#1}}			
\newcommand{\Hallf}[1]{\mathrm{Hall}_\finitary^\boxtimes\maybebrackets{#1}}	

\DeclareMathOperator{\Pol}{Pol}
\DeclareMathOperator{\Polone}{Pol_k^{(1)}}
\newcommand{\on}{{\otimes n}}
\newcommand{\twofib}[3][]{{{#2}\twotimes{#1}\{#3\}}}
\newcommand{\twotimes}[1]{\mathbin{\overset{(2)}\times_{#1}}}
\newcommand{\Ffree}[1]{\lfree{\Fone[{#1}]}}
\newcommand{\FGfree}{\Ffree G}
\newcommand{\cm}{\kappa}
\newcommand{\cmreg}{r}
\newcommand{\cmtriv}{t}

\newcommand{\otV}[1]{\bigotimes_\calV\maybebrackets{#1}}
\newcommand{\fotV}[2][]{\bigotimes\nolimits_{#2}^{#1}}
\newcommand{\pV}[1]{\calV^\star\maybebrackets{#1}}
\newcommand{\sV}[1]{\calV_\shriek\maybebrackets{#1}}
\newcommand{\VX}[1]{\calV(\calX_{\{#1\}})}
\newcommand{\pVX}[1]{\pV{\calX_{#1}}}
\newcommand{\sVX}[1]{\sV{\calX_{#1}}}
\newcommand{\sVY}[1]{\sV{\calY_{#1}}}

\newcommand{\xX}[1]{\X_{\{{#1}\}}}
\newcommand{\sX}[1]{\X_{\{#1\}}}

\newcommand{\hypo}{hypo} 
\newcommand{\Hypo}{Hypo} 
\newcommand{\Dop}{{\Delta^{\op}}}
\newcommand{\Deltap}{\Delta_\S}
\newcommand{\Dpop}{{\Delta_\S^{\op}}}
\newcommand{\Dnop}{{\Delta_\mindecor^{\op}}}
\newcommand{\Dmop}{{\Delta_\maxdecor^{\op}}}
\newcommand{\Dnmop}{{\Delta_\mmdecor^{\op}}}
\newcommand{\minimal}{{-\infty}}
\newcommand{\maximal}{{+\infty}}
\newcommand{\mindecor}{\geqslant}
\newcommand{\maxdecor}{\leqslant}
\newcommand{\mmdecor}{\diamond}

\newcommand{\XAQ}{\calX^{\calA,\calQ}}
\newcommand{\XAQb}[1]{\calX^{\calA,\calQ}_{\{#1\}}}

\newcommand{\n}{{[n]}}
\newcommand{\m}{{[m]}}
\renewcommand{\l}{{[l]}}
\newcommand{\0}{{[0]}}
\newcommand{\1}{{[1]}}
\newcommand{\2}{{[2]}}
\newcommand{\enu}[1]{{\{0,\ldots,{#1}\}}}

\newcommand{\X}{\calX}
\newcommand{\V}{\calV}

\newcommand{\C}{\calC}

\newcommand{\cococo}{componentwise cocontinuous }

\author{Tashi Walde}

\begin{document}

\pagestyle{fancy}
\renewcommand{\sectionmark}[1]{\markboth{#1}{}}

\fancyhf{}
\pagenumbering{arabic}

\fancyhead[L]{\ifthenelse{\isodd{\value{page}}}{\thepage}{}}
\fancyhead[R]{\ifthenelse{\isodd{\value{page}}}{\leftmark}{\thepage}}

\renewcommand{\thefootnote}{\arabic{footnote})}		

\def\preambleloaded{}	
\fi

\pagenumbering{roman}
\maketitle 
\begin{abstract}
\ifx\preambleloaded\undefined
   
   \def\isstandalone{}
   
\fi

We construct so called \buzzword{Hall monoidal categories} (and \buzzword{Hall modules} thereover) and exhibit them as a categorification of classical Hall and Hecke algebras (and certain modules thereover). The input of the (functorial!) construction are simplicial groupoids satisfying the $2$-Segal conditions (as introduced by Dyckerhoff and Kapranov~\cite{DyckerhoffKapranovHigherSegalSpacesI}), the main examples come from Waldhausen's S-construction. To treat the case of modules, we introduce a relative version of the $2$-Segal conditions.

Furthermore, we generalize a classical result about the representation theory of symmetric groups to the case of wreath product groups: We construct a monoidal equivalence between the category of complex $G\wr\Sn$-representations (for a fixed finite group $G$ and varying $n\in\BN$) and the category of \roughly{$G$-equivariant} polynomial functors; we use this equivalence to prove a version of Schur-Weyl duality for wreath products.\\

This paper is, up to minor modifications, the author's Master's thesis as submitted to the University of Bonn on July 22, 2016.

\ifx\isstandalone\undefined
\else
\newpage
\bibliographystyle{amsalpha}
\bibliography{literatur}
\end{document}
\fi

\end{abstract}

\newpage
\tableofcontents



\newpage\pagenumbering{arabic}
\section{Introduction}
\ifx\preambleloaded\undefined
   
   \def\isstandalone{}
\fi

The theory of Hall algebras began more than a century ago when Steinitz observed~\cite{SteinitzZurTheorieDerAbelschenGruppen} that (for each fixed prime $p$) there is an associative (!) product 
\[A\cdot B\coloneqq \sum_Ct(C,A,B)\cdot C\]
on (formal sums of) isomorphism classes $\{A,B,C,\dots\}$ of finite abelian ($p$-)groups, where $t(C,A,B)$ counts the number of subgroups $A'\subseteq  C$ such that $A'\cong A$ and $\rquot C{A'}\cong B$. This product yields an associative algebra with a basis consisting of partitions by identifying each finite abelian $p$-group $\bigoplus_{\lambda}\rquot{\BZ}{(p^{\lambda_i})}$ with its type $\lambda$.\\
The topic remained forgotten for more than fifty years until Hall rediscovered what is nowadays known as \buzzword{Hall's algebra of partitions}~\cite{HallTheAlgebraOfPartitions}. This most classical of Hall algebras was and is of great interest due to its close relationship to several fundamental objects in mathematics such as symmetric functions~\cite{MacDonaldSymmetricFunctionsHallPolynomials} and flag varieties.\\
Hall algebras came back into the spotlight in the early 90s due to Ringel's groundbreaking discovery~\cite{RingelHallAlgebrasAndQuantumGroups} that the positive part of the quantized enveloping algebra $\calU_v(\frakg)$ (of a simple complex Lie algebra $\frakg$) can be realized by applying Hall's construction to the category $\lrep{\overrightarrow Q}\Fq$ of finite dimensional representations over $\Fq$ of the corresponding Dynkin quiver.\\
Many variants of Hall algebras have since been introduced and studied; see Schiffmann's lecture notes \cite{SchiffmannLecturesOnHallAlgebras}\cite{SchiffmannLecturesOnCanonicalAndCrystalBasesOfHallAlgebras} for a detailed overview.

\subsection{A universal perspective}
Dyckerhoff and Kapranov~\cite{DyckerhoffKapranovHigherSegalSpacesI} propose the following perspective on Hall algebras:\\
The various collections
\[\calX^\calA_\n\coloneqq \{0= A_0\hra A_1\hra\dots\hra A_n\}\]
of flags of length $n\in \BN$ in a certain suitable category\footnote{%
In the classical case of Steinitz and Hall, for instance, $\calA$ would be the category of finite abelian $p$-groups. For us $\calA$ will always be a \buzzword{proto-abelian} category~\cite[Definition 1.2]{DyckerhoffHigherCategoricalAspectsOfHallAlgebras}.
}$\calA$ can be naturally organized into a simplicial object\footnote{For us $\calX^\calA$ will always be a simplicial \emph{groupoid}, i.e.\ $\calC\coloneqq \Grpd$. Depending on the precise nature of the category $\calA$ other types of simplicial objects (e.g.\ simplicial spaces or simplicial stacks) can be used~\cite{DyckerhoffHigherCategoricalAspectsOfHallAlgebras}.}
\[\calX^\calA\colon \Dop\lra \calC\]
which is known in algebraic K-theory as \buzzword{Waldhausen's S-construction}~\cite{WaldhausenAlgebraicKTheoryOfSpaces}.\\
We are supposed to view $\calX^\calA$ as the \roughly{universal} Hall algebra; all other Hall algebras are created out of $\calX^\calA$ by applying appropriate specializations $\calC\rsqa \?$ (called \buzzword{transfer theories}~\cite[Section 8.1]{DyckerhoffKapranovHigherSegalSpacesI}).
\begin{Expl}[see Section~\ref{classicalHallAlg}]~\label{introclassicalRingelHall}The classical Ringel-Hall algebra $\hall\calA$ (of a \buzzword{proto-abelian} category $\calA$) can be recovered~\cite[Proposition 2.19]{DyckerhoffHigherCategoricalAspectsOfHallAlgebras} by taking $\calC\coloneqq \Grpd$ to be the category of groupoids and by passing to (finitely supported) functions on isomorphism classes.
\end{Expl}
Associativity and unitality of the various Hall algebras are encoded universally in the simplicial object $\calX^\calA$ as the so called \buzzword{$2$-Segal conditions}\cite[Section 2.3]{DyckerhoffKapranovHigherSegalSpacesI}.\\ 
The $2$-Segal conditions can be formulated for any simplicial object; it is therefore tempting to consider other simplicial objects $\calX\colon \Dop\to\calC$ which do not necessarily arise through Waldhausen's S-construction.
\begin{Expl}[see Section~\ref{convolutionalgebras}]\label{introclassicalHecke}The classical Hecke algebra $\Hecke(G,H)$ associated to an inclusion $H\subseteq G$ of finite groups can be obtained~\cite[Example 8.2.11]{DyckerhoffKapranovHigherSegalSpacesI} by taking the $2$-Segal simplicial groupoid (called the \buzzword{Hecke-Waldhausen construction})
\[\calX^{G,H}\colon \n\longmapsto \B H\twotimes{\B G}\cdots\twotimes{\B G}\B H\]
(with $n+1$ many factors) and passing to (finitely supported) functions on isomorphism classes.
\end{Expl}

The Ringel-Hall algebra and the Hecke algebra can thus be put into the same context: both algebras can be obtained by applying the transfer theory 
\[V\colon \Grpd\to \lMod\BZ,\quad V\colon A\mapsto \BZ^{(\pi_0A)}\]
(of finitely supported functions on isomorphism classes) to an appropriate $2$-Segal simplicial groupoid $\calX$; in both cases we simply speak of the Hall algebra $\hall\calX$ of the corresponding simplicial groupoid $\calX$.

\subsection{Hall monoidal categories}
When reading about functions on isomorphism classes of a groupoid the following objection should immediately jump to mind:
\begin{Obj}\cite{DyckerhoffHigherCategoricalAspectsOfHallAlgebras} If we are to consider \roughly{maps} out of a groupoid $A$ then surely we should not be using \emph{functions} $\pi_0A\to\BZ$ on isomorphism classes but rather \emph{functors} $A\to\calV$ for some suitably chosen category $\calV$.\end{Obj}

This idea leads to a categorification of the Hall algebra.
\begin{Thm}[see Section~\ref{sectionmonconstruction}]\label{IntroHallMonCat} If $\calV$ is a cocomplete monoidal category (such that the monoidal product is cocontinuous in each variable) then $[\blank,\calV]\colon \Grpd\to\CAT$ induces a functor
\[\{\text{$2$-Segal simplicial groupoids}\}\xra{\HallV{}}\{\text{monoidal categories $\mid$ lax monoidal functors}\}\]
such that $\HallV\calX$ is a monoidal structure on the functor category $[\calX_\1,\calV]$ for every $2$-Segal simplicial groupoid $\calX\colon \Dop\to\Grpd$.
\end{Thm}

The construction in Theorem~\ref{IntroHallMonCat} relies on the fact~\cite[Section 1.1]{LurieDerivedAlgebraicGeometryII} that monoidal categories can be described as \buzzword{op-fibrations of categories} (see Appendix~\ref{AppendixGrothendieckConstruction}) over $\Dop$ satisfying the \buzzword{pointed $1$-Segal condition} (see Section~\ref{sectionmonoidalobjects}).\\

Observe how Theorem~\ref{IntroHallMonCat} upgrades the Hall \emph{algebra} to a \emph{monoidal category} and upgrades the \emph{set-map}
\[\{\text{$2$-Segal simplicial groupoids}\footnote{Fine print: satisfying some suitable finiteness conditions}\}\xra{\hall{}}\{\text{$\BZ$-algebras}\}\]
to a \emph{functor}.\\
In the case $\calV\coloneqq \Vect\BC$ and under suitable finiteness conditions it is also possible to go the other way (\buzzword{decategorify}) and turn the Hall monoidal category $\Hall\calX$ back into an algebra by first restricting to a smaller monoidal category $\Hallf\calX$ and then passing to the Grothendieck-ring $\K_0$. We will study the relationship between this \buzzword{fat Hall algebra}\footnote{%
We call it fat because it is significantly bigger than the Hall algebra itself: every basis element of the Hall algebra (corresponding to an isomorphism class $[a]$ in $\calX_\1$) gets replaced by the whole representation theory of the group $\Aut(a)$.
} $\K_0\Hallf\calX$ and the Hall algebra $\hall\calX$ (see Section~\ref{sectionfinitarycase}).
\begin{Rem}
The first example of a Hall monoidal category appearing in the literature is probably due to Joyal and Street. They introduced~\cite{JoyalStreetTheCategoryOfRepresentationsOfTheGeneralLinearGroupOverAFiniteField} what we call the \buzzword{finitary Hall monoidal category} of $\vect\Fq$ and even equipped it with a braiding.
\end{Rem}

\begin{Rem}We prove a slightly more general version of Theorem~\ref{IntroHallMonCat} by replacing the assignment $[\blank ,\calV]\colon \Grpd\to\CAT$ by an arbitrary \buzzword{\cococo left monoidal derivator\footnote{Warning: There are various definitions for the term \qquote{derivator}, starting from the original one by Grothendieck~\cite[\S I]{GrothendieckLesDerivateurs}. We use an ad hoc version of the one found in Groth's thesis~\cite{GrothOnTheTheoryOfDerivators} and we sweep most subtleties of monoidal derivators under the rug. Therefore some additional care needs to be taken if one wants to use the actual monoidal derivators occurring in the literature.}
of groupoids} (see Section~\ref{sectionmonoidalderivators}). The rich supply of monoidal derivators coming from monoidal model categories~\cite{GrothMonoidalDerivatorsAndAdditiveDerivators} should therefore provide new (and hopefully interesting) Hall monoidal categories.
\end{Rem}

\subsection{Can we also construct \roughly{Hall modules}?}
We think of a $2$-Segal simplicial object as a \roughly{mold} out of which unital associative structures (e.g.\ unital algebras, monoidal categories,...) can be formed by passing to specializations.
\begin{Qstn}Can we define a similar mold which gives rise to \emph{modules over} associative structures?\end{Qstn}
Dyckerhoff and Kapranov suggest~\cite{DyckerhoffKapranovHigherSegalSpacesII} that we consider morphisms ${\calY\to\calX}$ of simplicial objects and impose some suitable \roughly{relative} variant of the $2$-Segal conditions. We should then think of $\calY$ as an abstract $\calX$-module which gives rise to an actual module after specialization.
\begin{Qstn}\label{questionaboutmodulesandmorphisms}
What is the appropriate notion of a $2$-Segal morphism of simplicial objects? How do $2$-Segal morphisms give rise to modules?
\end{Qstn}

To answer these questions we introduce the language of \buzzword{\hypo-simplicial objects}.
\begin{Def}
A \introduce{\hypo-simplicial object} (in $\calC$) is a functor $\calX\colon \Dpop\to\calC$ where $\Deltap\subseteq\Delta$ is an \buzzword{admissible subcategory of $\Delta$} (see Section~\ref{sectionsubsimplicialobjects}).
\end{Def}
It is straightforward to generalize the $2$-Segal conditions from simplicial objects to \hypo-simplicial objects (see Definition~\ref{Definitionunital}).\\
The main observation which ties the theory of \hypo-simplicial objects to Question~\ref{questionaboutmodulesandmorphisms} is that a morphism $\Phi\colon \calY\to\calX$ of simplicial objects (in $\calC$) is the same thing as a \hypo-simplicial object $\Phi^\mindecor\colon \Dnop\to \calC$ (a \buzzword{relative simplicial object}) defined on a certain admissible subcategory $\Delta_\mindecor\subset\Delta$ (see Section~\ref{sectionrelativesimplicialobjects}).
\begin{Def}\label{introDefrel2Segal}
A morphism $\Phi\colon \calY\to\calX$ of simplicial objects is called \introduce{$2$-Segal} if the corresponding relative simplicial object $\Phi^\mindecor\colon \Dnop\to\calC$ is $2$-Segal.
\end{Def}
\begin{Rem}
The notion of relative $2$-Segal objects was independently discovered by Young~\cite{YoungRelative2SegalSpaces}. Proposition~\ref{criterionrelativeunital} shows that Young's definition is equivalent to ours.
\end{Rem}

\begin{Expl} [see Section~\ref{sectionboundedflags}]
We introduce a \roughly{bounded} version of Waldhausen's S-construction by considering the groupoids $\calY^{\calA,\calQ}_\n$ of those flags
\[0=A_0\hra A_1\hra \dots \hra A_n\hra A_\maximal\]
(in a proto-abelian category $\calA$) where all the quotients $\rquot{A_\maximal}{A_i}$ (and the morphisms between them) are contained in a suitable fixed subcategory $\calQ\subseteq\calA$ (a \buzzword{quotient datum}\footnote{%
Quotient data themselves are plentiful: Each pair $(X,S)$ consisting of an object $X\in\calA$ and a group $S\subseteq \Aut(X)$ of automorphisms gives rise to a quotient datum $\calQ_S(X)$, by restricting (roughly) to those morphisms induced by $S$ on quotients of $X$ in $\calA$.
}). These groupoids $\calY^{\calA,\calQ}_\n$ can naturally be assembled into a $2$-Segal morphism $\calY^{\calA,\calQ}\to\calX^\calA$ of simplicial groupoids.
\end{Expl}
\begin{Expl}[see Section~\ref{convolutionalgebras}]
For inclusions $H\subseteq G\supseteq P$ of finite groups, the projections
\[\calX_\n^{G,H;P}\coloneqq \calX^{G,H}_\n\twotimes{\B G}{\B P} \lra\calX^{G,H}_\n\]
assemble to a $2$-Segal morphism $\calX^{G,H;P}\to\calX^{G,H}$ of simplicial groupoids. We will see that the resulting module over the Hecke algebra $\Hecke(G,H)$ (after passing to finitely supported functions on isomorphism classes) is just the usual convolution module on $H$-$P$-biinvariant functions $G\to\BZ$.
\end{Expl}

Using this new language we can prove a refined version of Theorem~\ref{IntroHallMonCat}:
\begin{Thm}[see Section~\ref{sectionmonconstruction}]\label{IntroTheoremGenHallconstr}
Let $\Deltap\subseteq\Delta$ be an admissible subcategory. Every \cococo left monoidal derivator of groupoids $\calV\colon \Grpd\to\Cat$ induces a functor
\[\{\text{$2$-Segal $\Dpop\to\Grpd$}\}\xra{\HallV{}}\{\text{pointed $1$-Segal op-fibrations over $\Dpop$}\}^\lax\qedhere\]
\end{Thm}
We recover \autoref{IntroHallMonCat} because a pointed $1$-Segal op-fibration over $\Dop$ is nothing but a monoidal category. Similarly a pointed $1$-Segal op-fibration over $\Dnop$ is the same thing as a categorical right action of a monoidal category on a category (see Section~\ref{sectionmonoidalobjects}).

\begin{dCor}
Every $2$-Segal morphism $\Phi\colon \calY\to\calX$ of simplicial groupoids gives rise to a \buzzword{categorical right Hall module} $\HallV{\Phi^\mindecor}$ over the Hall monoidal category $\HallV\calX$.
\end{dCor}

Every $2$-Segal morphism $\Phi$ gives also rise to a similar categorical \emph{left} Hall module $\HallV{\Phi^\maxdecor}$ which has the same underlying category as $\HallV{\Phi^\mindecor}$. By using a certain admissible subcategory $\Delta_\mmdecor\subset\Delta$ we will see how these two modules are \roughly{dual} to each other (see Section~\ref{sectiondualitypairing}).

\subsection{A Schur-Weyl duality for wreath products via \roughly{G-equivariant} polynomial functors}
Most of the tools developed in this thesis will be tested on our pet example, the category $\FGfree$ of free representations of a finite group $G$ over $\Fone$ (see Section~\ref{HallalgebraofFGfree} and Section~\ref{HallFGfree}).\\
The Hall algebra $\hall\FGfree$ of the category $\FGfree$ is canonically isomorphic to the ring of divided powers. This ring is quite well understood. The same can not be said about the (finitary) Hall monoidal category $\Hallf\FGfree$ which turns out to be the category
\[\lrep{G\wr\SG\ast}\BC\coloneqq \bigoplus_{n\in\BN}\lrep{G\wr \Sn}\BC\]
equipped with the induction product (see Section~\ref{HallFGfree}).\\

One of the goals of this thesis is to understand the monoidal category $\lrep{G\wr\SG\ast}\BC$ in terms of \buzzword{polynomial functors} $\lrep G\BC\to\vect\BC$ and to use this understanding to obtain a new perspective on the representation theory of wreath products\footnote{%
The \buzzword{wreath product groups} $G\wr \Sn\coloneqq G^n\rtimes\Sn$ are generalizations of the symmetric group with a similar (but more complicated) representation theory~\cite[\S I.B]{MacDonaldSymmetricFunctionsHallPolynomials}. A special case is the Weyl group of type $B$ which is given as $\frac{\BZ}{(2)}\wr \Sn$.
}.
\begin{Thm}[see Section~\ref{sectionpolyfunctors}]\label{introEquivalencewreathprodpolyfunctors}
The formula
\[{G\wr\Sn}\actson E \longmapsto E\otimes_{G\wr\Sn}T^n(\blank)\]
(where $T^n(\blank)\colon \lrep G\BC\to\lrep{G\wr \Sn}\BC$ is the $n$-th tensor power functor) defines an equivalence of monoidal categories
\[\Xi\colon \lrep{G\wr \SG\ast}\BC\xra\simeq \{\text{polynomial functors}\footnote{Fine print: of bounded total degree}\text{ $\lrep G\BC\to\vect\BC$}\}.\]
If $G$ is abelian then $\Xi$ categorifies the classical ring isomorphism~\cite[\S I.B]{MacDonaldSymmetricFunctionsHallPolynomials}
\[\ch\colon \K_0(\lrep{G\wr \SG\ast}\BC)\xra{\cong}\Lambda(G)\]
to the ring of symmetric polynomial class functions.
\end{Thm}
\begin{Rem}Theorem~\ref{introEquivalencewreathprodpolyfunctors} makes precise the intuition that polynomial \emph{functors} $\lrep G\BC\to\vect\BC$ should be a categorification of $G$-equivariant (with respect to conjugation) symmetric polynomial \emph{functions} $\BC[G]^d\to\BC$ (with $d\to\infty$). In this sense we can think of polynomial functors $\lrep G\BC\to\vect\BC$ as a \roughly{$G$-equivariant} analogue of the classical polynomial functors $\vect\BC\to\vect\BC$.\end{Rem}

We can use Theorem~\ref{introEquivalencewreathprodpolyfunctors} to obtain the following Schur-Weyl duality for wreath products.
\begin{Thm}[Schur-Weyl duality for wreath products; see Section~\ref{sectionSchurWeylduality}]\label{introSchurWeylduality}
Let $G$ be a finite abelian group and let $V$ be a free $G$-representation (over $\BC$) of rank $d$.
Then $\Xi$ (from \autoref{introEquivalencewreathprodpolyfunctors}) induces an assignment (with $0\mapsto0$)
\begin{equation}\label{introSchurWeylassignment}
	\begin{tikzcd}
	{\{\text{irred.\ $G\wr \Sn$-reps (over $\BC$) (up to iso)}\}}\cup{\{0\}}\dar{\Phi}\\
	{\{\text{irred.\ polyn.\ $\Aut_G(V)$-reps.\ (over $\BC$) of homogeneous deg.\  $n$ (up to iso)}\}}\cup {\{0\}}
	\end{tikzcd}
\end{equation}
which is always surjective and is always injective outside of its kernel.\\
The assignment $\Phi$ is a bijection in the case $n\leqslant d$.
\end{Thm}

\begin{Rem}
The topic of polynomial $\GL_n$-representations and polynomial functors is an old one which appears in the literature in many different variants and flavors; see for instance~\cite{GreenPolynomialRepresentationsOfGLnWithAppendix},~\cite{FriedlanderSuslinCohomologyOfFiniteGroupSchemesOverAField},~\cite{KockNotesOnPolynomialFunctors}, and references therein.\\
Theorem~\ref{introEquivalencewreathprodpolyfunctors} and Theorem~\ref{introSchurWeylduality} reduce to very classical statements~\cite[\S I.A]{MacDonaldSymmetricFunctionsHallPolynomials} for $G=\triv$; most ideas in this case go back to Schur~\cite{SchurUeberEineKlasseVonMatrizen, SchurUeberDieRationalenDarstellungenDerAllgemeinenLinearenGruppe}.\\
After the completion of this thesis I was made aware that Theorem~\ref{introEquivalencewreathprodpolyfunctors} is proved implicitly (with essentially the same arguments) in a paper by Macdonald~\cite{MacDonaldPolynomialFunctorsAndWreathProducts}; our proof should therefore be seen as a repackaging of his ideas. A similar result was proved by Gal and Gal~\cite[Section 7]{GalGalSymmetricSelf-AdjointHopfCategoriesAndACategoricalHeisenbergDouble} by using a categorified version of Zelevinsky's theory of positive self-adjoint Hopf algebras.\\
What might be of interest in our proof of Theorem~\ref{introSchurWeylduality} is the approach via \roughly{G-equivariant} polynomial functors. Various other approaches to a Schur-Weyl duality for wreath products are known; see for instance~\cite{MazorchukStropppelGlkdModulesViaGroupoids} and references therein.
\end{Rem}

\subsection{Acknowledgements}
First and foremost I thank my advisor, Tobias Dyckerhoff, who first introduced me to the beautiful topic of Hall algebras and who was my main source of ideas, feedback and guidance.
Further I am very thankful to Catharina Stroppel for many inspiring discussions.
My thanks also go to Aras Ergus and Malte Leip for correcting many of my countless mistakes.\\
I am grateful to the \emph{DAAD} and the \emph{Studienstiftung des deutschen Volkes} who provided me with financial stability during my studies at the University of Bonn. Last, but certainly not least, I want to thank those dear people who gave me emotional support when I had trouble persevering.

\subsection{Notations and conventions}

\begin{itemize}
\item The symbols $\BN$, $\BZ$, $\BQ$, $\BC$ denote the natural numbers (including $0$), integers, rational numbers and complex numbers respectively.
\item The cardinality of a finite set (or group) $X$ is denoted by $\#X$ or by $|X|$.
\item The notations $\subseteq$ and $\subset$ are synonymous, but we preferably use the latter if equality clearly does not occur or is not interesting (e.g.\ $\BZ\subset\BQ$).
\item All rings and algebras are associative with $1$ but are usually not commutative.
\item The categories of \emph{all} vector spaces (over $k$), representations (of $G$ over $k$) and (left)-modules (over a ring $R$) are denoted by $\Vect k$, $\lRep G k$ and $\lMod R$ respectively. When considering \emph{finite dimensional} vector spaces/representations or \emph{finitely generated} modules we use lower-case letters (e.g.\ $\vect k$, $\lrep G k$ or \ $\lmod R$).
\item We write $a\in\calC$ if $a$ is an object in the category $\calC$.
\item We denote the initial and terminal object in some category by $\initial$ and $\terminal$ respectively.
\item The opposite of a category $\calC$ is denoted either by $\calC^\op$ or by $\calC_\op$ depending on the typographical needs.
\item If $f\colon a\ra b$ is a morphism (from $a$ to $b$) in $\calC$ and we view it as a morphism (from $b$ to $a$) in $\calC^\op$ then we use the notation $f\colon b\la a$; the arrow will always show the direction in the original category $\calC$.
\item The category of functors $\calC\to\calB$ (and natural transformations between them) is denoted by $[\calC,\calB]$ both in the enriched and in the non-enriched setting.
\item If $G$ is a monoid (resp.\ $R$-algebra) then we sometimes think of $G$ as the category (resp.\ $R$-linear category) $\B G$ with one element $\star_G$ and $\End_{\B G}(\star_G)\coloneqq G$.
\item The symbol $\Delta$ can denote both the simplex category (of finite linearly ordered sets) or a diagonal morphism (e.g.\ $\Delta\colon X\to X\times X$); the context will remove any ambiguity. See Section~\ref{preliminariessimpl} for more notation and conventions about $\Delta$.
\item If $n\in \BN$ is a natural number then the symbol $\n$ denotes the linearly ordered set $\{0<\dots<n\}$ (usually seen as an object in $\Delta$).
\item If $G$ is a group then $k[G]$ denotes the group algebra and $\lreg Gk$ denotes the regular representation of $G$ (over $k$).
\item For a finite group $G$ we denote by $G_\star$ the set conjugacy of classes in $G$ and by $G^\star$ the set of irreducible characters of $G$.
\item The size of a partition $\lambda=(\lambda_1\geqslant\lambda_2\geqslant\dots)$ is $|\lambda|\coloneqq \sum_i\lambda_i$. We denote by $\calP_n$ the set of partitions of size $n$.
\item For a set $X$ we denote by $\calP_n(X)$ the set of partition-valued maps on $X$ with total size $n$; in other words, an element of $\calP_n(X)$ is a family $\(\lambda(x)\)_{x\in X}$ of partitions such that $\|\lambda\|\coloneqq \sum_x |\lambda(x)|=n$.
\item For a partition $\lambda$ we denote by $s_\lambda$, $e_\lambda$ and $h_\lambda$ the corresponding Schur function, elementary symmetric function and complete symmetric function respectively.
\item The tuple $(1,\dots,1)$ with $d$ many $1$'s is denoted by $(1^d)$. If the $d$ is clear from the context we might just write $(1)$.
\end{itemize}

\ifx\isstandalone\undefined
\else
\newpage
\bibliographystyle{amsalpha}
\bibliography{literatur}
\end{document}
\fi


\newpage\section{The classical Hall algebra}\label{classicalHallAlg}
\ifx\preambleloaded\undefined
   
   \def\isstandalone{}
\fi

We start off by recalling the classical definition of the Hall algebra.

Let $\calA$ be a proto-abelian category~\cite[Definition 1.2]{DyckerhoffHigherCategoricalAspectsOfHallAlgebras}. In particular $\calA$ has a zero object, direct sums, as well as a notion of flags and short exact sequences.
\begin{Def}
A proto-abelian category $\calA$ is called \introduce{finitary} if it is essentially small and for each pair of objects $A,A'\in \calA$ the sets $\Hom_\calA(A,A')$ and $\Ext_\calA(A,A')$ (of extensions $0\to A'\to ?\to A\to0$ up to equivalence) have finite cardinality.
\end{Def}
\begin{Expl}
The prototypical example of a finitary proto-abelian category is the category $\vect\Fq$ of finite dimensional vector spaces over the finite field $\Fq$; more generally we can take the category of all finite abelian $p$-groups for some prime $p$.\\
Another example is the category $\vect\Fone$ of finite pointed sets and partial bijections which we study in more detail in Section~\ref{HallFGfree}.
\end{Expl}
\begin{War}
A proto-abelian category does not need to be enriched over abelian groups as the example $\vect\Fone$ shows.
\end{War}

\begin{Cstr}The \introduce{classical Hall algebra} $\hall\calA$ of a proto-abelian category $\calA$ is given as follows:
\begin{itemize}
\item The underlying abelian group is the free $\BZ$-module
\[\bigoplus\limits_{[M]\in \Iso(\calA)}\BZ\cdot[M]\]
on the set of isomorphism classes of objects in $\calA$.
\item The product of two basis elements is defined as 
\[[N]\cdot [L]\coloneqq \sum_{[M]\in \Iso}g_{N,L}^M[M]\]
with the coefficients $g_{N,L}^M$ given by
\[g_{N,L}^M\coloneqq \frac{\#\{\text{s.e.s. } L\hra M\thra N\}}{\#\Aut(L)\#\Aut(N)}\]
\end{itemize}
One can prove that $\hall\calA$ is an associative algebra~\cite[Theorem 1.6]{DyckerhoffHigherCategoricalAspectsOfHallAlgebras}. We denote the extension of scalars $R\otimes_\BZ\hall\calA$ of the Hall algebra by any ring $R$ by $\hall[R]\calA$.
\end{Cstr}

Following Dyckerhoff and Kapranov~\cite[Section 8.2]{DyckerhoffKapranovHigherSegalSpacesI} we subdivide the construction $\calA\leadsto \hall\calA$ of the Hall algebra into two steps:
\begin{enumerate}
\item From the proto-abelian category $\calA$ we can construct the simplicial groupoid $\calX^\calA\colon \Dop\to\Grpd$ of flags.
\item To each groupoid $A$ we can assign an abelian group $V(A)$ of finitely supported functions $\pi_0A\to \BZ$.
\end{enumerate}
We can then redefine the Hall algebra by taking $V(\calX^\calA_\1)$ as the underlying abelian group; the multiplicative structure will be given by the extra functoriality structure of $V$ and the higher groupoids $\calX^\calA_\n$.

\subsection{The simplicial groupoid of flags}\label{defgroupoidofflags}
Given a protoabelian category $\calA$ we can define the following simplicial groupoid $\X^\calA\colon \Dop\to \Grpd$:
\begin{itemize}
	\item For $\n\in \Dop$, the objects of the groupoid $\X^\calA_\n$ are diagrams of the shape
	\begin{equation}\label{flagtrianglen}
		\begin{tikzcd}
			0\ar[r, hookrightarrow]&	A_{01}\ar[r, hookrightarrow]\ar[rd, "{\square}"description, phantom]\ar[d, twoheadrightarrow]&	A_{02}\ar[r, hookrightarrow]\ar[d, twoheadrightarrow]&	\cdots\ar[r, hookrightarrow]&	A_{0,n-2}\ar[rd, "{\square}"description, phantom]\ar[d, twoheadrightarrow]\ar[r, hookrightarrow]&	A_{0,n-1}\ar[rd, "{\square}"description, phantom]\ar[d, twoheadrightarrow]\ar[r, hookrightarrow]	&A_{0,n}\ar[d, twoheadrightarrow]\\
			&	0\ar[r, hookrightarrow]&	A_{12}\ar[d, twoheadrightarrow]\ar[r, hookrightarrow]&	\cdots\ar[r, hookrightarrow]&	A_{1,n-2}\ar[rd, "{\square}"description, phantom]\ar[d, twoheadrightarrow]\ar[r, hookrightarrow]&	A_{1,n-1}\ar[rd, "{\square}"description, phantom]\ar[d, twoheadrightarrow]\ar[r, hookrightarrow]&	A_{1,n}\ar[d, twoheadrightarrow]\ar[d, twoheadrightarrow]\\ 
			&&	0\ar[r, hookrightarrow]&	\cdots\ar[r, hookrightarrow]&	A_{2,n-2}\ar[d, twoheadrightarrow]\ar[r, hookrightarrow]&	A_{2,n-1}\ar[d, twoheadrightarrow]\ar[r, hookrightarrow]&	A_{2,n}\ar[d, twoheadrightarrow]\ar[d, twoheadrightarrow]\\ 
			&&&\ddots&\vdots\ar[d, twoheadrightarrow]&\vdots\ar[d, twoheadrightarrow]&\vdots\ar[d, twoheadrightarrow]\ar[d, twoheadrightarrow]\\
			&&&&0\ar[r, hookrightarrow]&A_{n-2, n-1}\ar[rd, "{\square}"description, phantom]\ar[r, hookrightarrow]\ar[d, twoheadrightarrow]&A_{n-2,n}\ar[d, twoheadrightarrow]		\\	
			&&&&&0\ar[r, hookrightarrow]&A_{n-1,n}
		\end{tikzcd}
	\end{equation}
	where all the vertical maps are epis, all the horizontal maps are monos and all the square are bicartesian (i.e.\ both pullback and pushout).
	The morphisms of $\X^\calA_n$ are simply all isomorphisms in $\calA$ of such diagrams.
	\item A morphism $f\colon \n\la\m$ in $\Dop$ induces a morphism $\X^\calA_f\colon \X^\calA_\n\to\X^\calA_\m$ of groupoids by simultaneously omitting (in the case of face maps) and duplicating (in the case of degeneracy maps) rows and columns of such diagrams. It is easy to see that this construction is functorial.
\end{itemize}
We call $\calX^\calA$ the \introduce{simplicial groupoid of flags} in $\calA$.

\begin{Rem}
Clearly the groupoid $\X^\calA_\1$ is canonically isomorphic to the groupoid $\isoA$ obtained from $\calA$ by throwing away all non-isomorphisms. More generally, $\X^\calA_\n$ is equivalent to the groupoid
\[\{0\hra A_1\hra A_2\hra\dots\hra A_n\}^{\cong}\]
of flags of length $n$ in $\calA$. This is due to the fact that Diagram~\ref{flagtrianglen} is determined up to isomorphism by the first row and the fact that all squares are pushouts.
\end{Rem}

\begin{Rem}
The simplicial object of flags $\calX^\calA$ first appeared in the context of algebraic K-theory where it is known as the \buzzword{Waldhausen S-construction}~\cite[Section 1.3]{WaldhausenAlgebraicKTheoryOfSpaces}.
\end{Rem}

\subsection{Functions on isomorphism classes}

Consider the assignment $V\colon \Grpd\to \lMod\BZ$ given by sending a groupoid $A$ to the free $\BZ$-module $V(A)\coloneqq \BZ^{(\pi_0A)}$ of finitely supported functions on the isomorphism classes of $A$. This assignment can not be extended to a functor but only to a partial functor. We just give a sketch, see Dyckerhoff's lecture notes~\cite[Section~2.3]{DyckerhoffHigherCategoricalAspectsOfHallAlgebras} for more details and complete proofs.
\begin{Def}
	Let $f\colon A\to B$ be a map of groupoids.
	\begin{itemize}
		\item The map $f$ is called \introduce{$\pi_0$-finite} if the induced map of sets $\pi_0A\to\pi_0B$ has finite fibers.\\
		In this case the \introduce{pullback} $f^\star\colon  V(B)\to V(A)$ given by 
		\[f^\star\colon \varphi\mapsto (a\mapsto \varphi(f(a)))\]
		 is a well defined $\BZ$-module homomorphism.
		 \item We call $f$ \introduce{locally finite} if for every $a\in A$ and $b\in B$ the $2$-fiber $\twofib[{f\restriction A(a)}]{A(a)}b$ (where $A(a)$ denotes the connected component of $a$ in $A$) is finite, i.e.\ has finite automorphism groups and finitely many isomorphism classes.\\
		 In this case we can define the \introduce{pushforward} $f_\shriek\colon \BQ\otimes_\BZ V(A)\to \BQ\otimes_\BZ V(B)$ via
		 \begin{equation}\label{definepushforward}
			(f_\shriek\psi)(b) \coloneqq \int_{\twofib[f]Ab}\psi\coloneqq \sum\limits_{[x, f(x)\xra{\cong}b]\in \pi_0\(\twofib[f]Ab\)}\frac{\psi(x)}{\#\Aut([x, f(x)\xra{\cong}b])}
		 \end{equation}
		 for a function $\psi\colon \pi_0A\to\BQ$ of finite support.\qedhere
	\end{itemize}
\end{Def}

Sometimes we can define the pushforward without having to tensor with $\BQ$.
\begin{Lem}
Let $f\colon A\to B$ be faithful in addition to being locally finite. Then the pushforward $f_\shriek\colon \BQ\otimes_\BZ V(A)\to \BQ\otimes_\BZ V(B)$ restricts to a $\BZ$-module homomorphism $f_\shriek\colon V(A)\to V(B)$.
\end{Lem}
\begin{Prf}
If $f$ is faithful then all the $2$-fibers $\twofib[f]Ab$ are discrete, hence the denominators in Formula~\ref{definepushforward} disappear.
\end{Prf}

Using the fact that $2$-pullbacks are invariant under natural isomorphisms it is easy to see that both the pushforward and the pullback construction are invariant under natural isomorphisms, i.e.\ $f^\star=g^\star$ and $f_\shriek=g_\shriek$ if $f\cong g$ (whenever defined). Further, the obvious functoriality properties $(fg)^\star=g^\star f^\star$ and $(fg)_\shriek=f_\shriek g_\shriek$ can be shown to hold whenever both sides are defined. Finally for each $2$-pullback square of groupoids
	\[\begin{tikzcd}[cramped]
		X\rar{f'}\dar{g'}& B\dar{g}\\
		A\rar{f}&D\\
	\end{tikzcd}\]
we have the pull-push equation $(f')_\shriek\circ (g')^\star= g^\star\circ f_\shriek$ whenever both sides are defined (actually the left side will be automatically defined if the right side is defined).

\begin{Expl}
For a group homomorphism $f\colon \B H\to \B G$ the $2$-fiber $\twofib[f]{\B H}{\star_{\B G}}$ is nothing but the action groupoid $\qquot GH$ which has objects labeled by $g\in G$ and morphisms $g\xra{h} gh$ for $h\in H$.\\
Moreover $f$ is locally finite if and only if the index $\multipl{G}{\Im f}$ and the cardinality $\#\Ker f$ are finite. In this case the pushforward
\[f_\shriek\colon \BQ\cong {\BQ\otimes_\BZ V(\B H)}\lra \BQ\otimes_\BZ V(\B G)\cong \BQ\]
is given by the number $\frac{\multipl{G}{\Im f}}{\#\Ker f}$ which is just $\frac{\#G}{\#H}$ if both groups are finite. If $f$ is faithful then this number is just the integer $[G\colon H]$.
\end{Expl}

\subsection{The classical Hall algebra, revised}\label{classHallalgrevisited}
Let $\calA$ be a finitary proto-abelian category. As outlined before, we can use the pull-push-structure on $V$ and the groupoid $\calX^\calA_\2$ to define the multiplication on $V(\calX^\calA_\1)$.
Consider the \introduce{multiplication span}
\begin{equation}\label{multiplicationspan}
	\calX_\1^\calA\times\calX_\1^\calA\xla{c_2}\calX_\2^\calA\xra{\nu_2} \calX_\1^\calA
\end{equation}
induced by the inclusions $\{0,1\},\{1,2\}\hra\{0,1,2\}\hla\{0,2\}$. We call $c_2$ the \introduce{chopping map} and $\nu_2$ the \introduce{extremal map}.

\begin{Prop}\cite[Proposition 2.19]{DyckerhoffHigherCategoricalAspectsOfHallAlgebras}\label{HallAlgsrevisited}
	\begin{enumerate}
		\item The chopping map $c_2$ in the span~\ref{multiplicationspan} is $\pi_0$-finite.
		\item The extremal map $\nu_2$ in the span~\ref{multiplicationspan} is locally finite and faithful.
		\item The induced map
\[\mu\colon V(\calX_\1^\calA)\times V(\calX_\1^\calA)\xra{\cdot} V(\calX_\1^\calA\times\calX_\1^\calA)\xra{c_2^\star}V(\calX_\2^\calA)\xra{(\nu_2)_\shriek}V(\calX_\1^\calA),\]
endows the abelian group $V(\calX_\1^\calA)$ with the structure of an associative and unital $\BZ$-algebra.
		\item The canonical isomorphism $\hall\calA\xra{\cong}V(\isoA)\cong V(\calX_\1^\calA)$  (which sends an isomorphism class $[M]$ to the function $\delta_{[M]}$) of abelian groups is compatible with the multiplications, hence we can (and will) identify the classical Hall algebra $\hall\calA$ with $\(V(\calX_\1^\calA),\mu\)$.\qedhere
	\end{enumerate}
\end{Prop}

\subsection{Example: vector spaces and free group representations over $\Fone$}\label{HallalgebraofFGfree}
Let $\vect\Fone$ be the category of finite pointed sets where the morphisms are partial bijections, i.e.\ maps $f\colon X\to Y$ of finite pointed sets such that the restriction $\restr f{X\setminus f^\inv\{0\}}$ (of $f$ to the set of elements it doesn't kill) is injective. We always denote the distinguished element in a set by $0$.\\
More generally, if $G$ is a finite group we can consider the category $\FGfree$ of free $G$-representations over $\Fone$. The objects of this category are pointed sets with a free $G$-action (by $\Fone$-linear isomorphisms, i.e.\ pointed bijections) and the morphisms are $G$-equivariant partial bijections. Of course we recover $\vect\Fone$ as $\Ffree\triv$.\\
It is straightforward to show that $\FGfree$ is proto-abelian and finitary with kernels and cokernels given by $\Ker f=f^\inv\{0\}$ and $\Cok (f\colon X\to Y)= \rquot{Y}{\Im f}\coloneqq \{0\}\cup \(Y\setminus \Im f\)$ respectively. The direct sum in $\FGfree$ is given by the wedge sum of pointed sets with the induced $G$-action.

We denote by $\calX$ the simplicial groupoid of flags in $\FGfree$.\\

Each object $X\in \FGfree$ can be decomposed into $G$-orbits as 
\[X=\bigoplus_{\{0\}\neq G\actiondot x\in \lquot XG} \langle x\rangle\]
where all the $\langle x\rangle\coloneqq \{0\}\cup G\actiondot x$ are isomorphic via $x\mapsto 1_G$ (and $0\mapsto 0$) to the regular representation $\lreg G\Fone\coloneqq \{0\}\cup G$. Hence the only isomorphism-invariant in the category $\FGfree$ is the \introduce{rank} defined by $\rk X\coloneqq \frac{|X\setminus\{0\}|}{\#G}$.\\ 
If $X$ is of rank $n$ and $\{0\}\dot\cup\{x_1,\dots, x_n\}$ is a system of representatives of $\lquot {X}G$ (in this case we call $\{x_1,\dots,x_n\}$ a \introduce{$\Fone[G]$-basis} of $X$) then we have an isomorphisms of groups
\begin{equation}\label{describeautomgroupofn}
	G\wr \Sn\xra{\cong} \Aut (X)
\end{equation}
given by $(g_1,\dots g_n, \sigma)\actiondot x_i\coloneqq g_{\sigma(i)}\actiondot x_{\sigma(i)}$.\\
We conclude that the assignment $\star_{G\wr \Sn}\mapsto \spanby n\coloneqq \bigoplus_{i=1}^n \lreg G{\Fone}$ on objects extends to an equivalence of groupoids
\begin{equation}\label{groupoidequivalencetoSn}
	\coprod\limits_{n\in \BN}\B{G\wr \Sn}\xra{\simeq}\calX_\1\simeq {\FGfree}^{\cong}
\end{equation}
by the group isomorphism \ref{describeautomgroupofn} for the standard basis of $\spanby{n}$.\\
Therefore the underlying abelian group of $\hall{\FGfree}$ has a basis given by the stalk functions $\delta_\n$ for $n\in \BN$.

\begin{Claim}\label{multiplicationspanFG}
	The equivalence~\ref{groupoidequivalencetoSn} extends to an equivalence of spans
	\begin{equation}\label{describespan\Sn}
		\begin{tikzcd}
			\calX_\1^\calA\times\calX_\1^\calA&\calX_\2^\calA\rar{\nu_2}\ar[l,"{c_2}"']& \calX_\1^\calA\\
			\coprod\limits_{n\in \BN}\B{G\wr \Sn}\times \coprod\limits_{m\in \BN}\B{G\wr \Sm}\ar[u,"\simeq"]&\coprod\limits_{n,m\in \BN}\B{G\wr (\Sn\times \Sm)}\ar[l,"\cong"']\ar[r, ""]\ar[u,"\simeq"]&\coprod\limits_{l\in \BN}\B{G\wr \SG{l}}\ar[u,"\simeq"]
		\end{tikzcd}
	\end{equation}
	where the lower right horizontal map is given by the inclusions ${G\wr(\Sn\times \Sm)}\hra G\wr \SG{l}$ for $l=n+m$.
\end{Claim}
\begin{Prf}
	We define the middle vertical map by using the obvious maps
	\[G\wr(\Sn\times \Sm)\lra \Aut\Big(0\to\spanby n\hra\spanby{n+m}\thra\spanby m\to0\Big)\]
	clearly making Diagram~\ref{describespan\Sn} commutative. It is straightforward to check that we get an equivalence of groupoids by using the fact that every short exact sequence $0\to X\to Z\to Y\to 0$ is just a wedge-partition of $Z$ into two $G$-stable subsets $X$ and $Y$ with trivial intersection.
\end{Prf}

We are left to compute the multiplicative structure. To do this we fix $n,m,l\in \BN$ and compute $\mu(\delta_n,\delta_m)\spanby l$. By dropping the connected components on which $\delta_n\cdot\delta_m$ is zero and restricting to the fiber of $\star_{\B{G\wr \SG{l}}}$ (which corresponds to $\spanby l\in \calX_\1^\calA$), we are left with computing the pushforward of the function $\One\colon \B {G\wr {(\Sn\times \Sm)}}\to \BZ$ along the inclusion
${G\wr(\Sn\times \Sm)}\hra G\wr\SG l$ if $l=n+m$; we are left with nothing otherwise.

We conclude the identity
\[\delta_n\cdot\delta_m=\frac{|G\wr\SG {n+m}|}{|G\wr(\Sn\times \Sm)|}\delta_{n+m}=\binom{n+m}{n}\delta_{n+m}.\]
in $\hall{\FGfree}$.
\begin{dCor}\label{hallalgofFGfreedividedpowers}
	The $\BZ$-linear isomorphism
	\[\hall{\FGfree}\cong \BZ^{(\BN)}\xra{\cong}\Gamma_\BZ[x]\coloneqq \BZ\left\{\frac{x^n}{\fact n}\Bigmid n\in \BN\right\}\]
	given by $\delta_n\mapsto \frac{x^n}{\fact n}$ gives an algebra isomorphism from the Hall algebra of $\FGfree$ to the ring of divided powers.
\end{dCor}
We observe that the Hall algebra is independent of the group $G$. This can be seen as the first hint of the fact that the Hall algebra is not a sufficiently strong invariant and that we could maybe do better.

\begin{Rem}
Hall algebras of representation categories \roughly{over $\Fone$} have been studied in great detail by Szczesny. For instance, he considers (Dynkin) quivers~\cite{SzczensyRepresentationsOfQuiversOverF1} (where $\hall{\lrep Q\Fone}$ is related (but not necessarily isomorphic!) to the positive part of the universal enveloping algebra of the corresponding simple Lie algebra) and coherent sheaves on $\BP^1$~\cite{SzczesnyOnTheHallAlgebraOfCoherentSheavesOnP1OverF1}.
\end{Rem}

\ifx\isstandalone\undefined
\else
\newpage
\bibliographystyle{amsalpha}
\bibliography{literatur}
\end{document}
\fi


\newpage\section{Preliminaries on (\hypo-)simplicial objects}\label{preliminariessimpl}

\subsection{Notation and Conventions}
\ifx\preambleloaded\undefined
   
   \def\isstandalone{}
\fi

\begin{itemize}
\item We view the \introduce{simplex category} $\Delta$ as the category of finite linear orders; morphisms are weakly monotone maps.

\item When writing about objects in $\Delta$ we use $\{x_0,\dots,x_n\}$ and $\{x_0<\dots<x_n\}$ interchangeably.
\item We write $N^n\in \Delta$ if the set $N$ has $n+1$ elements (i.e.\ $N$ is of \introduce{dimension} $n$).
\item We call subset $M'\subseteq M$ (of an objects $M\in \Delta$) convex if for $x'<y<z'$ in $M$ with $x',z'\in M'$ it follows that $y\in M'$.
\item For a morphism $f\colon M\to N$ in $\Delta$ and a subset $M'\subseteq M$ we abbreviate $\ol f(M')\coloneqq \{x\in N\mid f(\min M')\leqslant x\leqslant f(\max M')\}$ for the \introduce{convex hull} of $f(M')$ in $N$.
\item For $\{x_0<\dots< x_n\}=N^n\in\Delta$  we define the subsets $N_i\coloneqq \{x_{i-1},x_i\}$ for all $1\leqslant i\leqslant n$.
\item A \introduce{simplicial object} in a category $\calC$ is a functor $\Dop\to\calC$.\\
A \introduce{pseudo-simplicial object} in a $2$-category $\calC$ is a pseudo-functor $\Dop\to\calC$.
\item If the objects in the category $\calC$ have a common name <name> then we say \qquote{(pseudo)-simplicial <name>} instead of \qquote{(pseudo)-simplicial object in $\calC$}. For instance we say \qquote{simplicial groupoid} for a functor $\Dop\to\Grpd$ or \qquote{pseudo-simplicial category} for a  pseudo-functor $\Dop\to\CAT$.

\end{itemize}
\begin{Rem}
When working in $\Delta$ (or $\Dop$) it is enough to consider the standard objects $\n\coloneqq \enu{n}$ since any other object $N^n\in \Delta$ of dimension $n$ is isomorphic to $\n$ via a \emph{unique} isomorphism.
\end{Rem}

\ifx\isstandalone\undefined
\else
\newpage
\bibliographystyle{amsalpha}
\bibliography{literatur}
\end{document}
\fi

\subsection{\Hypo-simplicial objects}\label{sectionsubsimplicialobjects}
\ifx\preambleloaded\undefined
   
   \def\isstandalone{}
\fi

The reason why we bloated our category $\Delta$ (which we could have defined on the skeleton of objects $\{\n\mid n\in\BN\}$ without losing anything) is that we want to consider certain subcategories of $\Delta$ and so we need some redundancy in the objects in order to have \roughly{enough space}.

\begin{Def}\label{admissiblecategory}
We call a subcategory $\Deltap\subseteq \Delta$ \introduce{admissible} if it satisfies the following conditions
\begin{enumerate}[label=\axiomlabelstyle{ad}, ref=\axiomlabelstyle{ad}]
\item If a morphism $f$ in $\Deltap$ is a bijection (i.e.\ $f$ is an isomorphism in $\Delta$) then the inverse of $f$ is also in $\Deltap$ (i.e.\ $f$ is an isomorphism in $\Deltap$).
\item The inclusion $M\hra N$ is a morphism of $\Deltap$ whenever $\{x_0,\dots,x_n\}=N^n\in\Deltap$ and $M\subseteq N$ is not one of the two extremal singletons $\{x_0\}$ and $\{x_n\}$.
\item Every morphism $f\colon M\to N$ in $\Deltap$ can be factorized (in $\Deltap$) as
\[f\colon M\thra f(M)\hra \ol f(M)\hra N\qedhere\]
\end{enumerate}
\end{Def}

\begin{Def}
 Let $\calC$ be a ($2$-)category. We call a (pseudo)-functor $\Dpop\to\calC$ defined on an admissible subcategory $\Deltap\subseteq \Delta$ a \introduce{(pseudo)-\hypo-simplicial object} in $\calC$.
\end{Def}

\begin{Expl}
The whole simplex category $\Delta\subseteq \Delta$ is clearly admissible. So we recover the notion of (pseudo)-simplicial objects by specialising $\Deltap\coloneqq \Delta\subseteq \Delta$.
\end{Expl}

\begin{Def}If $\calX\colon \Dpop\to\calC$ is a (pseudo-)simplicial object then we call the morphisms
$\xX{x_0,\dots x_n}\to \xX{x_0,\dots,x_{i-1},x_{i+1},\dots,x_n}$ induced by the inclusion \introduce{face maps} and the morphisms $\xX{x_0,\dots x_n}\to \xX{x_0,\dots,x_{i},x'_{i},\dots,x_n}$ induced by $x_i,x'_i\mapsto x_i$ \introduce{degeneracy maps} (if they exist).\end{Def}

\subsubsection{Relative simplicial objects}\label{sectionrelativesimplicialobjects}
We will now see the other main example which we will use later to construct modules over the Hall monoidal category.

\begin{Def}
Fix the formal symbol $\minimal$. The subcategory $\Delta_\mindecor\subset \Delta$ is defined as follows:
\begin{itemize}
\item Objects are those $\{\minimal\}\neq N\in\Delta$ such that $\minimal\in N$ implies that $\minimal$ is a minimal element of $N$.
\item Morphisms are weakly monotone maps $f$ such that for all $x\in\Dom f$ we have $f(x)=\minimal$ if and only if $x=\minimal$.
\end{itemize}
Similarly we can define (after fixing another formal symbol $\maximal$) the subcategory $\Delta_\maxdecor\subset\Delta$:
\begin{itemize}
\item Objects are those $\{\maximal\}\neq N\in\Delta$ such that $\maximal\in N$ implies that $\maximal$ is a maximal element of $N$.
\item Morphisms are weakly monotone maps $f$ such that for all $x\in\Dom f$ we have $f(x)=\maximal$ if and only if $x=\maximal$.
\end{itemize}
We call $\Delta_\mindecor$ and $\Delta_\maxdecor$ the \introduce{right-} and \introduce{left-relative simplex category} respectively.
\end{Def}
\begin{Rem}
Clearly both $\Delta_\mindecor$ and $\Delta_\maxdecor$ are admissible subcategories of $\Delta$.\\
For the rest of Section~\ref{preliminariessimpl} we will work with $\Delta_\mindecor$; everything we do works verbatim for $\Delta_\maxdecor$.
\end{Rem}

Note that $\Delta_\mindecor$ decomposes on objects as $\Delta_\mindecor=\Delta_{\in}\dotcup\Delta_{\not\in}$ where $\Delta_{\in}=\{N\in\Delta_\mindecor\mid \minimal\in N\}$ and $\Delta_{\not\in}=\{N\in\Delta_\mindecor\mid \minimal\not\in N\}$. Observe that both $\Delta_\in$ and $\Delta_{\not\in}$ are canonically equivalent to $\Delta$ via $N\mapsto N\setminus\{\minimal\}$ and $N\mapsto N$ respectively.\\
There are clearly no morphisms from $\Delta_\in$ to $\Delta_{\not\in}$ and the morphisms from $\Delta_{\not\in}\to\Delta_\in$ are generated by the inclusions $\iota_N\colon N\hra \{\minimal\}\cup N$ for $N\in \Delta_{\not\in}$, i.e.\
\[\Hom(\Delta_{\not\in},\Delta_\in)=\Hom(\Delta_\in,\Delta_\in)\{\iota_N\mid N\in \Delta_\mindecor\}\Hom(\Delta_{\not\in},\Delta_{\not\in}).\]

\begin{Prop}\label{relativeisgoodnotion}
Let $\calC$ be a category. We have an equivalence of categories
\begin{equation}
	\{\text{\hypo-simpl.\ objects $\Dnop\to\calC$}\} \lra \{\text{maps of simpl.\ objects $\Dop\to\calC$}\}
\end{equation}
sending a \hypo-simplicial object $\calX^\mindecor\colon \Dnop\to\calC$ to the natural transformation
\[\calX^\mindecor_{\iota_N}\colon \calX^\in_N\coloneqq \calX^\mindecor_{\{\minimal\}\sqcup N}\lra\calX^\mindecor_{N}\eqqcolon \calX^{\not\in}_N\]
between the two simplicial objects $\calX^{\in}\colon \Dop\simeq\Delta_{\in}^\op\xra{\calX^\mindecor}\calC$ and ${\calX^{\not\in}\colon \Dop\simeq\Delta_{\not\in}^\op\xra{\calX^\mindecor}\calC}$.
\end{Prop}

\begin{Rem}
It is not hard to see that $\Delta_\mindecor$ is equivalent to the Grothendieck construction (see Appendix~\ref{AppendixGrothendieckConstruction}) of the functor $[\Id]\colon \{\bullet\to\bullet\}\to\Cat$ which picks out the identity $\Id\colon \Delta\to\Delta$. In light of this description, Proposition~\ref{relativeisgoodnotion} becomes purely formal.
\end{Rem}

\begin{Prf}[\Proofof{Proposition~\ref{relativeisgoodnotion}}]
The inverse functor is constructed as follows:\\
Let $\Phi\colon \calP\to\calA$ be a morphism of simplicial objects $\calP,\calA\colon \Dop\to\calC$. Denote by $\hat\calP$ and $\hat\calA$ the simplicial objects $\Delta_{\in}^\op\simeq\Dop\xra{\calP}\calC$ and $\Delta_{\not\in}^\op\simeq\Dop\xra{\calA}\calC$ respectively; let $\hat\Phi\colon \hat\calP\to\hat\calA$ be the transformation induced by $\Phi$. We define the \hypo-simplicial object $\Phi^\mindecor\colon \Dnop\to\calC$ by $\Phi^\mindecor \equiv \hat\calP$ on $\Delta_{\in}^\op$, by $\Phi^\mindecor \equiv \hat\calA$ on $\Delta_{\not\in}^\op$ and by $\Phi^\mindecor_{\iota_N}\coloneqq \hat\Phi_N$. The functoriality of $\Phi^\mindecor$ corresponds to the naturality of $\Phi$ together with the functoriality of $\calP,\calA$.\\
It is clear that the two procedures described above are inverse to each other and it is easy to see that they are functorial.
\end{Prf}

\begin{Rem}
If $\calC$ is a $2$-category then we can replace all functors and natural transformations by pseudo-functors and pseudo-natural transformations and obtain an analogous statement for pseudo-(\hypo-)simplicial objects.
\end{Rem}

Proposition~\ref{relativeisgoodnotion} justifies the following definition:
\begin{Def}
We call (pseudo-)\hypo-simplicial objects defined on $\Dnop$ and $\Delta_{\maxdecor}^\op$ \introduce{right-} and \introduce{left-relative (pseudo-)simplicial objects} respectively.
\end{Def}

\ifx\isstandalone\undefined
\else
\newpage
\bibliographystyle{amsalpha}
\bibliography{literatur}
\end{document}
\fi

\subsection{The $1$-Segal condition}\label{sectionmonoidalobjects}
\ifx\preambleloaded\undefined
   
   \def\isstandalone{}
\fi

Let $\calC$ be a $2$-category with finite products and a terminal object $\terminal$.\\
In this section we denote by $\approxeq$ a map which admits an adjoint inverse (see Appendix~\ref{calculusofmates}). This condition is usually stronger than being an equivalence ($\simeq$) and weaker than being an isomorphism ($\cong$).
If a map $f$ admits an adjoint inverse, we fix one and denote it by $f^\inv$.\\

Let $\Deltap\subseteq\Delta$ be an admissibile subcategory.
\begin{Def} A (pseudo)-\hypo-simplicial object $\calX\colon \Dpop\to\calC$ is said to satisfy the \introduce{pointed $1$-Segal} condition if
\begin{enumerate}[label=\axiomlabelstyle{pS}, ref=\axiomrefstyle{pS}]
		\item\label{suboneSegal} for every $\{x_0<\dots<x_n\}=N^n\in\Dpop$ the chopping map
\[c_N\colon \calX_N\xra{\approxeq}\sX{x_0,x_1}\times\dots\times\sX{x_{n-1},x_n}\] 
induced by the inclusions $N_i\coloneqq \{x_{i-1},x_i\}\hra N$ for $1\leqslant i\leqslant n$ admits an adjoint inverse.
	\end{enumerate}
	If $\calC$ has some notion of homotopy fiber products (denoted by $\htimes$) then we can also define the \introduce{unpointed $1$-Segal} condition:
\begin{enumerate}[label=\axiomlabelstyle{uS}, ref=\axiomrefstyle{uS}]
		\item\label{unpointsuboneSegal} For every $\{x_0<\dots<x_n\}=N^n\in\Dpop$ the chopping map
\[c_N\colon \calX_N\xra{\approxeq}\sX{x_0,x_1}\htimes_{\xX{x_1}}\dots\htimes_{\xX{x_{n-1}}}\sX{x_{n-1},x_n}\] 
induced by the inclusions $N_i\coloneqq \{x_{i-1},x_i\}\hra N$ for $1\leqslant i\leqslant n$ admits an adjoint inverse.\qedhere
	\end{enumerate}
\end{Def}

\begin{Rem}
Condition~\ref{suboneSegal} for $\{x\}=N^1\in\Dpop$ implies that $\sX{x}$ is adjoint equivalent to the empty product which we interpret as the terminal object $\terminal\in\calC$, i.e.\ $\sX{x}$ is contractible. We call such \hypo-simplicial objects (where all the $\sX{x}$ are contractible) \introduce{pointed}. It is then clear that a \hypo-simplicial groupoid is pointed $1$-Segal if and only if it is unpointed $1$-Segal and pointed.
\end{Rem}

\subsubsection{Weak monoid objects and tensor products}\label{sectionmontens}
If we specialize to actual simplicial objects defined on the whole $\Delta$ then we obtain so called weak monoid objects:

\begin{Def}
	A pointed $1$-Segal pseudo-simplicial object $\calA\colon \Dop\to \calC$ is called a \introduce{(weak) monoid object in $\calC$} or, more precisely, a \introduce{(weak) monoidal structure} on the \introduce{underlying object} $\calA_\1$.
\end{Def}
\begin{Rem}
For brevity we will drop the adjective \emph{weak}; if we ever consider non-weak (strict) monoid objects, we will say so explicitly.
\end{Rem}

Let $\calA\colon \Dop\to\calC$ be a monoid object with underlying category $ A\coloneqq\calA_\1$. 
For any $n\geqslant 1$ we can define an $n$-ary tensor product on $ A$ by
\begin{equation}\label{definenarytensorproduct}
	\boxtimes_n\colon  A\times\dots\times A\cong \calA_{\{0,1\}}\times\dots\times\calA_{\{n-1,n\}}\xra{c_\n^{-1}}\calA_\n\lra\calA_{\{0,n\}}\cong  A.
\end{equation}

Note that condition~\ref{suboneSegal} asserts, for $n=0$, that the map $c_\0\colon \calA_\0\xra\approxeq \terminal$ admits an adjoint inverse. Hence the canonical morphism $\1\to\0$ gives a map $I:\terminal\xra{c_\0^\inv} \calA_0\to\calA_1= A$.\\
Next, we consider the diagram
\begin{equation}\label{prerighthalfassociator}
\begin{tikzcd}
	\calA_{\{0,1\}}\times\calA_{\{1,2\}}\times\calA_{\{2,3\}}&	|[alias=X]|\calA_{\{0,1\}}\times\calA_{\{1,2,3\}}\ar[l,"\approxeq"']\rar&	\calA_{\{0,1\}}\times\calA_{\{1,3\}}\\
	&	\calA_{\{0,1,2,3\}}\rar\uar{\approxeq}\ar[dr, ""{name=B, above}]\ar[lu, "\approxeq", ""{above, name=A}] \ar[ru, Leftrightarrow]&	\calA_{\{0,1,3\}}\uar{\approxeq}\dar\ar[Leftrightarrow,from=B]\\
	&&	\calA_{\{0,3\}}
	\ar[from=A, to=X, Leftrightarrow]
\end{tikzcd}
\end{equation}

which comes by applying $\calA$ to various subsets of $\{0,1,2,3\}$ and commutes up to the isomorphisms specified by the pseudo-functoriality of $\calA$. After passing to the adjoint inverses we obtain the following canonical mates (see Appendix~\ref{calculusofmates}):
\begin{equation}\label{righthalfassociatormate}
\begin{tikzcd}
	\calA_{\{0,1\}}\times\calA_{\{1,2\}}\times\calA_{\{2,3\}}&	|[alias=X]|\calA_{\{0,1\}}\times\calA_{\{1,2,3\}}\ar[from=l, "\approxeq"]\rar&	\calA_{\{0,1\}}\times\calA_{\{1,3\}}\\
	&	\calA_{\{0,1,2,3\}}\rar\ar[ur, Leftrightarrow]\ar[from=u, "\approxeq"]\ar[dr, ""{name=B, above}]\ar[from=lu, "\approxeq"', ""{above, name=A}]&	\calA_{\{0,1,3\}}\ar[from=u, "\approxeq"]\dar\ar[Leftrightarrow,from=B]\\
	&&	\calA_{\{0,3\}}
	\ar[from=A, to=X, Leftrightarrow]
\end{tikzcd}
\end{equation}
These natural isomorphisms can be pasted to a natural isomorphism in the big triangle, i.e.\ a natural isomorphism $\alpha_r\colon \boxtimes_3\lLRa \blank\boxtimes(\blank\boxtimes\blank)$.\\
Similarly, we construct a natural isomorphism $\alpha_l\colon \boxtimes_3\lLRa(\blank\boxtimes\blank)\boxtimes\blank$.\\
Composing these two natural isomorphisms gives the so called \introduce{associator}
\[\alpha\colon \blank\boxtimes(\blank\boxtimes\blank)\lLRa(\blank\boxtimes\blank)\boxtimes\blank.\]
\begin{Rem} Expressions of the form $\blank\boxtimes(\blank\boxtimes\blank)$ are just a more intuitive way of writing compositions like $\boxtimes_2\circ (\Id_ A\times\boxtimes_2)$. Similarly we might write $\blank\boxtimes\blank\boxtimes\blank$ for $\boxtimes_3$.
\end{Rem}

Following a similar idea we can pass to adjoint inverses in the pseudo-commutative diagram
\begin{equation*}
	\begin{tikzcd}[column sep=large]
		\terminal\times\calA_{\{0,1\}}&	|[alias=X]|\calA_{\{0\}}\times\calA_{\{0,1\}}\ar[l,"\approxeq"']\ar[r, "{0\mapsfrom 0,0'}", "1\mapsfrom 1"']&	\calA_{\{0,0'\}}\times\calA_{\{0',1\}}\\
		&\calA_{\{0,1\}}\ar[r, "{0\mapsfrom 0,0'}", "1\mapsfrom 1"']\ar[u, "\approxeq"]\ar[dr, ""{name=B, above}]\ar[lu, "\cong", ""{above, name=A}] \ar[ru, Leftrightarrow]&\calA_{\{0,0',1\}}\ar[u, "\approxeq"']\ar[d]\ar[Leftrightarrow,from=B]\\
		&&\calA_{\{0,1\}}
		\ar[from=A, to=X, Leftrightarrow]
	\end{tikzcd}
\end{equation*}
and paste the resulting mates to get a natural isomorphism
\[\lambda\colon {I\boxtimes\blank}\lLRa\Id_ A,\]
called the \introduce{left unitor}.
Similarly we obtain the \introduce{right unitor}
\[\rho\colon {\blank\boxtimes I}\lLRa \Id_ A.\]

\subsubsection{The Mac Lane pentagon}\label{MacLanepent}
Consider the following diagram of functors $ A^4\to  A$ and natural isomorphisms (depicted by single arrows) in which the dashed arrows are yet to be constructed:
	\begin{equation}\label{maclanepentagon}
	\begin{tikzcd}[cramped, column sep=65]
		&(\blank\boxtimes\blank)\boxtimes\blank\boxtimes\blank\ar[ld, "\alpha_l(\boxtimes_2\times\Id\times\Id)"', leftrightarrow]\ar[r, "\alpha_r(\boxtimes_2\times\Id\times\Id)", leftrightarrow]&(\blank\boxtimes\blank)\boxtimes(\blank\boxtimes\blank)\\%
		((\blank\boxtimes\blank)\boxtimes\blank)\boxtimes\blank&&\blank\boxtimes\blank\boxtimes(\blank\boxtimes\blank)\ar[d, leftrightarrow, "\alpha_r(\Id\times\Id\times\boxtimes_2)"]\ar[u, leftrightarrow, "\alpha_l(\Id\times\Id\times\boxtimes_2)"']\\%
		(\blank\boxtimes\blank\boxtimes\blank)\boxtimes\blank\ar[u, leftrightarrow, "\alpha_l\boxtimes\blank"]\ar[d, leftrightarrow, "\alpha_r\boxtimes\blank"']&\boxtimes_4\ar[l, leftrightarrow, dashed, "\beta^3_l"]\ar[dd, leftrightarrow, dashed, "\beta^2_m"]\ar[uu, leftrightarrow, dashed, "\beta^2_l"]\ar[rd, leftrightarrow, dashed, "\beta^3_r"]\ar[ru, leftrightarrow, dashed, "\beta^2_r"]&\blank\boxtimes(\blank\boxtimes(\blank\boxtimes\blank))\\%
		(\blank\boxtimes(\blank\boxtimes\blank))\boxtimes\blank&&\blank\boxtimes(\blank\boxtimes\blank\boxtimes\blank)\ar[u, leftrightarrow, "\blank\boxtimes\alpha_r"']\ar[d, leftrightarrow, "\blank\boxtimes\alpha_l"]\\%
		&\blank\boxtimes(\blank\boxtimes\blank)\boxtimes\blank\ar[lu, leftrightarrow, "\alpha_l(\Id\times\boxtimes_2\times\Id)"]\ar[r, leftrightarrow, "\alpha_r(\Id\times\boxtimes_2\times\Id)"']&\blank\boxtimes((\blank\boxtimes\blank)\boxtimes\blank)\\%
	\end{tikzcd}
	\end{equation}

The dashed arrow $\beta^3_l$ is constructed by passing to adjoint inverses in the pseudo-commutative diagram
\begin{equation*}
\begin{tikzcd}
	\calA_{\{0,1\}}\times\calA_{\{1,2\}}\times\calA_{\{2,3\}}\times\calA_{\{3,4\}}\ar[from=dr, ""{name=A},"\approxeq"']&&\\
	\calA_{\{0,1,2,3,4\}}
	\ar[dr, ""{name=B, very near end}, ""{name=C, below}]\ar[r, "\approxeq"]\ar[u, "\approxeq"]\ar[d]\ar[dr]&	\calA_{\{0,1,2,3\}}\times\calA_{\{3,4\}}\ar[dr, ""{name=A, very near start, left}]&\\
	\calA_{\{0,4\}}
	&\calA_{\{0,3,4\}}\ar[l]\ar[r, "\approxeq"]&	\calA_{\{0,3\}}\times\calA_{\{3,4\}}
\end{tikzcd}
\end{equation*}
and pasting the resulting mates. The natural isomorphism $\beta_r^3$ is constructed in the same way from a similar diagram.\\
The pseudo-commutative diagram for constructing $\beta^2_m$ is
\begin{equation*}
\begin{tikzcd}
	\calA_{\{0,1\}}\times\calA_{\{1,2\}}\times\calA_{\{2,3\}}\times\calA_{\{3,4\}}&	\calA_{\{0,1\}}\times\calA_{\{1,2,3\}}\ar[l,"\approxeq"']\rar\times\calA_{\{3,4\}}&	\calA_{\{0,1\}}\times\calA_{\{1,3\}}\times\calA_{\{3,4\}}\\
	\calA_{\{0,1,2,3,4\}}\ar[ur, "\approxeq"]\rar\dar\ar[u, "\approxeq"]&\calA_{\{0,1,3,4\}}\ar[ur,"\approxeq"]\ar[ld]&\\
	\calA_{\{0,4\}}&&
\end{tikzcd}
\end{equation*}
and $\beta^2_l$ and $\beta^2_r$ arise from a similar diagram.\\
Using the compatibility of canonical mates with pasting it is now easy to show that each of the five pieces of Diagram~\ref{maclanepentagon} commutes.

\begin{Rem}
Diagram~\ref{maclanepentagon} is known as the Mac Lane pentagon and is part of the usual axiomatization~\cite[\S VII.1]{MacLaneCategoriesForTheWorkingMathematician} of monoidal categories. This axiomatization also includes two smaller diagrams involving the left and right unitors; these can be proved by the same techniques.
\end{Rem}
\begin{Rem}
The Mac Lane coherence theorem~\cite[\S VII.2]{MacLaneCategoriesForTheWorkingMathematician} essentially states the converse (at least in the case $\calC=\CAT$): Every monoidal object given as $(A,\boxtimes, I, \alpha,\lambda,\rho)$ satisfying the MacLane pentagon (and a certain coherence statement about $\lambda$ and $\lambda$) can be made into a pointed $1$-Segal (pseudo-)simplicial object. 
More precisely, the assignment described in Section~\ref{sectionmontens} can be upgraded to an equivalence of categories
\[\{\text{$1$-Segal (pseudo)-simpl. $\Dop\to\calC$}\}\xra{\simeq}\{\text{$(A,\boxtimes, I,\alpha,\lambda,\rho)$ in $\calC$ (+ coherence)}\}\qedhere\]
\end{Rem}

\subsubsection{Modules over monoid objects}
Fix a monoid object $\calA\colon \Dop\to\calC$ in $\calC$. We want to define the notion of a module over $\calA$ in the language of pseudo-simplicial objects.

\begin{Def}
	A \introduce{(weak) right $\calA$-module structure} on an object $P\in \calC$ is a pseudo-simplicial object $\calP\colon \Dop\to\calC$ with $\calP_\0=P$ together with a pseudo-natural transformation $\Phi\colon \calP\to\calA$ which satisfies the \introduce{right-relative pointed $1$-Segal} condition:
	\begin{enumerate}[label=\axiomlabelstyle{rpS}, ref=\axiomrefstyle{rpS}]
		\item\label{reloneSegal} For every $n\in\BN$ the map
\[(\iota_0^\star,\Phi_\n)\colon \calP_\n\xra{\approxeq}\calP_\0\times\calA_\n\] 
induced by $\Phi$ and $\iota_0\colon \0=\{0\}\hra\enu n=\n$ admits an adjoint inverse.\qedhere
	\end{enumerate}
\end{Def}

Recall that morphisms $\Phi\colon \calP\to\calA$ of pseudo-simplicial objects $\calA,\calP\colon \Dop\to\calC$ are the same thing as relative pseudo-simplicial objects $\Phi^\mindecor\colon \Dnop\to\calC$. So we can rephrase the definition of modules in terms of relative pseudo-simplicial objects.

\begin{Prop} The morphism $\Phi\colon \calP\to\calA$ is a right module over the monoid object $\calA$ (i.e.\ $\calA$ is pointed $1$-Segal and $\Phi$ satisfies the right relative pointed $1$-Segal condition) if and only if the corresponding right-relative pseudo-simplicial object $\Phi^\mindecor\colon \Dnop\to\calC$ satisfies the pointed $1$-Segal condition~\ref{suboneSegal}.
\end{Prop}
\begin{Prf}
Recall that $\Phi^\mindecor\colon \Dnop\to\C$ is constructed on $N^n=\{x_0<x_1<\dots<x_n\}\in \Dnop$ by $\Phi^\mindecor_N\coloneqq \calP_{\{x_1,\dots,x_n\}}$ if $x_0=\minimal$ and by $\Phi^\mindecor_N\coloneqq \calA_N$ if $x_0\neq \minimal$; the missing maps are induced by the morphism $\Phi$.\\
Hence the case $x_0\neq \minimal$ of condition~\ref{suboneSegal} for $\Phi^\mindecor$ corresponds precisely to the pointed $1$-Segal condition for $\calA$; the case $x_0=\minimal$ corresponds precisely to the relative pointed $1$-Segal condition for $\Phi$.
\end{Prf}

From now on we will use $\Phi\colon \calP\to\calA$ and $\Phi^\mindecor$ interchangeably when talking about right modules over monoid objects.\\

Given a right module $\Phi\colon \calP\to\calA$ over a monoid object $\calA$, we can define an $n$-ary action map as the composition
\[\rcaction_n\colon P\times A^n\cong\Phi^\mindecor_{\{\minimal,0\}}\times\Phi^\mindecor_{\{0,1\}}\times\dots\times\Phi^\mindecor_{\{n-1,n\}}\xra{c_{\{\minimal\}\cup\n}^\inv} \Phi^\mindecor_{\{\minimal,0,\dots,n\}}\xra{}\Phi^\mindecor_{\{\minimal,n\}}\cong P.\]
We can then play the same game as in Section~\ref{sectionmontens} and Section~\ref{MacLanepent} to construct the data of a categorical action with an associator $(\blank \rcaction\blank )\rcaction\blank \lLRa \blank \rcaction(\blank \boxtimes \blank )$ and an unitor $\blank \rcaction I\lLRa \Id$ which satisfy the usual coherence diagrams.
\if{0==0}
Passing to adjoint inverses in the diagrams
\[insert diagram\]
\begin{equation}\label{prerighthalfmoduleassociator}
\begin{tikzcd}
	\calP_{\{0\}}\times\calA_{\{0,1\}}\times\calA_{\{1,2\}}&	|[alias=X]|\calP_{\{0\}}\times\calA_{\{0,1,2\}}\ar[l,"\approxeq"']\rar&	\calP_{\{0\}}\times\calA_{\{0,2\}}\\
	&	\calP_{\{0,1,2\}}\rar\uar{\approxeq}\ar[dr, ""{name=B, above}]\ar[lu, "\approxeq", ""{above, name=A}] \ar[ru, Leftrightarrow]&	\calP_{\{0,2\}}\uar{\approxeq}\dar\ar[Leftrightarrow,from=B]\\
	&&	\calP_{\{2\}}
	\ar[from=A, to=X, equal]
\end{tikzcd}
\end{equation}
we obtain an \introduce{associator} $\blank .(\blank \boxtimes\blank )\lLRa (\blank .\blank ).\blank $ and an \introduce{unitor} $\blank .I\lLRa\Id_M$, where the dot indicates the $1$-ary action $P\times A\to A$ defined above. The compatibility of canonical mates with pasting gives again all the desired coherence diagrams.
\fi
%


\subsubsection{Monoidal categories and monoidal modules}\label{moncatmonmod}
A monoidal category is nothing but a monoid object in the $2$-category of categories. Since we can talk interchangeably about pseudo-functors $\Dop\to\CAT$ (respectively $\Dnop\to\CAT$) and op-fibrations over $\Dop$ (resp.\ $\Dnop$) via the Grothendieck construction (see Appendix~\ref{AppendixGrothendieckConstruction}) we obtain the following definition of a monoidal category and of monoidal modules:
\begin{Def}
\begin{enumerate}
\item A \introduce{monoidal structure} on the underlying category $A$ is an op-fibration of categories $\calA\to\Dop$ (Condition~\ref{Defopfibration}) with $\calA_\1\eqqcolon A$ such that the corresponding pseudo-simplicial category $\Dop\to \CAT$ satisfies the pointed $1$-Segal condition~\ref{suboneSegal}.
\item A \introduce{right monoidal module} is an op-fibration of categories $\Phi^\mindecor\to\Dnop$ such that the corresponding relative pseudo-simplicial category $\Dnop\to\CAT$ satisfies the pointed $1$-Segal condition~\ref{suboneSegal}.\qedhere
\end{enumerate}
\end{Def}

\begin{Rem}In the $2$-category $\CAT$ every equivalence admits an adjoint inverse~\cite[Theorem IV.4.1]{MacLaneCategoriesForTheWorkingMathematician}. Hence in the (pointed or unpointed) $1$-Segal conditions (resp.\ relative $1$-Segal condition) we can just require the chopping maps $c_N$ (resp.\ the maps $(\iota_0,\Phi_\n)$) to be equivalences of categories.
\end{Rem}

\ifx\isstandalone\undefined
\else
\newpage
\bibliographystyle{amsalpha}
\bibliography{literatur}
\end{document}
\fi

\subsection{Lax and oplax morphisms}\label{laxandoplaxmorph}
\ifx\preambleloaded\undefined
   
   \def\isstandaloneDLFOSE{}
\fi

There is also an obvious notion of morphisms (pseudo)-\hypo-simplicial objects.
\begin{Def}
A morphism $\calX\to\calY$ between pseudo-\hypo-simplicial objects ${\calX,\calY\colon \Dpop\to\C}$ in a $2$-category $\calC$ is a pseudo-natural transformation, i.e.\ a collection of morphisms $\beta_N\colon \calX_N\to\calY_N$ in $\calC$ (for $N\in\Dpop$) which is natural up to specified $2$-morphisms
\begin{equation}
	\begin{tikzcd}
	\calX_N\ar[r,"\calX_f"]\ar[d,"\beta_N"]&\calX_M\ar[d,"\beta_M"]\\
	\calY_N\ar[r,"\calY_f"]\ar[ru,Leftrightarrow]&\calY_M
\end{tikzcd}
\end{equation}
(for $f\colon N\la M$ in $\Dpop$) which are compatible with horizontal pasting.
\end{Def}

Often it is necessary to weaken this notion.

\begin{Def}
A \introduce{lax transformation} $\calX\to\calY$ is a collection of maps ${\beta_N\colon \calX\to\calY}$ in $\calC$ together with $2$-morphisms
\begin{equation}\label{laxsquareforbeta}
	\begin{tikzcd}
	\calX_N\ar[r,"\calX_f"]\ar[d,"\beta_N"]&\calX_M\ar[d,"\beta_M"]\\
	\calY_N\ar[r,"\calY_f"]\ar[ru,"\beta_f"{description},Rightarrow]&\calY_M
\end{tikzcd}
\end{equation}
which are compatible with horizontal pasting and such that $\beta_f$ is invertible whenever $f\colon \{x_0,\dots, x_n\}\hla \{x_i,\dots,x_j\}$ is the inclusion of a convex subset.
\end{Def}
\subsubsection{(Lax) monoid homomorphisms}
If $\beta\colon \calA\to\calA'\colon \Dop\to\calC$ is a (lax) morphism of monoid objects (i.e.\ pointed $1$-Segal pseudo-simplicial objects) in $\calC$, then we call $\beta$ a \introduce{(lax) monoid homomorphism}.\\
Part of the structure of a lax monoid homomorphism is a homomorphism
\[b=\beta_\1\colon A=\calA_\1\to\calA'_\1=A'\]
on underlying objects.\\
By taking adjoint inverses in the diagram
\begin{equation}
	\begin{tikzcd}
		\calA_\1\times\dots\times  \calA_\1\ar[d,"\beta_\1\times\dots\times  \beta_\1"'] &\calA_\2\ar[d,"\beta_\2"]\ar[l,"{\approxeq}"']\ar[r]&\calA_\1\ar[d,"\beta_\1"]\\
		\calA'_\1\times \dots\times \calA'_\1 &\calA'_\2\ar[l,"{\approxeq}"']\ar[lu,"\beta_c"description,Leftrightarrow]\ar[ru,"\beta_\nu" description, Rightarrow]\ar[r]&\calA'_\1
	\end{tikzcd}
\end{equation}
(where $\beta_c$ is invertible because the maps $\2\hla\{0,1\},\dots,\{n-1,n\}$ are convex) we can construct  natural transformations
\[\zeta_n\colon b(\blank)\boxtimes' \dots\boxtimes' b(\blank)\lRa b(\blank \boxtimes\dots\boxtimes\blank)\]
which are isomorphisms in the non-lax case.\\
Using the compatibility of $\beta$ and canonical mates with pasting, we can prove that $\zeta$ satisfies the obvious coherence axioms which we will not spell out here.\\

By taking the case $\Deltap\coloneqq \Delta_\mindecor$ we can similarly define (lax) module homomorphisms $\beta\colon \calP\to\calP'$ which come equipped with transformations
\[\zeta_n\colon b(\blank)\rcaction b(\blank)\boxtimes' \dots\boxtimes' b(\blank)\lRa b(\blank \rcaction\blank \boxtimes\dots\boxtimes\blank).\]

For the situation of $\calC=\CAT$ (i.e.\ monoidal functors and categorical modules) we can again translate using the Grothendieck construction.

\begin{Def}
A (lax) monoidal functor $\calA\to\calA'$ between monoidal categories (seen as op-fibrations $p,p'\colon \calA,\calA'\to\Dop$) is a functor $\beta\colon \calA\to\calA'$ which commutes with the fiber functors (i.e. $p'\circ \beta =p$) together with transformations~\ref{laxsquareforbeta} which make the induced collection $\beta_n\colon \calA_n\to\calA'_n$ into a (lax) monoid homomorphism.\\
The definition of a (lax) module functor between categorical modules is the same with $\Dop$ replaced by $\Dnop$.
\end{Def}
Of course we can copy everything in Section~\ref{laxandoplaxmorph} verbatim and just change the direction of the transformation inhabiting the square~\ref{laxsquareforbeta}. We then get the dual notion of \introduce{oplax homomorphisms}, \introduce{oplax monoidal functors} and \introduce{oplax module functors}.

\ifx\isstandaloneDLFOSE\undefined
\else
\newpage
\bibliographystyle{amsalpha}
\bibliography{literatur}
\end{document}
\fi

\subsection{The $2$-Segal condition}\label{sectionTwoSegalCondition}
\ifx\preambleloaded\undefined
   
   \def\isstandalone{}
\fi

Recall the abbreviation $M_i\coloneqq \{x_{i-1},x_i\}$ for an object $\{x_0,\dots,x_m\}=M^m\in\Delta$ and the notation $\ol f(M')\coloneqq \{x\in N\mid f(\min M')\leqslant x\leqslant f(\max M')\}$ for a morphism $f\colon M\to N$ in $\Delta$ and a subset $M'\in M$.

\begin{Def}\label{Definitionunital}
Let $\Deltap\subseteq \Delta$ be an admissible subcategory.
A \hypo-simplicial groupoid $\calX\colon \Dpop\to\Grpd$ is called \introduce{$2$-Segal} if
\begin{enumerate}[label=\axiomlabel S2, ref=\axiomref S2]
	\item\label{generalunitalsubsimplicial} for every map $f\colon N\la M$ in $\Dpop$ the following square is a $2$-pullback:
	\begin{equation}\label{subsimplicialunitalcd}\begin{tikzcd}
		\calX_{\ol f(M)}\dar{\calX_f}\rar& \calX_{\ol f(M_1)}\times\dots \times \calX_{\ol f(M_m)}\dar{\calX_{f\restriction M_1}\times\cdots\times \calX_{f\restriction M_m}}\\
		\calX_M\rar & \calX_{M_1}\times\dots\times\calX_{M_m}
	\end{tikzcd}\qedhere\end{equation}
\end{enumerate}
\end{Def}

Of course this definition does not just make sense for \hypo-simplicial groupoids but also for \hypo-simplicial objects in any category $\calC$ which has some notion of homotopy limits, i.e.\ homotopy fiber products.

\subsubsection{The $2$-Segal condition for simplicial groupoids}
In the case where $\Deltap=\Delta$ we can give an easier description of the $2$-Segal condition which agrees with the one given by Dyckerhoff and Kapranov~\cite[Proposition 2.3.2, Definition 2.5.3]{DyckerhoffKapranovHigherSegalSpacesI}

\begin{Prop}\label{definition2segal}
A simplicial groupoid $\calX\colon\Dop\ra \Grpd$ is $2$-Segal if and only if the following two commutative diagrams of groupoids are 2-pullback squares for all $0\leqslant i<j\leqslant  n$ and all $0\leqslant k\leqslant n$ respectively:
\begin{equation}\label{cd2segalgroupoid}
\begin{tikzcd}
\calX_{\{0,\dots,n\}}\rar\dar &\calX_{\{i,\dots,j\}}\dar\\
\calX_{\{0,\dots,i,j,\dots,n\}}\rar&\calX_{\{i,j\}}
\end{tikzcd}
\end{equation}\begin{equation}\label{cdunitalgroupoid}
\begin{tikzcd}
\calX_{\{0,\dots,n\}}\rar\dar{\sigma_k} &\calX_{\{k\}}\dar{\sigma_k}\\
\calX_{\{0,\dots,k,k',\dots n\}}\rar&\calX_{\{k,k'\}}
\end{tikzcd}
\end{equation}
where $\sigma_k$ is the degeneracy map induced by $k,k'\mapsto k$.
\end{Prop}

\begin{dCor}\label{oneSegalgivesunital}\cite[Proposition 2.3.3, Proposition 2.5.3]{DyckerhoffKapranovHigherSegalSpacesI}
If a simplicial groupoid $\calX\colon \Dpop\to\Grpd$ satisfies the unpointed $1$-Segal condition then $\calX$ is $2$-Segal.
\end{dCor}

For the proof of Proposition~\ref{definition2segal} we will need an easy Lemma about $2$-pullbacks of groupoids.
\begin{Lem}\label{lemmaaddisotopullback}
Let \[\begin{tikzcd}
A\ar[d,"f"']\rar&B\ar[d,"f'"]&&A\ar[d,"f"']\rar&Y\dar{\cong}\\
C\ar[ur,Leftrightarrow]\rar&D&&C\ar[ur,Leftrightarrow]\rar&Y'
\end{tikzcd}\]
be two pseudo-commuting squares. Then the induced square
\[\begin{tikzcd}
A\ar[d,"f"']\rar&B\dar{{f'}\times{\cong}}\times Y\\
C\ar[ur,Leftrightarrow]\rar&D\times Y'
\end{tikzcd}\]
is a 2-pullback square if and only if the first square was $2$-pullback.
\end{Lem}
\begin{Prf}
For strict pullbacks the result is clear. Hence it is also true if $f'$ (hence $f'\times{\cong}$) is an iso-fibration because in this case we can compute $2$-pullbacks by strict pullbacks~\cite[Proposition 2.5]{DyckerhoffHigherCategoricalAspectsOfHallAlgebras}.\\
We can replace $B$ by the equivalent groupoid $\ol B$ which has additional copies $b_\alpha$ of every object $b\in B$ labeled by the isomorphisms $\alpha\colon f'(b)\xra{\cong} d$ in $\calD$. We can also replace $f'$ by the iso-fibration $\ol{f'}\colon \ol B\to D$ which maps the \roughly{identity} $b\to b_\alpha$ to the morphism $\alpha$. Hence we can deduce the general case of the Lemma from the case where $f'$ is an iso-fibration.
\end{Prf}

\begin{Prf}[\Proofof{Proposition~\ref{definition2segal}}]
Let $\calX\colon \Dop\to\Grpd$ be a simplicial groupoid satisfying the two $2$-pullback conditions~\ref{cd2segalgroupoid} and \ref{cdunitalgroupoid}.
We call a map $f\colon \n \la\m$ \emph{hc-good}, if the corresponding square \ref{subsimplicialunitalcd} is a 2-pullback; then we proceed in three steps:
\begin{enumerate}
\item\label{injishcgood} Every injective map is hc-good.
\item Every degeneracy map $\{0,\dots,k,\dots,m\}\thla\{0,\dots,k,k',\dots,m\}$ is hc-good (for $0\leqslant k\leqslant m$).
\item\label{compofhcgood} If $\n\xla{g}\m$ and $\m\xla{f}\l$ are hc-good and $f$ is surjective then the composition $gf$ is hc-good.
\end{enumerate}
Statements \ref{injishcgood}-\ref{compofhcgood} imply that every map $\bullet\la\bullet$ in $\Dop$ is hc-good, since it can be factored as a composition $\bullet\la\bullet\thla\cdots\thla\bullet$ where the leftmost map is injective and the others are degeneracies (which are surjective).
\begin{enumerate}

\item Let $f$ be injective. We can assume that $f$ is the inclusion $\enu n \hla\{f(0),\dots, f(m)\}$ of an \m-indexed subset.\\
Consider the following squares for $i=1,\dots, m$.
\begin{equation}\label{cdpartialsubdivision}\begin{tikzcd}
\sX{f(0),f(1),\dots,f(i-1),f(i-1)+1,\dots,f(m)-1, f(m)}\rar\dar&	\prodl_{k=1}^{i-1}\sX{f(k-1),f(k)}\times\prodl_{k=i}^{m}\sX{f(k-1),\dots,f(k)}\dar\\
\sX{f(0),f(1),\dots,f(i-1),f(i),f(i)+1,\dots, f(m)}\rar&	\prodl_{k=1}^{i}\sX{f(k-1),f(k)}\times\prodl_{k=i+1}^{m}\sX{f(k-1),\dots,f(k)}
\end{tikzcd}\end{equation}
By removing the factors of the rightmost vertical map which appear on both sides (using Lemma~\ref{lemmaaddisotopullback}) we can reduce this square to
\begin{equation}\begin{tikzcd}
\sX{f(0),f(1),\dots,f(i-1),f(i-1)+1,\dots,f(i),\dots,f(m)-1, f(m)}\rar\dar&	\sX{f(i-1),\dots, f(i)}\dar\\
\sX{f(0),f(1),\dots,f(i-1),f(i),f(i)+1,\dots, f(m)}\rar&	\sX{f(i-1),f(i)}
\end{tikzcd}\end{equation}
which is the instance of the square \ref{cd2segalgroupoid} corresponding to the line $\{f(i-1),f(i)\}$ in the polygon
\[\{f(0),f(1),\dots,f(i-1),f(i-1)+1,\dots, f(i), f(i)+1,\dots,f(m)-1, f(m)\}.\] Hence we know that the square \ref{cdpartialsubdivision} is a 2-pullback.\\
Finally, we can vertically paste the squares~\ref{cdpartialsubdivision} for the various $i=1,\dots,m$ to obtain the desired 2-pullback square~\ref{subsimplicialunitalcd}.

\item If $\sigma_k\colon \enu{m}\thla\{0,\dots,k,k',\dots m\}$ is the $k$-th degeneracy for $0\leqslant k\leqslant m$ then the Diagram~\ref{subsimplicialunitalcd} will look as follows:
\[\begin{tikzcd}
\calX_\enu {m}\dar{\calX_{\sigma_k}}\rar&\sX{0,1}\times\dots\times\sX{k-1,k}\times\sX{k=k'}\times\sX{k',k+1}\times\sX{m-1,m}\dar{\Id\times\dots\times \calX_{\sigma_k}\times \dots\times \Id}\\
\sX{0,\dots,k,k',\dots,m}\rar&\sX{0,1}\times\dots\times\sX{k-1,k}\times\sX{k,k'}\times\sX{k',k+1}\times\dots\times\sX{m-1,m}
\end{tikzcd}\]
This square arises from Diagram~\ref{cdunitalgroupoid} (which is a 2-pullback square because $\calX$ is $2$-Segal) by adding a finite number of isomorphisms to the right map; hence it is still a 2-pullback by Lemma \ref{lemmaaddisotopullback}
\item Diagram~\ref{subsimplicialunitalcd} for $gf$ can be subdivided as follows:
\[\begin{tikzcd}
\sX{gf(0), gf(0)+1,\dots, gf(l)}\rar\dar{\calX_g}&	\sX{gf(0),\dots,gf(1)}\times\dots\times\sX{gf(l-1),\dots,gf(l)}\dar{\calX_{g\restriction\{f(0),\dots,f(1)\}}\times\cdots}\\
\sX{f(0), f(0)+1,\dots, f(l)}\dar{\calX_f}\rar&	\sX{f(0),\dots,f(1)}\times\dots\times\sX{f(l-1),\dots,f(l)}\dar{\calX_{f\restriction\{0,1\}}\times\cdots}\\
\sX{0,\dots,l}\rar&	\sX{0,1}\times\dots\times\sX{l-1,l}
\end{tikzcd}\]
The bottom square is a 2-pullback since $f$ is hc-good. Moreover, since $f$ is surjective, the upper square reduces to
\[\begin{tikzcd}
\sX{gf(0), gf(0)+1,\dots, gf(l)}\rar\dar{\calX_g}&	\sX{gf(0),\dots,gf(1)}\times\dots\times\sX{gf(l-1),\dots,gf(l)}\dar{\calX_{g\restriction\{f(0),f(1)\}}\times\cdots}\\
\sX{0,\dots,m}\rar&	\sX{0=f(0),f(1)}\times\dots\times\sX{f(l-1),f(l)=m}
\end{tikzcd}\]
Note, that if $f(i-1)=f(i)$ for some $1\leqslant i\leqslant l$ then the corresponding map 
\[\calX_{g\restriction\{f(i-1),f(i)\}}\colon \sX{f(i)}\to\sX{gf(i)}\]
is an isomorphism; hence we can without loss remove those factors from the diagram using Lemma \ref{lemmaaddisotopullback}. What is left after removing those factors will be precisely the instance of the square \ref{subsimplicialunitalcd} for $g$, which is a 2-pullback since $g$ is hc-good.
\end{enumerate}

Conversely, if $\calX$ is $2$-Segal then we can specialize $f\colon N\la M$ to $\enu{n}\hla\{0,\dots,i,j,\dots,n\}$ and to $\enu{n}\thla\{0,\dots,k,k',\dots,n\}$ in the $2$-pullback square~\ref{subsimplicialunitalcd}. This gives us precisely (after dropping some vertical identities using Lemma~\ref{lemmaaddisotopullback}) that \ref{cd2segalgroupoid} and \ref{cdunitalgroupoid} are $2$-pullback squares.
\end{Prf}

The concept of $2$-Segal simplicial groupoids wouldn't be any use to us if it did not apply to our main object of interest $\calX^\calA$. So here comes the obligatory example:
\begin{Prop}\cite[Theorem 2.10]{DyckerhoffHigherCategoricalAspectsOfHallAlgebras}
	If $\calA$ is a proto-abelian  category then the simplicial groupoid $\calX^\calA$ of flags in $\calA$ is $2$-Segal.
\end{Prop}
\begin{Prf} We will later prove a very similar statement in the case of bounded flags (see Proposition~\ref{boundedflags2Segal}), hence we omit this proof.
\end{Prf}

\subsubsection{The relative $2$-Segal condition}

We want to define the notion of a $2$-Segal morphism of simplicial groupoids.
\begin{Def}\label{morphism2Segal}
A morphism $\Phi\colon \calY\to\calX$ of simplicial objects is called $2$-Segal if the right-relative simplicial groupoid $\Phi^\mindecor$ is $2$-Segal.
\end{Def}
\begin{Rem}
Definition~\ref{morphism2Segal} seems arbitrarily biased towards the right-relative point of view. It will turn out, however, that this is only an illusion; the left-relative point of view would have given the same definition (see Remark~\ref{leftorrightmattersrelativelylittle}).
\end{Rem}

We can mimic the proof of Proposition~\ref{definition2segal} to show an analogous statement for a large class of \hypo-simplicial groupoids
\begin{Prop}\label{sub2segalcrit}
Let $\Deltap\subseteq\Delta$ be an admissible subcategory with \roughly{enough degeneracies}, i.e.\ assume that every map in $\Deltap$ can be factorized (in $\Deltap$) as a composition $\bullet\hla\bullet\thla\dots\thla\bullet$ of some degeneracies followed by an inclusion.\\
Then a \hypo-simplicial groupoid $\calX\colon\Dnop\ra \Grpd$ is $2$-Segal if and only if for $\{x_0<\dots<x_n\}=N^n\in\Dnop$ the following two commutative diagrams of groupoids are 2-pullback squares (for all $0\leqslant i<j\leqslant  n$ and all $0\leqslant k\leqslant n$ respectively)
\begin{equation}\label{sub2segalcriterium}
\begin{tikzcd}
\calX_{\{x_0,\dots,x_n\}}\rar\dar &\calX_{\{x_i,\dots,x_j\}}\dar\\
\calX_{\{x_0,\dots,x_i,x_j,\dots,x_n\}}\rar&\calX_{\{x_i,x_j\}}
\end{tikzcd}
\end{equation}
\begin{equation}\label{sub2unitalcriterium}
\begin{tikzcd}
\calX_{\{x_0,\dots,x_n\}}\rar\dar{\sigma_k} &\calX_{\{x_k\}}\dar{\sigma_k}\\
\calX_{\{x_0,\dots,x_k,x_k',\dots n\}}\rar&\calX_{\{x_k,x_k'\}}
\end{tikzcd}
\end{equation}
whenever the degeneracy map $\sigma_k$ is defined (i.e. whenever the map $x_k,x_k'\mapsto x_k$ lies in $\Deltap$).
\end{Prop}

\begin{Rem}
Clearly $\Delta_\mindecor$ has enough degeneracies; the only degeneracies which are not defined are those corresponding to $\minimal,x\mapsto x$.
\end{Rem}

We can further refine the criterium of Proposition~\ref{sub2segalcrit} in the case where $\Deltap=\Delta_\mindecor$:

\begin{Prop}\label{criterionrelativeunital}
	A relative simplicial groupoid $\calX\colon \Dnop\to\Grpd$ (which we can think of as a morphism $\Phi\colon \calX^\in\to\calX^{\not\in}$ of groupoids) is $2$-Segal if and only if all of the following conditions are satisfied:
	\begin{enumerate}
		\item\label{indprooftargetunital} The target simplicial groupoid $\calX^{\not\in}\colon \Dop\simeq\Delta_{\not\in}^\op\subset\Dnop\to\Grpd$ is $2$-Segal.
		\item\label{indproofsourcesegal} The source simplicial groupoid $\calX^{\in}\colon \Dop\simeq \Delta_\in^\op\subset\Dnop\to\Grpd$ satisfies the unpointed $1$-Segal condition~\ref{unpointsuboneSegal}.
		\item The morphism $\Phi\colon \calX^\in\to\calX^{\not\in}$ satisfies the following \introduce{relative $2$-Segal conditions}:
		\begin{enumerate}[label=\axiomlabel{rS}{2\roman*}, ref=\axiomref{rS}{2\roman*}]
		\item\label{indproofrel2segal} For every $0\leqslant i<j\leqslant n$ the following square is $2$-pullback:
		\begin{equation}\label{cdrel2segalgroupoid}
			\begin{tikzcd}
				\calX^\in_{\{0,\dots,n\}}\rar\dar &\calX^{\not\in}_{\{i,\dots,j\}}\dar\\
				\calX^\in_{\{0,\dots,i,j,\dots,n\}}\rar&\calX^{\not\in}_{\{i,j\}}
			\end{tikzcd}
		\end{equation}

		\item\label{indproofrelunital} For every $0\leqslant k\leqslant n$ the following square is $2$-pullback:
		\begin{equation}\label{cdrelunitalgroupoid}
			\begin{tikzcd}
				\calX^\in_{\{0,\dots,n\}}\rar\dar{\sigma_k} &\calX^{\not\in}_{\{k\}}\dar{\sigma_k}\\
				\calX^\in_{\{0,\dots,k,k',\dots n\}}\rar&\calX^{\not\in}_{\{k,k'\}}
			\end{tikzcd}\qedhere
		\end{equation}
		\end{enumerate}
	\end{enumerate}
\end{Prop}
\begin{Prf}
We use Proposition~\ref{sub2segalcrit} and identify the various cases with the conditions of Proposition~\ref{criterionrelativeunital}
\begin{enumerate}
\item corresponds to the case $\minimal\not\in N$
\item corresponds to Diagram~\ref{sub2segalcriterium} in the case $\minimal\in N$, $i=0$.
\item[\ref{indproofrel2segal}] corresponds to Diagram~\ref{sub2segalcriterium} in the case $\minimal\in N$, $i\neq 0$. 
\item[\ref{indproofrelunital}] corresponds to Diagram~\ref{sub2unitalcriterium} in the case $\minimal\in N$ and $k\neq 0$.
\end{enumerate}
Note that the case $\minimal\in N$ and $k=0$ cannot appear since that is the situation where the degeneracy $\sigma_k$ is not defined.
\end{Prf}
\begin{Rem}
A definition of relative $2$-Segal objects was independently proposed by Young~\cite{YoungRelative2SegalSpaces}. Proposition~\ref{criterionrelativeunital} shows that our notion of relative $2$-Segal simplicial objects is the same as his notion of \emph{unital} relative $2$-Segal simplicial objects. We obtain his definition of (not necessarily unital) relative $2$-Segal simplicial objects by dropping condition \ref{indproofrelunital}.
\end{Rem}

\begin{Rem}\label{leftorrightmattersrelativelylittle}
Observe that the conditions of Proposition~\ref{criterionrelativeunital} do not change if we view the morphism $\Phi\colon \calX^\in\to\calX^{\not\in}$ as a \emph{left}-relative simplicial groupoid $\calX\colon \Dmop\to\Grpd$ and require it to be $2$-Segal. Hence we can speak of a $2$-Segal morphism $\calY\to\calX$ of groupoids without having to specify left or right.
\end{Rem}

\ifx\isstandalone\undefined
\else
\newpage
\bibliographystyle{amsalpha}
\bibliography{literatur}
\end{document}
\fi


\newpage\section{The generalized Hall construction}

\subsection{Monoidal left derivators of groupoids}\label{sectionmonoidalderivators}
\ifx\preambleloaded\undefined
   
   \def\isstandalone{}
\fi

In this section we want to introduce what we need as a second ingredient to construct the Hall monoidal category. We want to give an axiomatic characterization of what the assignment $\calV\colon \Grpd\to\CAT$ should satisfy: It should be a componentwise cocontinuous monoidal left derivator.\\
We will introduce all these terms in a rather ad hoc manner without delving into the abstract theory of derivators. We will make some simplifying assumptions which are probably not necessary because in the end we are mostly concerned with the case $\calV\coloneqq [\blank ,\Vect\BC]$ and the additional complications wouldn't be worth it.\\
See Groth's thesis~\cite{GrothOnTheTheoryOfDerivators} for a systematic (and much more careful) development of the theory of derivators and monoidal derivators.

\subsubsection{Ad hoc definitions}
Consider a $2$-functor $\calV\colon \Grpd\to\CAT$ which is contravariant on $1$-cells (but covariant on $2$-cells). Such a functor is called a \introduce{pre-derivator (of groupoids)}. To stress the contravariance we denote the induced map on $1$-cells by $\pV{}$.
\begin{Def} We call $\calV$ a \introduce{left derivator (of groupoids)} if the following conditions hold:
\begin{enumerate}[label=\axiomlabelstyle{D}, ref=\axiomrefstyle{D}]
	\item\label{derivatortakescoprodtoprod} The $2$-functor $\calV$ takes coproducts to products.
	\item\label{derivatorhasleftkanextensions} For each map of groupoids $f\colon A\to B$ the induced map $\pV f\colon \calV(B)\to\calV(A)$ admits a left adjoint $\sV f\colon \calV(A)\to\calV(B)$
	\item\label{2pullbacksarehc} Assume we have a $2$-pullback square of groupoids:
	\begin{equation}
		\begin{tikzcd}
			X\rar{f'}\dar{g'}& B\dar{g}\\
			A\ar[ru, Leftrightarrow]\rar{f}& D
		\end{tikzcd}
	\end{equation}
	Then this square is required to be \introduce{$\calV$-exact}. This means that the canonical mates (see Appendix~\ref{calculusofmates})
	\begin{gather}
	\sV{f'}\circ \pV{g'} \Longrightarrow \pV{g}\circ \sV{f}\label{canonicalmateforhcbis}\\
	\sV{g'}\circ \pV{f'}\Longrightarrow \pV{f}\circ \sV{g}\label{canonicalmateforhc}
	\end{gather}
	of the induced diagram
	\begin{equation}
		\begin{tikzcd}
			\calV(X)\ar[from=r,"\pV{f'}"']\ar[from=d,"\pV{g'}"]& \calV(B)\ar[from=d,"{\pV g}"']\\
			\calV(A)\ar[ru, Leftrightarrow]\ar[from=r,"\pV{f}"]& \calV(D)
		\end{tikzcd}
	\end{equation}
	are required to be isomorphisms.\qedhere
\end{enumerate}
\end{Def}

\begin{Rem}
	In a derivator~\cite[Definition 1.11]{GrothOnTheTheoryOfDerivators} the maps $\pV f\colon \calV(B)\to\calV(A)$ are required to have both a left and a right adjoint. In case we only need left adjoints, we talk about \emph{left} derivators.\\
	Usually (pre-)derivators are defined as $2$-functors $\Cat\to\CAT$ defined on all small categories and not just on groupoids. Moreover the usual axioms are a bit different (and stronger); what we call Axiom~\ref{2pullbacksarehc} can then be proved as a proposition. For us Condition~\ref{2pullbacksarehc} is really the key statement; hence it makes sense to take it as an axiom.\\
	We will not always mention that our left derivators are only defined on groupoids and by \qquote{left derivator} we will always mean \qquote{left derivator of groupoids}.
\end{Rem}

\begin{PropDef}
Let $\calV$ be a cocomplete category. Then $\calV$ defines a left derivator
\[\calV= [\blank ,\calV]\colon \Grpd\to\CAT\]
called the \introduce{left derivator represented by $\calV$}.
\end{PropDef}
\begin{Prf}\begin{itemize}
\item [\ref{derivatortakescoprodtoprod}] Clear.
\item [\ref{derivatorhasleftkanextensions}] For a map $f\colon A\to B$ of groupoids, the induced functor $\pV f\colon [B,\calV]\to[A,\calV]$ has a left adjoint $\sV f\colon [A,\calV]\to[B,\calV]$ given by left Kan extension~\cite[\S X.]{MacLaneCategoriesForTheWorkingMathematician} along $f$, which exists because $\calV$ is cocomplete.
\item[\ref{2pullbacksarehc}] We have to show that the mate~\ref{canonicalmateforhc} is an isomorphism evaluated at any object $a\in A$; the proof for the mate~\ref{canonicalmateforhcbis} is similar.\\
Fix $a\in A$ and consider the following diagram where the two small squares (hence the big rectangle) are $2$-pullbacks:
	\begin{equation}
		\begin{tikzcd}
			P\ar[d,"\pi"]\rar{p}&X\rar{f'}\dar{g'}& B\dar{g}\\
			\terminal\ar[ru,Leftrightarrow]\ar[r,hookrightarrow,"a"]&A\ar[ru, Leftrightarrow]\rar{f}& D
		\end{tikzcd}
	\end{equation}
Applying $\calV=[\blank ,\calV]$ and passing to left adjoint gives the following pasting of canonical mates
	\begin{equation}\label{pastingdiagramtocheckmateof2pullback}
		\begin{tikzcd}[column sep=50]
			{[P,\calV]}\ar[from=rr,"\pV{f'p}"',bend right=20]\ar[dr,Rightarrow]\ar[d,"\colim_P=\sV\pi"']\ar[from=r,"\pV p"']&{[X,\calV]}\ar[dr,Rightarrow]\ar[from=r,"\pV{f'}"']\ar[d,"\sV{g'}"'near start]& {[B,\calV]}\ar[d,"{\sV g}"]\\
			\calV\ar[from=rr,bend left,"\ev_{f(a)}=\pV{f(a)}"]\ar[from=r,"\ev_a=\pV{a}"]&{[A,\calV]}\ar[from=r,"\pV{f}"]& {[D,\calV]}
		\end{tikzcd}
	\end{equation}
Observe that the $2$-fiber $P\simeq\twofib [g']Xa\simeq  \twofib [g]B{f(a)}$ of groupoids coincides with the comma categories $\overcat X a$ and $\overcat B {f(a)}$. Hence the pointwise formula for Kan extensions~\cite[Theorem X.3.1]{MacLaneCategoriesForTheWorkingMathematician} states precisely that the canonical mates
\begin{align*}
\colim_P\circ \pV p &\lRa\ev_a\circ \sV{g'}\\
\colim_P\circ \pV {f'p} &\lRa \ev_{f(a)}\circ \sV{g}
\end{align*}
living in the left square and in the big rectangle are isomorphisms.\\
It follows that the canonical mate~\ref{canonicalmateforhc} living in the right square of Diagram~\ref{pastingdiagramtocheckmateof2pullback} is an isomorphism after composition with $\ev_a$ (i.e.\ after evaluation at $a$). Because the obect $a \in A$ was arbitrary, we conclude that all the components of the mate transformation are isomorphisms, hence the mate itself is an isomorphism.\qedhere
\end{itemize}\end{Prf}

\begin{Def}A \introduce{pseudo-natural transformation} between two pre-derivators $\calV,\calW\colon \Grpd\to \CAT$ is a collection of functors $\calV(A)\xra{F_A} \calW(A)$ for $A\in \Grpd$ together with specified naturality isomorphisms living in the squares
\begin{equation}\label{pseudonatural}
	\begin{tikzcd}
		\calV(A)\ar[r, "F_A"]&\calW(A)\ar[dl, Leftrightarrow]\\
		\calV(B)\ar[u,"\pV f"]\ar[r, "F_B"']&\calW(B)\ar[u, "\calW^\star(f)"']
	\end{tikzcd}
\end{equation}
such that some coherence conditions are satisfied.
\end{Def}

\begin{Def}
A \introduce{monoidal left derivator} is a left derivator $\calV\colon \Grpd\to\CAT$ together with pseudo-natural transformations $\otV{}\calV\times\calV\to \calV$ and $\terminal\to\calV$ satisfying the usual associativity and unitality conditions..
\end{Def}
\begin{Rem}\label{cheatingoncoherence}
Of course this is not the correct definition. What we should really be saying is that the associativity and unitality conditions are satisfied \emph{up to coherent $2$-isomorphisms}. In other words we would have to introduce the $2$-category of left derivators and define a monoidal left derivator as a weak monoid object in this $2$-category~\cite{GrothMonoidalDerivatorsAndAdditiveDerivators}. The technical details are quite burdensome, so we shall make do with this naive definition. In fact we shall do even worse: we will assume that our pseudo-natural transformations are strict, i.e.\ that all the instances of Square~\ref{pseudonatural} commute on the nose; in other words we talk about \emph{natural} and not pseudo-natural transformations.\\
We can then define $n$-ary tensor products $\otV{}\colon \calV\times \dots\times \calV\to \calV$ (including the unit $J\colon \terminal\to\V$) without having to track all the involved isomorphisms.
\end{Rem}

\begin{Expl} If $\otimes\colon \calV\times\calV\to \calV$ is a cocomplete (strict) monoidal category, then the corresponding represented left derivator is monoidal with tensor product given by the composition
\[\otV{}\colon [\blank ,\calV]\times[\blank ,\calV]\xla{\cong}[\blank ,\calV\times \calV]\xra{[\blank ,\otimes]}[\blank ,\calV].\qedhere\]
\end{Expl}
\begin{Rem}
Following the same spirit as Remark~\ref{cheatingoncoherence} we will pretend that all our monoidal categories are strict so that we don't have to track the coherence isomorphisms. So for instance in $\Vect\BC$ we want to pretend that $(U\otimes V)\otimes W$ is \emph{equal} to $U\otimes (V\otimes W)$.
\end{Rem}

\subsubsection{Multi-valued tensor products}\label{multivaluedtensorprod}
Fix a monoidal left derivator $\calV\colon \Grpd\to\CAT$.
\begin{Cstr}Let $\n\xla{f}\m\xla{g}\l$ be a diagram in $\Dop$ and let $A_i$, $B_j$ and $C_k$ be groupoids for $1\leqslant i\leqslant n$, $1\leqslant j\leqslant m$ and $1\leqslant k\leqslant l$ together with maps $y_{jk}\colon C_k\to B_j$ for  $j=g(k-1)+1,\dots, g(k)$ and $x_{ij}\colon B_j\to A_i$ for $i=f(j-1)+1,\dots, f(j)$. It is sometimes convenient to collect all the $x_{ij}$'s for a fixed $j$ into the map $x_j\colon B_j\to \prod_{i=f(j-1)+1}^{f(j)}A_i$\\
We define the \introduce{multi-valued tensor product}
\begin{equation*}
	\fotV[x_{ij}]f\colon \prod_{i=1}^n\calV(A_i)\xra{\pr}\prod_{j=1}^m\prod_{i=f(j-1)+1}^{f(j)}\calV(A_i)\xra{\prod_j\prod_i\pV{x_{ij}}}\prod_{j}^m\prod_i^{f(j)}\calV(B_j)\xra{\prod_j\otV{B_j}}\prod_{j=1}^m\calV(B_j)\qedhere
\end{equation*}
\end{Cstr}
We call the structure maps $y_{jk}$ and $x_{ij}$ \introduce{composable} if all the expressions $z_{ik}\colon C_k\xra{y_{jk}}B_j\xra{x_{ij}}A_i$ are independent from the choice of $j$, as the notation suggests.
In this case we have $\fotV[y_{jk}]g\circ\fotV[x_{ij}]f=\fotV[z_{ik}]{fg}$ as can be seen by considering the following commutative diagram:
\begin{equation*}
\begin{tikzcd}[column sep=65]
		\prod_i^n\calV(A_i)\ar[dd,"\pr"]\ar[r, "\pr"]&\prod_j^m\prod_i^{f(j)}\calV(A_i)\ar[d, "\prod_j^m\prod_i^{f(j)}\pV{x_{ij}}"]&\\
	&\prod_j^m\prod_i^{f(j)}\calV(B_j)\ar[d, "\pr"]\ar[r, "\prod_j^m\otV{B_j}"]&\prod_j^m\calV(B_j)\ar[d,"pr"]\\
	\prod_k^l\prod_j^{g(k)}\prod_i^{f(j)}\calV(A_i)\ar[rd, "\prod_k\prod_j\prod_i\pV{z_{ik}}"']&\prod_k^l\prod_j^{g(k)}\prod_i^{f(j)}\calV(B_j)\ar[d, "\prod_k\prod_j\prod_i\pV{y_{jk}}"]\ar[r,"\prod_k\prod_j\otV{B_j}"]&\prod_k^l\prod_j^{g(k)}\calV(B_j)\ar[d,"\prod_k\prod_j\pV{y_{jk}}"]\\
	&\prod_k^l\prod_j^{g(k)}\prod_i^{f(j)}\calV(C_k)\ar[r,"\prod_k\prod_j\otV{C_k}"]\ar[rd, "\prod_k\otV{C_k}"']&\prod_k^l\prod_j^{g(k)}\calV(C_k\ar[d,"\prod_k\otV{C_k}"])\\
	&&\prod_k^l\calV(C_k)
\end{tikzcd}
\end{equation*}
Here the second square on the right commutes because $\otV{}$ is a natural transformation, the triangle in the bottom right is associativity of $\otV{}$, and the other two pieces commute trivially (using $z_{ik}=x_{ij}y_{jk}$).

\begin{Rem}
In the case where $n=0$ and $x_j\colon B_j\to\terminal$ are the unique maps to the point, the multi-valued tensor product
\[\fotV{f}\colon \terminal=\prod_{j=1}^m \terminal(B_j) \to \prod_{j=1}^m \calV(B_j)\]
(induced by $f\colon \0\la\m$) is just a product of unit maps $J_{B_j}\colon \terminal(B_j)\to\calV(B_j)$ coming from the monoidal structure on $\calV$. Here $\terminal(\blank)$ denotes the terminal derivator (a.k.a.\ the terminator).\\
The multi-valued tensor product induced by $f\colon \{0,\dots,k,\dots,n\}\thla\{0,\dots,k,k',\dots,n\}$ inserts a unit at position $k$.
\end{Rem}

\subsubsection{Naturality of multi-valued tensor products}
Next consider maps $B_j\xra{\beta_j}B_j'$ of groupoids (for $1\leqslant j\leqslant m$) and a bunch of commutative diagrams in $\Grpd$ as follows (for $f\colon \n\la\m$, $1\leqslant i\leqslant n$ as above):
\begin{equation}\label{squarefortensorproductchange}
\begin{tikzcd}
	A_i&	A'_i\ar[l, "\alpha_i"']\\
	B_j\ar[u, "x_{ij}"']&	B'_j\ar[u, "x'_{ij}"']\ar[l, "\beta_j"']
\end{tikzcd}
\end{equation}
Applying the above procedure we obtain the commutative diagram
\begin{equation}\label{naturalityofgentensorproductjj}
\begin{tikzcd}[column sep=75]
	\prod_i^n\calV(A_i)\ar[r, "\prod_i\pV{\alpha_i}"]\ar[d, "\pr"]& 	\prod_i^n\calV(A'_i)\ar[d, "\pr"]\\
	\prod_j^m\prod_i^{f(j)}\calV(A_i)\ar[r, "\prod_j\prod_i\pV{\alpha_i}"]\ar[d, "\prod_j\prod_i\pV{x_{ij}}"]& 	\prod_j^m\prod_i^{f(j)}\calV(A'_i)\ar[d, "\prod_j\prod_i\pV{x'_{ij}}"]\\
	\prod_j^m\prod_i^{f(j)}\calV(B_j)\ar[r, "\prod_j\prod_i\pV{\beta_j}"]\ar[d, "\prod_j\otV{B_j}"]& 	\prod_j^m\prod_i^{f(j)}\calV(B'_j)\ar[d, "\prod_j\otV{B'_j}"]\\
	\prod_j^m\calV(B_j)\ar[r, "\prod_j\pV{\beta_j}"]& 	\prod_j^m\calV(B'_j)
\end{tikzcd}
\end{equation}
where the last square commutes by naturality of $\otV{}$. Hence we get the following commutative naturality square for the multi-valued tensor product:
\begin{equation}\label{naturalityofmultivaltensor}
\begin{tikzcd}[column sep=75]
		\prod_i^n\calV(A_i)\ar[r, "\prod_i\pV{\alpha_i}"]\ar[d, "{\fotV[x_{ij}]f}"]& 	\prod_i^n\calV(A'_i)\ar[d, "{\fotV[x'_{ij}]f}"]\\
	\prod_j^m\calV(B_j)\ar[r, "\prod_j\pV{\beta_j}"]& 	\prod_j^m\calV(B'_j)
\end{tikzcd}
\end{equation}
By applying the above to the commutative diagrams

\[\begin{tikzcd}
	A_i&	A_i\ar[l, "="']\\
	\prod_{i=f(j-1)+1}^{f(j)}A_i\ar[u, "\pr"']&	B_j\ar[u, "x_{ij}"']\ar[l, "x_j"']
\end{tikzcd}\]
we immediately obtain:

\begin{dCor}\label{altdescrofmultivaltensor}
We can compute the multi-valued tensor-product corresponding to $f\colon \n\la\m$ and $x_{ij}\colon B_j\to A_i$ as the composition
\begin{equation}
	\fotV[x_{ij}]f\colon \prod_{i=1}^n\calV(A_i)\xra{\fotV f}\prod_{j=1}^m\calV\(\prod_{i=f(j-1)}^{f(j)} A_i\)\xra{\prod_j\pV{x_j}}\prod_{j=1}^m\calV(B_j)\qedhere
\end{equation}
\end{dCor}

\subsubsection{Componentwise compatibility with homotopy colimits}
A particular case of the above discussion is when $f\colon \n\la\m$ is just the inclusion $\enu n\hla\{0,n\}\cong \1$ and Diagram~\ref{squarefortensorproductchange} is chosen to be
\begin{equation*}
	\begin{tikzcd}
		A_i&	A'_i\ar[l, "\alpha_i"']\\
		\prod_{i=1}^nA_i\ar[u, "\pr_i"']&	\prod_{i=1}^nA'_i\ar[u, "\pr_i"']\ar[l, "\prod_i \alpha_i"']
	\end{tikzcd}
\end{equation*}
In this case the naturality square of the multi-valued tensor product looks as follows
\begin{equation}\label{standardnaturalityofmultivaltensor}
\begin{tikzcd}[column sep=75]
		\prod_i^n\calV(A_i)\ar[r, "\prod_i\pV{\alpha_i}"]\ar[d, "{\fotV{n}}"]& 	\prod_i^n\calV(A'_i)\ar[d, "{\fotV{n}}"]\\
	\calV(\prod_i^nA_i)\ar[r, "\pV{\prod_i\alpha_i}"]& 	\calV(\prod_i^nA_i)
\end{tikzcd}
\end{equation}

\begin{Rem}\label{multitensorwithdiagonal} It is a straightforward calculation using the definitions that we can recover the $n$-ary tensor product from this multi-valued tensor product as the composition
\[\otV{A}\colon \prod_{i=1}^n\calV(A)\xra{\fotV{n}}\calV(\prod_{i=1}^n A)\xra{\pV {\Delta_n}} \calV(A)\]
where $\Delta_n\colon A\to \prod_{i=1}^n A$ is the diagonal.
\end{Rem}

\begin{Def}
	A monoidal left derivator is said to be \textbf{componentwise cocontinuous (cococo)} if the mate
	\begin{equation}\label{cococomate}
\begin{tikzcd}[column sep=75]
		\prod_i^n\calV(A_i)\ar[from=r, "\prod_i\sV{\alpha_i}"']\ar[d]\ar[d, "{\fotV{n}}"]& 	\prod_i^n\calV(A'_i)\ar[d, "{\fotV{n}}"]\\
	\calV(\prod_i^nA_i)\ar[from=r, "\sV{\prod_i\alpha_i}"']& 	\calV(\prod_i^nA_i)\ar[ul, Rightarrow]
\end{tikzcd}
\end{equation}
(corresponding to Diagram~\ref{standardnaturalityofmultivaltensor} after replacing the horizontal arrows with their left adjoint) is an isomorphism for all $n\in \BN$ and $A_1,\dots,A_n\in \Grpd$.
\end{Def}

\begin{Expl}
	What does being componentwise cocontinuous mean for the represented left derivator $[-,\calV]$?\\
	It is easy to compute that the map $\fotV{n}\colon \prod_i^n[A_i,\calV]\to[\prod_i^nA_i,\calV]$ is in this case given by
	\begin{align*}(F_1,\dots, F_n)&\mapsto& F_1\otimes \dots\otimes F_n\colon \prod\nolimits _i^nA_i\xra{\prod_iF_i}\prod_i^n\calV\xra{\bigotimes}\calV\end{align*}
	Hence the component at $(F'_1,\dots,F'_n)\in \prod_{i=1}^n[A_i',\calV]$ of the mate~\ref{cococomate} is nothing but the canonical natural transformation (of functors $A_1\times \dots\times A_n\to V$)
	\[(\alpha_1\times \dots\times \alpha_n)_\shriek(F'_1\otimes \dots\otimes F'_n)\lra {\alpha_1}_\shriek(F'_1)\otimes \dots\otimes {\alpha_n}_\shriek(F'_n)\]
	Evaluating this transformation at an element $(a_1,\dots,a_n)$ of $A_1\times \dots\times A_n$ we obtain the map
	\begin{equation*}
		\colim \({F'_1(a'_1)\otimes\dots\otimes  F'_n(a'_n)}\)\lra (\colim F'_1(a'_1))\otimes \dots\otimes (\colim F'_n(a'_n))
	\end{equation*}
	in $\calV$, where the colimits are taken over all arrows $\alpha_i(a'_i)\to a_i$ in $A_i$. That all these maps in $\calV$ are isomorphisms means precisely that the tensor product in $\calV$ preserves colimits in each variable.\\
	In particular we see that our main example $(\Vect\BC,\otimes)$ gives rise to a componentwise cocontinuous monoidal left derivator.
\end{Expl}

We conclude this ad hoc introduction of monoidal left derivators by proving an easy lemma.

\begin{Lem}\label{Lemhccriterionfortensor}
	Let $\calV$ be componentwise cocontinuous. Then the mate corresponding to Diagram~\ref{naturalityofmultivaltensor} after replacing the horizontal arrows by their left adjoint is an isomorphism if we assume that the square
	\begin{equation}\label{hccriterionfortensor}
\begin{tikzcd}[column sep=60]
	\prod_{i=f(j-1)+1}^{f(j)}A_i&	\prod_{i=f(j-1)+1}^{f(j)}A'_i\ar[l, "\prod_i\alpha_i"']\\
	B_j\ar[u, "x_{j}"']&	B'_j\ar[u, "x'_{j}"']\ar[l, "\beta_j"']
\end{tikzcd}
\end{equation}
	is $\calV$-exact.
\end{Lem}
\begin{Prf}
	We can factor $f\colon \n\la\m$ as $\n\hla\ol f(\m)\xla{f}\m\xla{\Id}\m$ and $x_{ij}\colon A_i\la B_j$ as $A_i\xla{\Id}A_i\xla{\pr_i}\prod\limits_{i=f(j-1)+1}^{f(j)} A_i\xla{x_j}B_j$ and similarly for $x'_{ij}$. (Recall that $\ol f(\m) =\{f(0),f(0)+1,\dots,f(m)\}=\{f(0)\}\sqcup\coprod_{j=1}^m\{f(j-1)+1,\dots, f(j)\}$ is the convex hull of $\Im f$.) This gives a decomposition of Diagram~\ref{naturalityofmultivaltensor} as
\begin{equation*}
	\begin{tikzcd}[column sep=75]
		\prod_i^n\calV(A_i)\ar[r, "\prod_i\pV{\alpha_i}"]\ar[d, "{\fotV[\Id]{\hla}}=\pr"]& 	\prod_i^n\calV(A'_i)\ar[d, "{\fotV[\Id]{\hla}}=\pr"]\\
		\prod_j^m\prod_i^{f(j)}\calV(A_i)\ar[r, "\prod_j\prod_i\pV{\alpha_i}"]\ar[d]\ar[d, "{\fotV f}"]& 	\prod_j^m\prod_i^{f(j)}\calV(A'_i)\ar[d, "{\fotV f}"]\\
		\prod_j^m\calV(\prod_i^{f(j)}A_i)\ar[d, "{\fotV[x_j]{\Id}}=\prod_j\pV{x_j}"]\ar[r, "\prod_j\pV{\prod_i \alpha_i}"]&\prod_j^m\calV(\prod_i^{f(j)}A'_i)\ar[d, "{\fotV[x'_j]{\Id}}=\prod_j\pV{x_j}"]\\
		\prod_j^m\calV(B_j)\ar[r, "\prod_j\pV{\beta_j}"]& 	\prod_j^m\calV(B'_j)
	\end{tikzcd}
\end{equation*}
	After replacing each horizontal arrow $\pV{?}$ with its left adjoint $\sV{?}$ the mate of the first square is always an isomorphism and the mate of the second square (which is just a product of instances of Diagram~\ref{cococomate}) is an isomorphism because $\calV$ is componentwise cocontinuous. Hence the mate of the total rectangle being an isomorphism follows from the mate of the lower square being an isomorphism, i.e.\ Diagram~\ref{hccriterionfortensor} being $\calV$-exact.
\end{Prf}

\ifx\isstandalone\undefined
\else
\newpage
\bibliographystyle{amsalpha}
\bibliography{literatur}
\end{document}
\fi

\subsection{The generalized Hall construction}\label{sectionmonconstruction}
\ifx\preambleloaded\undefined
   
   \def\isstandalone{}
\fi

Recall that our goal is to construct the Hall monoidal category out of the $2$-Segal simplicial groupoid $\calX^\calA\colon \Dop\to \Grpd$ of flags (in some proto-abelian category $\calA$) by applying the represented left derivator $[\blank, \Vect\BC]$. More generally we can use any $2$-Segal simplicial groupoid and any componentwise cocontinuous monoidal left derivator of groupoids and still get a monoidal category (i.e.\ a pointed $1$-Segal op-fibration of categories over $\Dop$).\\
It turns out that with no additional effort we can take care not only of simplicial groupoids but of any \hypo-simplicial groupoid. This might seem like a cheap improvement but it pays off: as we have seen in Section~\ref{moncatmonmod} we can formulate the notion of modules over a monoidal category in the language of \hypo-simplicial groupoids, hence this more general construction will provide us with a way to construct categorical modules over the Hall monoidal categories.\\
Let $\V\colon \Grpd\to\CAT$ be a componentwise cocontinuous monoidal left derivator of groupoids and fix an admissible subcategory $\Deltap\subseteq \Delta$.\\
In this section we construct a contravariant functor
\begin{equation*}
	\{\text{$2$-Segal \hypo-simpl. $\Dpop\to\Grpd$}\}\xra{\HallV{}}\{\text{pointed $1$-Segal op-fibr. over $\Dpop$}\}^{\lax}
\end{equation*}
which we call the \introduce{generalized Hall construction}. The superscript $\lax$ on the right side indicates that we only consider those functors of op-fibrations ${?\to\Dpop}$ which induce lax morphisms when we pass to pseudo-\hypo-simplicial categories ${\Dpop\to\CAT}$.

\subsubsection{Construct the data on objects}
We start by constructing the assignment $\HallV{}$ on objects. Let $\calX\colon \Dpop\to\Grpd$ be a \hypo-simplicial groupoid, we construct an op-fibration $\HallV\calX\colon \monF[\calX]\to\Dpop$ and prove that it is pointed $1$-Segal.

We construct the category $\monF=\monF[\calX]$:
\begin{itemize}
\item For each element $N^n\in\Dpop$ the objects of $\monF$ over $N$ are finite tuples $F=(F_1,\dots,F_n)$ of objects $F_i\in \calV(\calX_{N_i})$ (recall that if $N=\{x_0,\dots,x_n\}$ then $N_i\coloneqq \{x_{i-1},x_i\}$).

\item To define the Hom-sets and composition over arrows $N^n\xla{f}M^m\xla{g}L^l$ in $\Dpop$ we use the following cospans which assemble to a commutative diagram:
	\begin{equation}\label{compositiondiagram}
	\begin{tikzcd}[cramped,column sep=50]
		&&F\in \prod\limits_{i=1}^n\calV(\calX_{N_i})\ar[d, "\fotV f"]\ar[dd, "\fotV{fg}", bend left=70]\\
		&G\in \prod\limits_{j=1}^m\calV(\calX_{M_j})\ar[r, bend left=5, "\prod\limits_j\pVX{f\restriction M_j}"]\ar[d,"\fotV g"]&\prod\limits_{j=1}^m\calV(\calX_{\ol f(M_j)})\ar[d, "\fotV g"]\\
		H\in \prod\limits_{k=1}^l\calV(\calX_{L_k})\ar[r, bend left=5, "\prod\limits_k\pVX{g\restriction L_k}"]\ar[rr, bend right=15, "\prod\limits_k\pVX{(fg)\restriction L_k}"']&\prod\limits_{k=1}^l\calV(\calX_{\ol g(L_k)})\ar[r, bend left=5, "\prod\limits_k\pVX{f\restriction \ol g(L_k)}"]&\prod_{k=1}^l\calV(\calX_{\ol {fg}(L_k)})
	\end{tikzcd}
	\end{equation}
	The vertical arrows are the multi-valued tensor product defined in Section~\ref{multivaluedtensorprod} using the structure maps coming from various inclusions e.g. $N_i\hra \ol f(M_j)$ or $\ol f(M_j)\hra \ol{fg}(L_k)$.
	Note that the square in the lower right is an instance of the naturality square~\ref{naturalityofmultivaltensor} for $\fotV{g}$.\\
	More explicitly, we define
	\[\Hom_f(F,G)\coloneqq \Hom_{\prod\limits_{j=1}^m\calV(\calX_{\ol f(M_j)})}\(\fotV{f}(F), \,\,\(\prod\limits_j\pVX{f\restriction M_j}\)(G)\)\]
	and composition $\Hom_g(G,H)\times\Hom_f(F,G)\to\Hom_{fg}(F,H)$ by pushing the homsets forward to the category $\prod_{k=1}^l\calV(\calX_{\ol {fg}(L_k)})$ (in the bottom right corner) and using the composition there.
\item Finally, we define the Hom-sets $\Hom(F,G)\coloneqq\coprod_{f\colon N\la M}\Hom_f(F,G)$. Clearly the above composition law extends. Moreover we obtain the obvious functor $\monF\to \Dpop$.
\end{itemize}
	Note that the fiber $\monF[N]$ over $N\in \Dpop$ is clearly $\prod_{i=1}^n\calV(\calX_{N_i})$.

\subsubsection{Check the properties}\label{checkpropertiesofgenHall}
It is clear that the above construction of $\mon F$ really gives a category: associativity of composition is obtained by assembling three cospans to a commutative diagram and identities are easy since the maps $\fotV{f}$ and $\prod\limits_j\pVX{f\restriction M_j}$ are both just the identity if $f=\Id$. By the construction of the composition in $\monF$ the map $ \monF\to \Dpop$ respects it; hence it can rightfully be called a functor.\\
We are left to check that this functor $\HallV\calX\colon \monF\to\Dpop$ is an op-fibration of categories and that the corresponding \hypo-simplicial category $\monF\colon \Dpop\to\CAT$ is pointed $1$-Segal.
\begin{itemize}
\item[\ref{opfibration}] We start with an object $F$ in the fiber $\monF[N]$ and a map $f\colon N\la M$ in $\Dpop$. We have to construct an object $G$ over $ M$ and a morphism $\varphi\colon F\to G$ in $\monF$ over $f$ such that for every $g\colon M\la L$ in $\Dpop$ the induced morphism
\begin{equation}
	\Hom_g(G,H)\xra{\varphi^\star}\Hom_{fg}(F,H)
\end{equation}
is an isomorphism.\\
We take Diagram~\ref{compositiondiagram} and replace all horizontal arrows by their left adjoints filling the formerly commuting areas with mates:
	\begin{equation}\label{compositiondiagramreversed}
	\begin{tikzcd}[cramped,column sep=50]
		&&F\in \prod\limits_{i=1}^n\calV(\calX_{N_i})\ar[d, "\fotV f"]\ar[dd, "\fotV{fg}", bend left=70]\ar[ld, ""{name=bla, above=31, right=31}, ""{name=blub, below=21, left=101}, phantom]\\
		&G\in \prod\limits_{j=1}^m\calV(\calX_{M_j})\ar[from=r, bend right=5, "\prod\limits_j\sVX{f\restriction M_j}"', ""{name=C, below=10, very near end}]\ar[d,"\fotV g"]&\prod\limits_{j=1}^m\calV(\calX_{\ol f(M_j)})\ar[d, "\fotV g", ""{name=B, left=50}]\\
		H\in \prod\limits_{k=1}^l\calV(\calX_{L_k})\ar[from=r, bend right=5, "\prod\limits_k\sVX{g\restriction L_k}"']\ar[from=rr, bend left=25, "\prod\limits_k\sVX{(fg)\restriction L_k}", ""{name=A, above}]&\prod\limits_{k=1}^l\calV(\calX_{\ol g(L_k)})\ar[from=r, bend right=5, "\prod\limits_k\sVX{f\restriction \ol g(L_k)}"']\ar[from=A, Leftrightarrow]&\prod_{k=1}^l\calV(\calX_{\ol {fg}(L_k)})
		\ar[from=B, to=C, Rightarrow]
		\ar[from=bla, to=blub, bend right, mapsto, dashed]
	\end{tikzcd}
	\end{equation}
	
	The lower mate is always an isomorphism (see Appendix~\ref{calculusofmates}). By Lemma~\ref{Lemhccriterionfortensor} the mate in the upper square is an isomorphism if all the squares (for $1\leqslant k\leqslant l$)
	\begin{equation}\label{criterioncheckingblabla}
	\begin{tikzcd}[column sep=70]
		\prod_{j}\calX_{M_j}&\prod_j\calX_{\ol f(M_j)}\ar[l, "\prod_j\X_{f\restriction M_j}"']\\
		\calX_{\ol g(L_k)}\ar[u]&\calX_{\ol{fg}(L_k)}\ar[u]\ar[l, "\X_{f\restriction \ol g(L_k)}"']
	\end{tikzcd}
	\end{equation}
	are $\calV$-exact; here the products in the upper row are indexed by those $1\leqslant j\leqslant m$ such that $M_j\subseteq \ol g(L_k)$ (hence $\ol f(M_j)\subseteq \ol{fg}(L_k)$), i.e.\ those $j$ such that the vertical maps are defined. The squares~\ref{criterioncheckingblabla} are instances of the $2$-Segal condition~\ref{generalunitalsubsimplicial} (applied to the map $\restr f{\ol g(L_k)}\colon N\la\ol g(L_k)$) for the \hypo-simplicial groupoid $\X\colon \Dpop\to\Grpd$, hence they are $2$-pullback, hence $\calV$-exact by Axiom~\ref{2pullbacksarehc}.\\
	Therefore all mates in Diagram~\ref{compositiondiagramreversed} are isomorphisms.\\
	As the diagram suggests, we define $G\coloneqq \(\prod_j\sVX{f\restriction}\circ \fotV{f}\)(F)$. Then we calculate
	\begin{equation}\label{calculationforinducedonfibers}\begin{aligned}
		\Hom_g(G,H)&=\Hom_g\(\(\prod\nolimits _j\sVX{f\restriction}\circ \fotV{f}\)(F), H\)\\
		&\coloneqq \Hom\(\(\fotV{g}\circ \prod\nolimits _j\sVX{f\restriction}\circ \fotV{f}\)(F), \(\prod_k\pVX{g\restriction}\)H\)\\
		&\cong \Hom\(\(\prod\nolimits_k\sVX{f\restriction}\circ \fotV{g}\circ \fotV{f}\)(F), \(\prod_k\pVX{g\restriction}\)(H)\)\\
		&\cong \Hom\(\(\fotV{g}\circ \fotV{f}\)(F), \(\prod\nolimits_k\pVX{f\restriction}\circ \prod_k\pVX{g\restriction}\)(H)\)\\
		&\cong  \Hom\(\fotV{fg}(F), \prod_k\pVX{(fg)\restriction}(H)\)\\
		&\eqqcolon \Hom_{fg}(F,H)
	\end{aligned}\end{equation}
	naturally in $H$. The map $\varphi\colon F\to G$ is now obtained as the image of $\Id_G$ under this isomorphism $\Hom_\Id(G,G)\xra{\cong}\Hom_f(F,G)$ specializing $H$ to $G$. It is straightforward to check, using the naturality of the above isomorphism, that this isomorphism $\Hom_g(G,H)\xra{\cong}\Hom_{fg}(F,H)$ is indeed given as $\varphi^\star$. This completes the proof that $p$ is an op-fibration of categories.
	\end{itemize}
\begin{Rem}\label{descriptioninducedonfibers}
Let $f\colon N^n\la M^m$ be a map $\in \Dpop$. The isomorphism we obtained from the calculation~\ref{calculationforinducedonfibers} is not only natural in $H$ but also in $F$. Therefore the functor 
\begin{equation}
	\prod_{i=1}^n\calV(\calX_{N_i})\xra{\fotV{f}}\prod_{j=1}^m\calV(\calX_{\ol f(M_j)})\xra{\prod_j\sVX{f\restriction M_j}} \prod_{j=1}^m\calV(\calX_{M_j})
\end{equation}
does not only describe the assignment $F\mapsto G$ on objects but in fact describes the functor $f^\star\colon \monF[N]\to\monF[M]$ induced by $f$ on fibers.
\end{Rem}
	\begin{itemize}
	\item[\ref{suboneSegal}] By construction the fiber over (the identity on) $N^n\in \Dpop$ is
	\[\monF[N]\coloneqq \prod_{i=1}^n\calV(\calX_{N_i})\eqqcolon\prod_{i=1}^n \monF[N_i]\]
	so we have to show that for $1\leqslant i\leqslant n$ the projection functor $\pr_i\colon\monF[N]\to\monF[N_i]$ is indeed induced by the inclusions $f_i\colon N\hla N_i$. This is immediate from the description in Remark~\ref{descriptioninducedonfibers}\end{itemize}

\subsubsection{Functoriality}\label{sectionHallfunctoriality}
Using the naturality of multi-valued tensor products it is immediate that Diagram~\ref{compositiondiagram} is (contravariantly) functorial in $\calX$ by sending a morphism $\alpha\colon \calY\to\calX$ to the induced morphisms $\prod\pV{\alpha_?}$. All the definition of objects, morphisms, composition in $\monF[\calX]$ as well as the fiber functor $\HallV\calX\colon \monF[\calX]\to\Dpop$ are encoded in Diagram~\ref{compositiondiagram}; hence our construction is automatically equipped with a functoriality $\alpha\mapsto \pV\alpha\eqqcolon\HallV\alpha$.\\
All in all we obtain a functor
\[\HallV{}\colon \{\text{$2$-Segal \hypo-simpl. $\Dpop\to\Grpd$}\}\lra\{\text{pointed $1$-Segal op-fibr. over $\Dpop$}\}\]
We still need to see how the functors $\HallV{\alpha}$ induce lax morphisms when we pass from opfibrations $\monF[\calX]\to\Dpop$ to pseudo-\hypo-simplicial objects $\monF[\calX]\colon \Dpop\to\CAT$.\\
Using Remark~\ref{descriptioninducedonfibers} can build the naturality square~\ref{laxsquareforbeta} (associated to a morphism $f\colon N\la M$ of $\Dpop$) as
\begin{equation}\label{compositenaturalitysquareforfiber}
	\begin{tikzcd}[column sep=huge]
		\prod_{i=1}^n\calV(\calX_{N_i})\rar{\fotV{f}}\ar[d,"\prod_{i}\pV{\alpha}"]&\prod_{j=1}^m\calV(\calX_{\ol f(M_j)})\rar{\prod_j\sVX{f\restriction M_j}}\ar[d,"\prod_j\pV{\alpha}"'] &\prod_{j=1}^m\calV(\calX_{M_j})\ar[d,"\prod_j\pV{\alpha}"]\\
				\prod_{i=1}^n\calV(\calY_{N_i})\rar{\fotV{f}}&\prod_{j=1}^m\calV(\calY_{\ol f(M_j)})\rar{\prod_j\sVY{f\restriction M_j}}\ar[ru,Rightarrow] &\prod_{j=1}^m\calV(\calY_{M_j})
	\end{tikzcd}
\end{equation}
where the first square commutes and the second square is inhabited by the canonical mate associated to the naturality square of $\alpha$.\\
If the morphism $f$ is convex then the maps $f_{M_j}\colon M_j\to\ol f(M_j)$ are isomorphisms; this implies that the mate in the right square of Diagram~\ref{compositenaturalitysquareforfiber} is an isomorphism.\\
To see the compatibility of Diagram~\ref{compositenaturalitysquareforfiber} with horizontal pasting (for two composable morphisms $N\xla{f}M\xla{g}L$ in $\Dop$) we take the commutative cube
	\begin{equation}\label{blablablsdkfjharwe}
	\begin{tikzcd}
	&\prod_{k=1}^l\calV(\calX_{\ol {fg}(L_k)})\ar[dd,"\prod_{j}\pV{\alpha}"near end]\ar[from=dl,"\fotV g"']&&\prod\limits_{k=1}^l\calV(\calX_{\ol g(L_k)})\ar[ll, "\prod\limits_k\pVX{f\restriction \ol g(L_k)}"']\ar[from=dl,"\fotV g"']\ar[dd,"\prod_{j}\pV{\alpha}"]\\
	\prod\limits_{j=1}^m\calV(\calX_{\ol f(M_j)})\ar[dd,"\prod_{j}\pV{\alpha}"]&&\prod\limits_{j=1}^m\calV(\calX_{M_j})\ar[ll, "\prod\limits_j\pVX{f\restriction M_j}"near end]\ar[dd,"\prod_{i}\pV{\alpha}" near start]\\
	&\prod_{k=1}^l\calV(\calY_{\ol {fg}(L_k)})\ar[from=dl,"\fotV g"']&&\prod\limits_{k=1}^l\calV(\calY_{\ol g(L_k)})\ar[ll, "\prod\limits_k\pVX{f\restriction \ol g(L_k)}"near end]\ar[from=dl,"\fotV g"']\\
	\prod\limits_{j=1}^m\calV(\calY_{\ol f(M_j)})&&\prod\limits_{j=1}^m\calV(\calY_{M_j})\ar[ll, "\prod\limits_j\pVX{f\restriction M_j}"']
	\end{tikzcd}
	\end{equation}
	and replace the horizontal arrows by their left adjoints and pass to mates.\\
	Pasting two copies of Diagram~\ref{compositenaturalitysquareforfiber} (one for $f$ and one for $g$) corresponds to taking the front and right face of Diagram~\ref{blablablsdkfjharwe}; taking Diagram~\ref{compositenaturalitysquareforfiber} for the composition $fg$ corresponds to taking the left and back face of Diagram~\ref{blablablsdkfjharwe}. Hence Diagram~\ref{compositenaturalitysquareforfiber} is compatible with pasting up to the top and bottom faces of Diagram~\ref{blablablsdkfjharwe}.\\
	The mates of the top and bottom faces are precisely those isomorphisms which we use to identify $\monF[g]\circ \monF[f]\cong \monF[fg]$, thus the proof is concluded.

\ifx\isstandalone\undefined
\else
\newpage
\bibliographystyle{amsalpha}
\bibliography{literatur}
\end{document}
\fi


\newpage\section{Hall monoidal categories and categorical modules}
\ifx\preambleloaded\undefined
   
   \def\isstandalone{}
\fi

We are mainly interested to specializing the generalized Hall construction from Section~\ref{sectionmonconstruction} to $\Deltap\coloneqq \Delta$ and $\Deltap\coloneqq \Delta_\mindecor$ (or $\Deltap\coloneqq \Delta_\maxdecor$ ). In these cases pointed $1$-Segal op-fibrations over $\Dpop$ are nothing but monoidal categories and monoidal right (or left) modules respectively. We can describe the monoidal product and the action product more explicitly.\\

Let $\calX\colon \Dop\to\Grpd$ be a $2$-Segal simplicial groupoid and let $\calX^\mindecor\colon\Dnop\to\Grpd$ be a $2$-Segal relative simplicial groupoid such that $\calX$ agrees with the simplicial groupoid
\[\calX^\mindecor_{\not\in}\colon \Dop\simeq\Delta^\op_{\not\in}\subset\Dnop\xra{\calX^\mindecor}\Grpd.\]

\begin{ThmDef}\label{descriptionHallmon} Let $\calV\colon \Grpd^\op\to\CAT$ be a componentwise cocontinuous monoidal left derivator of groupoids.
\begin{enumerate}
	\item\label{descriptionHallmonpartHallcat} The structure $\HallV\calX$ makes the underlying category $\calV(\calX_\1)$ into a monoidal category. The monoidal product $\boxtimes$ is given by the composition
		\begin{equation}\label{explicitmonoidalproduct}
			\boxtimes\colon \calV(\calX_\1)\times\calV(\calX_\1)\xra{\fotV{2}}\calV(\calX_\1\times\calX_\1)\xra{\pV{c_2}}\calV(\calX_\2)\xra{\sV{\nu_2}}\calV(\calX_\1)
		\end{equation}
		where the first functor is the multi-valued tensor product of $\calV$ and the rest arises by pull-pushing along the multiplication span
	\[\calX_{\{0,1\}}\times\calX_{\{1,2\}}\xla{c_2}\calX_\2\xra{\nu_2}\calX_{\{0,2\}}.\]
	The monoidal unit is given by the composition
	\begin{equation}\label{explicitmonoidalunit}
		I\colon \terminal\xra{\fotV{0}}\calV(\terminal)\xra{\pV{c_0}}\calV(\calX_\0)\xra{\sV{\nu_0}}\calV(\calX_\1).
	\end{equation}
	We call $\HallV\calX$ the \introduce{$\calV$-Hall monoidal category} of $\calX$.
	\item The structure $\HallV{\calX^\mindecor}$ makes the underlying category $\calV(\calX^\mindecor_{\{\minimal,0\}})$ into a categorical (right-)module over the monoidal category $\HallV{\calX}$ (which has underlying category $\calV(\calX_\1)=\calV(\calX^\mindecor_\1)$). The action product is given by the composition
	\begin{equation}\label{explicitactionproduct}
		\rcaction\colon \calV(\calX^\mindecor_{\{\minimal,0\}})\times\calV(\calX_\1)\xra{\fotV{}}\calV(\calX^\mindecor_{\{\minimal,0\}}\times\calX_\1)\xra{\pV{c}}\calV(\calX^\mindecor_{\{\minimal,0,1\}})\xra{\sV{\nu}}\calV(\calX^
		\mindecor_{\{\minimal,1\}})
	\end{equation}
	We call $\HallV{\calX^\mindecor}$ the \introduce{categorical $\calV$-Hall (right-)module} (over $\HallV \calX$) corresponding to $\calX^\mindecor$.
\end{enumerate}
In the special case $\calV=[\blank, \Vect\BC]$ we drop the $\calV$ from the notation and speak of the Hall monoidal category $\Hall\calX$ and categorical Hall (right-)module $\Hall{\calX^\mindecor}$.
\end{ThmDef}
\begin{Prf}
We will focus on the monoidal product in part~\ref{descriptionHallmonpartHallcat}. The rest is completely analogous.\\
We already constructed the monoidal category $\HallV\calX$ as a pointed $1$-Segal op-fibration $\monF=\monF[\calX]\to\Dop$; all we need to do is determine the monoidal product
\[\boxtimes_2\colon \monF[\1]\times \monF[\1]\lra\monF[\1]\] and see that it agrees with the one described by the composition~\ref{explicitmonoidalproduct}.
The monoidal product $\boxtimes_2$ arises by inverting the first map in the multiplication span
\[\monF[\1]\times\monF[\1]\xla{\approxeq}\monF[\2]\xra {}\monF[\1]\]
and we already know that this first map can be taken to be the identity if we identify $\monF[\2]$ with $\monF[\1]\times\monF[\1]$ as in Section~\ref{checkpropertiesofgenHall}.\\
As we saw in Remark~\ref{descriptioninducedonfibers}, the functor $\monF[\2]\to\monF[\1]$ (induced by $\2\hla{\{0,2\}}$) is given by the composition
\[\VX{0,1}\times\VX{1,2}\xra{\fotV[c_2]{}} \VX{0,1,2}\xra{\sV{\nu}}\VX{0,2},\]
which is equal to Composition~\ref{explicitmonoidalproduct} by Corollary~\ref{altdescrofmultivaltensor}.
\end{Prf}

\begin{Def}
	If $\calA$ is a proto-abelian category then we define its \introduce{Hall monoidal category} as $\Hall\calA\coloneqq \Hall{\calX^\calA}$, where $\calX^\calA$ is the $2$-Segal simplicial groupoid of flags in $\calA$ defined in Section~\ref{defgroupoidofflags}.
\end{Def}

\subsection{Examples}

\begin{Expl}[The regular action]
If $\calX\colon \Dop\to\Grpd$ is a $2$-Segal simplicial groupoid then we can restrict it to a right relative simplicial groupoid
\[\calX^\mindecor\coloneqq \restr\calX\Dnop\colon \Dnop\subset\Dop\to\Grpd,\]
which is still $2$-Segal.\\
If we view $\calX^\mindecor$ as a morphism of simplicial groupoids $\restr[\in]\calX\Dnop\to\restr[\not\in]\calX\Dnop$ then we obtain the canonical map $P^\triangleleft\calX\to\calX$ from the \buzzword{initial path space}~\cite[Section 6.2]{DyckerhoffKapranovHigherSegalSpacesI}.\\
Since in this case the action span for $\calX^\mindecor$ is canonically isomorphic to the multiplication span for $\calX$ (via the isomorphisms $\calX_{\{\minimal,0,\dots,n-1\}}\cong\calX_\n$ which we had forgotten when passing from $\Dop$ to $\Dnop$) we see that the resulting categorical module $\HallV{\calX^\mindecor}$ is nothing but the right-regular action of the monoidal category $\HallV\calX$ on itself.
\end{Expl}

\begin{Expl}
Let $\calX\colon \Dop\to\Grpd$ be an unpointed $1$-Segal groupoid which is hence $2$-Segal by Corollary~\ref{oneSegalgivesunital}. The identity map $\calX\to\calX$ is easily seen to be $2$-Segal using Proposition~\ref{criterionrelativeunital}, since squares~\ref{cdrel2segalgroupoid} and~\ref{cdrelunitalgroupoid} are in this case simply the $2$-Segal squares for $\calX$. Thus we get a categorical action of $\HallV{\calX}$ on the category $\calX_\0$.\\
\end{Expl}

\ifx\isstandalone\undefined
\else
\newpage
\bibliographystyle{amsalpha}
\bibliography{literatur}
\end{document}
\fi

\subsection{The finitary case and decategorification}\label{sectionfinitarycase}
\ifx\preambleloaded\undefined
   
   \def\isstandalone{}
\fi

In many cases it is important to impose some finiteness conditions; for instance we might want to pass to Grothendieck groups/rings later and these groups/rings might just vanish if the monoidal category is too big.

\begin{Def}
Fix $n\in\BN$. A \hypo-simplicial groupoid $\calX\colon \Dpop\to\Grpd$ is called \introduce{$n$-finitary} if for every object $\{x_0<\dots<x_i\}=N^i\in \Dpop$ of dimension $i\leqslant n$
\begin{itemize}
\item the chopping map $c_N\colon \calX_{N_1}\times\dots\times \calX_{N_i}\la\calX_N$ is $\pi_0$-finite
\item the extremal map $\nu_N\colon\calX_N\to\calX_{\{x_0,x_i\}}$ is locally finite.
\end{itemize}
If $\calX$ is $n$-finitary and the extremal maps $\nu_N$ (for $\dim N\leqslant n$) are faithful in addition to being locally finite then we say $\calX$ is \introduce{$n$-integral}.\\
If additionally the groupoid $\calX_{\{x_0,x_1\}}$ is locally finite (i.e.\ has finite automorphism groups) for all $\{x_0,x_1\}=N^1\in \Dpop$ then we call $\calX$ \introduce{$n$-regular}. \\
We will write \introduce{finitary}, \introduce{integral} and \introduce{regular} when we mean $2$-finitary, $2$-integral and $2$-regular respectively.
\end{Def}

\begin{Rem}
In the case $N^0=\{x_0\}$ the extremal map $\nu_N$ has to be read as the degeneracy $\nu_{\{x_0\}}\colon \xX{x_0}\to\xX{x_0,x_0'}$ and the chopping map $c_N$ is just the unique map $c_{\{x_0\}}\colon \calX_{\{x_0\}}\to\terminal$ to the point. Hence $\calX_0$ is $\pi_0$-finite if $\calX$ is finitary.
\end{Rem}

\begin{Rem}
If $\calX$ is $n$-regular and $i\leqslant n$, then for any $\{x_0<\dots<x_i\}=N^i\in \Dpop$ the automorphism groups of $\calX_N$ embed into the automorphism groups of $\calX_{\{x_0,x_i\}}$ via the faithful map $\nu_N$; hence the groupoid $\calX_N$ is also locally finite.
\end{Rem}

\begin{Expl}
We have seen in Proposition~\ref{HallAlgsrevisited} that the simplicial groupoid $\calX^\calA$ of flags in a finitary proto-abelian category $\calA$ is integral (recall that $\calX^\calA_\0\cong\terminal$). It is clearly also regular since $\calX^\calA_\1\cong\isoA$.
\end{Expl}

\begin{Def}
Let $A$ be a groupoid. A functor $A\to \Vect\BC$ is called \introduce{finitary} if it takes values in \emph{finite dimensional} vector spaces and if it is non-zero only on finitely many isomorphism classes of objects in $A$.
We denote by $[A,\vect\BC]_\finitary\subset [A,\Vect\BC]$ the full subcategory of finitary functors.
\end{Def}

\begin{PropDef}\label{restrictiontofinitary}
\begin{enumerate}
\item If $\calX$ is finitary then the monoidal unit $I\in [\calX_\1,\Vect\BC]$ is a finitary functor and the product
\[\boxtimes\colon [\calX_\1,\Vect\BC]\times[\calX_\1,\Vect\BC]\to[\calX_\1,\Vect\BC]\]
restricts to a monoidal product on $[\calX_\1,\vect\BC]_\finitary$.\\
We call the induced monoidal subcategory $\Hallf{\calX}$ of $\Hall{\calX}$ the \textbf{finitary Hall monoidal category} of $\calX$.
\item If $\calX^\mindecor$ is finitary then the action product
\[\rcaction\colon [\calX^\mindecor_{\{\minimal,0\}},\Vect\BC]\times[\calX_\1,\Vect\BC]\to[\calX^\mindecor_{\{\minimal,0\}},\Vect\BC]\]
restricts to a categorical right action of $[\calX_\1,\vect\BC]_\finitary$ on $[\calX^\mindecor_{\{\minimal,0\}},\vect\BC]_\finitary$.\\
We call the resulting $\Hall{\calX}$-module $\Hall{\calX^\mindecor}$ the \textbf{finitary categorical Hall (right-)module} associated to $\calX^\mindecor$.\qedhere
\end{enumerate}
\end{PropDef}

\begin{Prf}Proposition~\ref{restrictiontofinitary} follows immediately from the following Lemma~\ref{Lemmafinitaryfunctorspreserved} by using the explicit description for $I$, $\boxtimes$ and $\rcaction$ given in Proposition~\ref{descriptionHallmon}.
\end{Prf}
\begin{Lem}\label{Lemmafinitaryfunctorspreserved}
Let $f\colon A\to B$ be a map of groupoids.
\begin{enumerate}
	\item If $f$ is $\pi_0$-finite then the induced pullback map $f^\star\colon [B,\Vect\BC]\to[A,\Vect\BC]$ sends finitary functors to finitary functors.
	\item If $f$ is locally finite then the induced left Kan extension functor (pushforward) $f_\shriek\colon [A,\Vect\BC]\to[B,\Vect\BC]$ sends finitary functors to finitary functors.\qedhere
\end{enumerate}
\end{Lem}

\begin{Prf}
The first part is obvious so let us focus on the second part. Assume that $f\colon A\to B$ is locally finite and that $F\colon A\to\vect\BC$ is finitary. We can compute $f_\shriek F$ on an object $b\in B$ via the pointwise formula for Kan extensions:
	\begin{align*}
	(f_\shriek F)(b)&=\colim\(\twofib[f]Ab\to A\xra{F}\vect\BC\)\\
		&=\bigoplus_{[a]\in \pi_0 A}\colim\(\twofib[f]{A(a)}b\to A(a)\xra{F}\vect\BC\)
	\end{align*}
	Each of the summands is a \emph{finite} colimit (since the $2$-fiber $\twofib[f]{A(a)}b$ is finite by assumption) of finite dimensional vector spaces, hence finite dimensional. Since $F$ is non-zero only on finitely many isomorphism classes $[a_1],\dots,[a_r]$, we are left with finitely many summands, hence the vector space $(f_\shriek F)(b)$ is finite dimensional.\\
	Moreover, the $2$-fiber $\twofib[f]{A(a)}b$ is only non-trivial if $b\cong f(a)$, hence $f_\shriek F$ is supported on the (finitely many!) isomorphism classes $[f(a_1)],\dots,[f(a_r)]$.
\end{Prf}

Recall the assignment $V\colon \Grpd\to\lMod\BZ$ given by $A\mapsto \BZ^{(\pi_0A)}$ which came equipped with pullbacks $f^\star$ or pushforwards $f_\shriek$ whenever the map $f$ of groupoids was $\pi_0$-finite or locally finite and faithful respectively.\\
Proposition~\ref{restrictiontofinitary} and Proposition~\ref{descriptionHallmon} should be compared to the following analogous proposition for the Hall algebra.

\begin{PropDef}\label{descriptionHallalg}
\begin{enumerate}
	\item If $\calX$ is integral then the multiplication
\[\mu\colon V(\calX_\1)\times V(\calX_\1)\xra{\cdot} V(\calX_\1\times\calX_\1)\xra{c^\star}V(\calX_\2)\xra{\nu_\shriek}V(\calX_\1),\]
endows the abelian group $V(\calX_\1)$ with the structure of an associative and unital $\BZ$-algebra, where the unit is given by the composition
\[\varepsilon\colon \terminal\xra{1}\BZ\cong V(\terminal)\xra{c^\star}V(\calX_\0)\xra{\nu_\shriek}V(\calX_\1).\]
We call this algebra $\hall\calX$ the \introduce{Hall algebra} of $\calX$.
\item If $\calX^\mindecor$ is integral then the action
\[V(\calX^\mindecor_{\{\minimal,0\}})\times V(\calX_\1)\xra{\cdot} V(\calX^\mindecor_{\{\minimal,0\}}\times\calX_\1)\xra{c^\star}V(\calX^\mindecor_{\{\minimal,0,1\}})\xra{\nu_\shriek}V(\calX^\mindecor_{\{\minimal,1\}})
\]
endows the abelian group $V(\calX^\mindecor_{\{\minimal,0\}})\cong V(\calX^\mindecor_{\{\minimal,1\}})$ with a right-$\hall\calX$-module structure. We call this module $\hall{\calX^\mindecor}$ the \introduce{Hall module} associated to $\calX^\mindecor$.\qedhere
\end{enumerate}
\end{PropDef}

\begin{Prf}
We have already seen the first part in Section~\ref{classHallalgrevisited}, the rest is similar.
\end{Prf}

\subsubsection{The fat Hall algebra and the dimension map}
Let $\calX\colon \Dop\to\Grpd$ be a finitary $2$-Segal simplicial groupoid.
The category $[\calX_\1,\vect\BC]_\finitary$ is additive (since $\vect\BC$ is). Since the functor $\otimes\colon [\calX_\1,\vect\BC]\times [\calX_\1,\vect\BC]\to[\calX_\1\times\calX_\1,\vect\BC]$ preserves direct sums in each variable and Left Kan extensions also commute with direct sums we see from the description in Proposition~\ref{descriptionHallmon} that the monoidal product $\boxtimes\colon [\calX_\1,\vect\BC]_\finitary\times[\calX_\1,\vect\BC]_\finitary\to[\calX_\1,\vect\BC]_\finitary$ preserves direct sums in each variable.\\
Hence passing to the Grothendieck group we obtain an associative algebra
\[\fHall\calX\coloneqq \K_0\Hallf{\calX}\]
which we call the \introduce{fat Hall algebra} of $\calX$. If $R$ is any commutative ring we can consider an $R$-linear version and define $\fHall[R]{\calX}\coloneqq R\otimes_\BZ\fHall{\calX}$\\

The obvious question arising now is: What is the relationship between the fat Hall algebra and the classical Hall algebra? To answer this question we construct the \introduce{dimension map}
\[\cm_\calX= \cm_{\calX_\1}\colon \K_0\Hallf\calX\lra\hall\calX.\]
which is given by mapping a finitary functor $\rho\colon \calX_\1\to\vect\BC$ to the function on $\pi_0\calX_\1$ given by $[a]\mapsto \dim_\BC\rho(a)$. The map $\kappa_{\calX_\1}$ is a well defined $\BZ$-module homomorphism because $\dim_\BC$ is additive.\\

We can go through the same motions again for a relative simplicial groupoid $\calX^\mindecor\colon \Dnop\to\Grpd$ and define the \introduce{fat Hall (right-)module} $\fHall[R]{\calX^\mindecor}\coloneqq R\otimes_\BZ K_0\Hallf{\calX^\mindecor}$ and the \introduce{dimension map} $\cm_{\calX^\mindecor}=\cm_{\calX^\mindecor_{\{\minimal,0\}}}\colon \fHall{\calX^\mindecor}\lra \hall{\calX^\mindecor}$.

\begin{Prop}\label{Propcomparisonmap}
\begin{enumerate}
\item\label{cmisalghomo} If $\calX$ is integral then the dimension map $\cm_\calX\colon \fHall{\calX}\to\hall{\calX}$ is an unital $\BZ$-algebra homomorphism.
\item\label{cmismodhomo} If $\calX^\mindecor$ is integral then the dimension map $\cm_{\calX^\mindecor}\colon \fHall{\calX^\mindecor}\to\hall{\calX^\mindecor}$ is a $\fHall{\calX^\mindecor_{\not\in}}$-module homomorphism, where the algebra $\fHall{\calX^\mindecor_{\not\in}}$ acts on $\hall{\calX^\mindecor}$ via the algebra homomorphism $\cm\colon \fHall{\calX^\mindecor_{\not\in}}\to \hall{\calX^\mindecor_{\not\in}}$.\qedhere
\end{enumerate}
\end{Prop}

Before we can prove Proposition~\ref{Propcomparisonmap} we need to define the dimension map more generally: For any groupoid $A$ we define $\cm_A\colon [A,\vect\BC]_\finitary\to V(A)$ by the same formula $\cm_A(\rho)[a]\coloneqq \dim_\BC\rho(a)$. The next Lemma shows that the collection $\{\cm_A\}_{A\in \Grpd}$ is \roughly{natural enough} for our purposes.

\begin{Lem}\label{dimisnaturalenough}
	Let $f\colon A\to B$ be a map of groupoids.
	\begin{enumerate}
		\item\label{comparisonpullnatural} If $f$ is $\pi_0$-finite then we have the following commutative diagram of abelian groups
			\begin{equation}
				\begin{tikzcd}
					K_0[A,\vectC]_\finitary\ar[r,"\cm_A"]&V(A)\\
					\K_0[B,\vectC]_\finitary\ar[r,"\cm_B"]\ar[u,"\pV f"]&V(B)\ar[u,"f^\star"']
				\end{tikzcd}
			\end{equation}
		\item\label{comparisonpushnatural} If $f$ is locally finite and faithful then we have the following commutative diagram of abelian groups
			\begin{equation}\label{pushforwardcompsquare}
				\begin{tikzcd}
					\K_0[A,\vectC]_\finitary\ar[d,"\sV f"']\ar[r,"\cm_A"]&V(A)\ar[d,"f_\shriek"]\\
					\K_0[B,\vectC]_\finitary\ar[r,"\cm_B"]&V(B)
				\end{tikzcd}\qedhere
			\end{equation}
	\end{enumerate}
\end{Lem}

\begin{Rem}
If $H$ is a finite group then the map $f\colon \B H\to \triv$ is locally finite but Diagram~\ref{pushforwardcompsquare} reads (after tensoring with $\BQ$)
	\[\begin{tikzcd}
					\BQ\otimes_\BZ\K_0\lrep H\BC\ar[d,"\BQ\otimes_\BZ\K_0\(\BC\otimes_{\BC[H]}\blank\)"']\ar[r,"\dim_\BC"]&\BQ\ar[d,"\frac{1}{\#H}"]\\
					\BQ\otimes_\BZ\K_0\vect\BC\ar[r,"\dim_\BC"]&\BQ
	\end{tikzcd}\]
	which does not commute (try the trivial representation of $H$). Hence we really need the faithfulness assumption in the second part of Lemma~\ref{dimisnaturalenough} even if $f_\shriek$ would be well defined without it (after tensoring with $\BQ$).
\end{Rem}

\begin{Prf}[\Proofof{Lemma~\ref{dimisnaturalenough}}]
	The first part is obvious
	 so assume that $f\colon A\to B$ is locally finite and faithful. Let $[b]\in \pi_0B$ and $\rho\colon A\to \vect\BC$; we want to show the identity
	 \begin{equation}\label{pushforwardcomparison}
		 \dim_\BC\((\sV{f}\rho)(b)\) = f_\shriek(\cm_A(\rho))(b)\coloneqq \int\limits_{\twofib[f]Ab}[a]\mapsto\dim_\BC\rho(a)
	 \end{equation}
	 By the pointwise formula for Kan extensions, we can express the left side of Equation~\ref{pushforwardcomparison} in terms of a certain colimit over the $2$-fiber $\twofib[f]Ab$.
	 Since the $2$-fiber $\twofib[f]Ab$ only depends on the connected component $B(b)$ of $b$ in $B$ and the restriction of $f$ to the preimage of $B(b)$ we may assume that $B$ is connected.\\
	 Since $\rho$ is finitary we may further assume that it is nonzero only on one isomorphism class $[a]\in \pi_0A$. In this case both sides of Equation~\ref{pushforwardcomparison} only depend on the connected component of $a$ in $A$ so we may assume that $A$ is connected as well.
	If $A$ and $B$ are both connected, we may assume without loss of generality that $A=\B H$ and $B= \B G$.\\
	Now using that $f$ is faithful we reduce to the case where $f$ is just the inclusion of a supgroup $H\leqslant G$ (of finite index). In this case Diagram~\ref{pushforwardcompsquare} reads
	\[\begin{tikzcd}
					\K_0\lrep H\BC\ar[d,"\K_0\({\lreg G\BC}\otimes_{\BC[H]}\blank\)"']\ar[r,"\dim_\BC"]&\BZ\ar[d,"\cdot\multipl GH"]\\
					\K_0\lrep G\BC\ar[r,"\dim_\BC"]&\BZ
	\end{tikzcd}\]
	and commutes since $\lreg G\BC$ is a free $\BC[H]$-module of rank $\multipl GH$.
\end{Prf}

\begin{Prf}[\Proofof{Proposition~\ref{Propcomparisonmap}}]
We can compare the multiplications on the fat Hall algebra (\ref{explicitmonoidalproduct}) and the classical Hall algebra (Proposition~\ref{descriptionHallalg}) in the following diagram:
\begin{equation}\label{multiplicationfatclassical}
	\begin{tikzcd}
			\K_0\calV(\calX_\1)\times\K_0\calV(\calX_\1)\rar{\otimes }\ar[d,"\cm_{\calX_\1}\times\cm_{\calX_\1}"]&\K_0\calV(\calX_\1\times\calX_\1)\rar{\pV{c_2}}\ar[d,"\cm_{(\calX_\1\times\calX_\1)}"]&\K_0\calV(\calX_\2)\rar{\sV{\nu_2}}\ar[d,"\cm_{\calX_\2}"]&\K_0\calV(\calX_\1)\ar[d,"\cm_{\calX_\1}"]\\
		V(\calX_\1)\times V(\calX_\1)\rar{\cdot}& V(\calX_\1\times\calX_\1)\rar{c_2^\star}&V(\calX_\2)\rar{(\nu_2)\shriek}&V(\calX_\1)
	\end{tikzcd}
\end{equation}
The first square commutes because the dimension of a tensor product is the product of the dimensions and the other two squares commute by Lemma~\ref{dimisnaturalenough} (recall that the map $\calX_\2\xra{\nu} \calX_\1$ is faithful because $\calX$ is integral).\\
Therefore Diagram~\ref{multiplicationfatclassical} commutes, i.e.\ the dimension map $\cm\colon \fHall{\calX}\to\hall{\calX}$ preserves multiplication. Doing the same thing for the unit (using the description~\ref{explicitmonoidalunit}) we obtain that $\cm$ is a homomorphism of (unital) algebras.\\
This concludes Part~\ref{cmisalghomo} of the proposition. Part~\ref{cmismodhomo} is essentially the same using the explicit descriptions of the action maps in Theorem~\ref{descriptionHallmon} and Proposition~\ref{descriptionHallalg}.
\end{Prf}

\subsubsection{Sections of the dimension map}
We want to construct two extremal sections of the dimension map.
\begin{Def}
\begin{enumerate}
	\item The \introduce{trivial section} of the dimension map is the $\BZ$-linear map
	\[\cmtriv_\calX=\cmtriv_{\calX_\1} \text{ (resp.\ $\cmtriv_{\calX^\mindecor}=\cmtriv_{\calX^\mindecor_{\{\minimal,0\}}}$)}\colon \hall{\calX^{(\mindecor)}}\lra \fHall{\calX^{(\mindecor)}}\coloneqq \K_0\Hallf{\calX^{(\mindecor)}}\]
	given by mapping the basis element $\delta_{[a]}\in \hall\calX^{(\mindecor)}$ (for $[a]\in \pi_0\calX_\1$ resp.\ $[a]\in \pi_0\calX^\mindecor_{\{\minimal,0\}}$) to the class in $\K_0$ of the following finitary functor:
	\begin{equation}\label{reversecomparisonmaptrivial}
			\cmtriv(\delta_{[a]})\colon \calX_\1\xthra{}\pi_0\calX_\1\xra{\BC\Iso_{\pi_0\calX_\1}([a],\blank)}\vect\BC
	\end{equation}
	(resp.\ $\cmtriv(\delta_{[a]})\colon \calX^\mindecor_{\{\minimal,0\}}\xthra{}\pi_0\calX^\mindecor_{\{\minimal,0\}}\xra{\BC\Iso_{\pi_0\calX^\mindecor_{\{\minimal,0\}}}([a],\blank)}\vect\BC$)
	\item Assume $\calX^{(\mindecor)}$ is regular. The \introduce{regular section} of the dimension map is the $\BQ$-linear map
	\[\cmreg_\calX=\cmreg_{\calX_\1}\text{ (resp.\ $\cmreg_{\calX^\mindecor}=\cmreg_{\calX^\mindecor_{\{\minimal,0\}}}$)}\colon \hall[\BQ]{\calX^{(\mindecor)}}\lra \fHall[\BQ]{\calX^{(\mindecor)}}\coloneqq\BQ\otimes_\BZ \K_0\Hallf{\calX^{(\mindecor)}}\]
	given by mapping the basis element $\delta_{[a]}\in \hall\calX$ (for $[a]\in \pi_0\calX_\1$ resp.\ $[a]\in \pi_0\calX^\mindecor_{\{\minimal, 0\}}$) to the element
	\begin{equation}
			\cmreg(\delta_{[a]})\coloneqq \frac{1}{\#\Aut(a)}\left[\calX_\1\xra{\BC\Iso_{\calX_\1}(a,\blank)}\vect\BC\right].
	\end{equation}
(resp.\ $\cmreg(\delta_{[a]})\coloneqq \frac{1}{\#\Aut(a)}\left[\calX^\mindecor_{\{\minimal,0\}}\xra{\BC\Iso_{\calX^\mindecor_{\{\minimal,0\}}}(a,\blank)}\vect\BC\right]$)
	in the Grothendieck group.\qedhere
\end{enumerate}
\end{Def}
\begin{Rem}
	Here we view $\pi_0\calX_\1$ as a discrete category, hence Diagram~\ref{reversecomparisonmaptrivial} is just a fancy way of describing the functor which maps an arrow $x\xra{\cong} y$ in $\calX_\1$ to $\Id_\BC\colon \BC\to\BC$ if $x\cong a\cong y$ and all other arrows to zero.\\
	The functor $\BC\Iso(a,\blank)$ depends on the chosen representative $a\in [a]$ only up to isomorphism, hence $\cmreg(\delta_{[a]})$ is well defined in the Grothendieck group.\\
	Note that to define the regular section we need that the automorphism groups of $\calX_\1$ (resp.\ $\calX_{\{\minimal,0\}}^\mindecor$) are finite so that the denominator is finite and so that $\BC\Iso(a,\blank)$ does indeed take values in \emph{finite dimensional} vector spaces.
\end{Rem}

\begin{Prop}\label{sectionstocompmap}
	\begin{enumerate}
	\item\label{cmtrivissection} The map $\cmtriv$ is a $\BZ$-linear section of the dimension map $\cm\colon \fHall\calX\to\Hall\calX$. In particular the dimension map is always surjective.
	\item\label{cmregisquasisection} The map $\cmreg$ is a $\BQ$-linear section of the dimension map $\BQ\otimes\cm\colon \fHall[\BQ]\calX\to\hall[\BQ]\calX$
	\item\label{cmregisalghomo} Assume that the chopping map $c\colon \calX_\1\times\calX_\1\la\calX_\2$ (resp.\ $c\colon \calX^\mindecor_{\{\minimal,0\}}\times\calX_\1\la\calX_{\{\minimal,0,1\}}^\mindecor$) is faithful and $\calX_\0$ that is discrete. Then $\cmreg\colon \hall[\BQ]{\calX^{(\mindecor)}}\hra\fHall[\BQ]{\calX^{(\mindecor)}}$ is an algebra-(resp.\ module-)homomorphism.\qedhere
	\end{enumerate}
\end{Prop}

We call $\cmtriv$ and $\cmreg$ the \introduce{trivial section} and the \introduce{regular section} of the dimension map respectively.

\begin{Prf}[\Proofof{Proposition~\ref{sectionstocompmap}}]
Parts~\ref{cmtrivissection} and~\ref{cmregisquasisection} of Proposition~\ref{sectionstocompmap} are immediate from the definitions.\\
To prove that $\cmreg$ is an algebra- (resp.\ module-) homomorphism we proceed analogously to the proof for $\cm$: we define $\cmreg_{A}\colon \BQ\otimes_\BZ V(A)\to \BQ\otimes_\BZ \K_0[A,\vect\BC]_\finitary$ for every locally finite groupoid $A$ by the same formula $\cmreg(\delta_{[a]})\coloneqq \frac{1}{\#\Aut(a)}[\BC\Iso_A(a,\blank)]$. Then part~\ref{cmregisalghomo} of Proposition~\ref{sectionstocompmap} follows from the following Lemma which says that this collection is \roughly{natural enough} since all groupoids in the multiplication span (resp.\ action span) are locally finite by regularity.
\end{Prf}

\begin{Lem}\label{cmregisnaturalenough}
	Let $f\colon A\to B$ be a map of locally finite groupoids.
	\begin{enumerate}
		\item\label{regcomparisonpullnatural} If $f$ is $\pi_0$-finite and faithful then we have the following commutative diagram of abelian groups
			\begin{equation}\label{regpullbackcompsquare}
				\begin{tikzcd}
					\BQ\otimes_\BZ\K_0[A,\vectC]_\finitary\ar[from=r,"\cmreg_A"']&\BQ\otimes_\BZ V(A)\\
					\BQ\otimes_\BZ\K_0[B,\vectC]_\finitary\ar[from=r,"\cmreg_B"']\ar[u,"\pV f"]&\BQ\otimes_\BZ V(B)\ar[u,"f^\star"]
				\end{tikzcd}
			\end{equation}
		\item\label{regcomparisonpushnatural} If $f$ is locally finite (we do \emph{not} assume faithful) then we have the following commutative diagram of abelian groups
			\begin{equation}\label{regpushforwardcompsquare}
				\begin{tikzcd}
					\BQ\otimes_\BZ \K_0[A,\vectC]_\finitary\ar[d,"\sV f"]\ar[from=r,"\cmreg_A"']&\BQ\otimes_\BZ V(A)\ar[d,"f_\shriek"]\\
					\BQ\otimes_\BZ \K_0[B,\vectC]_\finitary\ar[from=r,"\cmreg_B"']&\BQ\otimes_\BZ V(B)
				\end{tikzcd}\qedhere
			\end{equation}
	\end{enumerate}
\end{Lem}
\begin{Prf}
	By making similar reduction steps as in the proof of Lemma~\ref{dimisnaturalenough} (and using that $f$ is faithful) we may assume that $f$ is just the inclusion $\B H\hra \B G$ of a subgroup $H\leqslant G$.
	In this case Diagram~\ref{regpullbackcompsquare} reads
		\[\begin{tikzcd}[column sep=100]
			\BQ\otimes_\BZ\K_0\lrep H\BC\ar[from=r,"{\frac{1}{\# H}\lreg H\BC}\longmapsfrom 1"']&\BQ\\
			\BQ\otimes_\BZ\K_0\lrep G\BC\ar[from=r,"{\frac{1}{\# G}\lreg G\BC}\longmapsfrom 1"']\ar[u,"\BQ\otimes_\BZ\K_0\Res^G_H"]&\BQ\ar[u,"="]
		\end{tikzcd}\]
		and commutes since $\Res^G_H\lreg G\BC\cong\lreg H\BC^{\oplus \multipl GH}$.\\
		In the second case we reduce to a group homomorphism $f\colon \B H\to \B G$ (which is not necessarily faithful); Diagram~\ref{regpushforwardcompsquare} then reads
				\[\begin{tikzcd}[column sep=100]
			\BQ\otimes_\BZ\K_0\lrep H\BC\ar[from=r,"{\frac{1}{\# H}\lreg H\BC}\longmapsfrom 1"']&\BQ\\
			\BQ\otimes_\BZ\K_0\lrep G\BC\ar[from=r,"{\frac{1}{\# G}\lreg G\BC}\longmapsfrom 1"']\ar[from=u,"\BQ\otimes_\BZ\K_0\({\lreg G\BC}\otimes_{\BC[H]}\blank\)"']&\BQ\ar[from=u,"\frac{\#G}{\# H}"]
		\end{tikzcd}\]
		which commutes since $\BC[G]\otimes_{\BC[H]}\lreg H\BC=\lreg G\BC$ and $\frac{\#G}{\# H}\frac{1}{\#G}=\frac{1}{\#H}$.
\end{Prf}

\ifx\isstandalone\undefined
\else
\newpage
\bibliographystyle{amsalpha}
\bibliography{literatur}
\end{document}
\fi

\subsection{Idempotent decomposition of the Hall monoidal category}
\ifx\preambleloaded\undefined
   
   \def\isstandaloneODLE{}
\fi

Let us denote by $J_T$ the unit $\terminal\to\calV(\terminal)\xra{\pV{}}\calV(T)$ of the (cococo) monoidal left derivator $\calV$.
Let $\calX\colon \Dop\to\Grpd$ be a $2$-Segal simplicial groupoid and let $\HallV\calX$ be its associated Hall monoidal category.\\
Theorem~\ref{descriptionHallmon} tells us how to compute the monoidal unit $I\in\HallV\calX$:\\
We take the monoidal unit $J_{\calX_\0}\colon \terminal\to \calV(\calX_\0)$ of $\calV$ (which in the case $\calV=\Vect\BC$ is explicitly given as the constant functor $\calX_\0\to\Vect\BC$ with value $\BC$) and push it forward to $\calV(\calX_\1)$ via the degeneracy map $\nu_0\colon \calX_\0\to\calX_\1$.\\

We can refine this description by constructing some orthogonal idempotents which add up to the monoidal unit.
\begin{Prop}\label{idempotentsinHallmoncat}
We have a decomposition
\[I\cong \coprod_{T\in\pi_0\calX_0} I_T\]
of the monoidal unit $I\in\Hall\calX$ into pairwise orthogonal idempotents (this means that $I_T\boxtimes I_{T'}\cong\initial$ if $T\neq T'$ and $I_T\boxtimes I_T\cong I_T$) defined by the compositions
\[I_T\colon \terminal\xra{J_T}\calV(T)\xra{\sV t}\calV(\calX_\0)\xra{\sV{\nu_0}}\calV(\calX_\1)\]
for each connected component $t\colon T\hra \calX_0$ of $\calX_0$. 
\end{Prop}
\begin{Rem}
Note that we do not need any finiteness assumptions since in $\calV(\calX_\1)$ we can take arbitrary coproducts.\\
In the finitary case the groupoid $\calX_\0$ is $\pi_0$-finite, hence we have a \emph{finite} set of idempotents $I_T$ which add up to the identity, so in the fat Hall algebra $\K_0\Hallf\calX$ we get a well defined decomposition $1=\sum_{T\in\pi_0\calX_0}[I_T]$ of the unit into pairwise orthogonal idempotents. If $\calX$ is integral then we can use the dimension map to push this decomposition down to the Hall algebra $\hall\calX$.
\end{Rem}

\begin{Prf}[\Proofof{Proposition~\ref{idempotentsinHallmoncat}}]
Let $t\colon T\hra \calX_0$ and $t'\colon T'\hra\calX_0$ be inclusions of two connected components. We define $(P,p,\Delta')$ via the (strict) pullback diagram
\begin{equation}\begin{tikzcd}\label{sdklafjlhgurweasdf}
	P\ar[r,"\Delta'"]\ar[d,"p"]&T\times T'\ar[d,"t\times t'"]\\
	\calX_0\ar[r,"\Delta"]&\calX_0\times\calX_0
	\end{tikzcd}\end{equation}
where $\Delta\colon \calX_0\to\calX_0\times \calX_0$ is the diagonal.\\
For $T=T'$ the morphism $\Delta'\colon P\to T\times T'$ is just the diagonal $T\to T\times T$ (and in this case $p=t$). If $T$ and $T'$ are different (hence disjoint) then $\Delta'$ is the empty morphism $\initial\to T\times T'$.\\
Observe that $(t\times t')\colon T\times T'\hra\calX_0$ is an iso-fibration, hence Diagram~\ref{sdklafjlhgurweasdf} is also a $2$-pullback square.
By applying $\calV$ to the whole setting including the multiplication span we obtain the diagram (where $J_{T,T'}$ is defined by the diagram)
\[\begin{tikzcd}
	\terminal\ar[ddd,"{(I_T,I_{t'})}"', bend right=75]\ar[d,"{(J_T,J_T')}"]\ar[rrd,"J_{T,T'}", dashed,bend left=15]\\
	\calV(T)\times\calV(T')\ar[d,"\sV t\times \sV {t'}"']\ar[r,"\fotV2"]&	\calV(T\times T')\ar[r,"\pV {\Delta'}"]\ar[d,"\sV{t\times t'}"]& \calV(P)\ar[d,"\sV p"]\\
	\calV(\calX_\0)\times\calV(\calX_\0)\ar[d,"\sV{}\times \sV{}"']\ar[r,"\fotV2"]\ar[ru,Leftrightarrow]&	\calV(\calX_\0\times \calX_\0)\ar[r,"\pV \Delta"]\ar[d,"\sV{}"]\ar[ru,Leftrightarrow]& \calV(\calX_\0)\ar[r,equal]\ar[d,"\sV{}"]& \calV(\calX_\0)\ar[d,"\sV{}"]\\
	\calV(\calX_\1)\times\calV(\calX_\1)\ar[rrr,"\boxtimes_2",bend right=15]\ar[r,"\fotV2"]\ar[ru,Leftrightarrow]&	\calV(\calX_\1\times \calX_\1)\ar[r,"\pV{c_2}"]\ar[ru,Leftrightarrow]& \calV(\calX_\2)\ar[r,"\sV{\nu_2}"]\ar[ru,Leftrightarrow]& \calV(\calX_\1)\\	
\end{tikzcd}\]
of categories and functors. The mates in the first column are isomorphisms because $\calV$ is componentwise cocontinuous. The first mate in the second column is an isomorphism since Diagram~\ref{sdklafjlhgurweasdf} is $2$-pullback. The second mate in the second column is an isomorphisms because the square
\[\begin{tikzcd}
	\xX{0}\dar\rar&\xX{0}\times \xX{0}\dar\\
	\xX{0,1,2}\rar&\xX{0,1}\times \xX{1,2}
\end{tikzcd}\]
is an instance of the $2$-Segal square~\ref{generalunitalsubsimplicial} (for the map $\{0\}\thla\{0,1,2\}$), hence is $2$-pullback. The last mate is always an isomorphism.\\
The composition along the left and bottom boundary always gives $I_T\boxtimes I_T'$
\begin{itemize}
\item If $T=T'$ then
\[J_{T,T}\coloneqq \pV{\Delta}\circ \fotV2(J_T,J_T)\cong J_T\otimes J_T\cong J_T\]
(see Remark~\ref{multitensorwithdiagonal} for the first isomorphism, for the second we use that $J_T$ is the monoidal unit), hence the composition along the upper right gives $I_T$. Hence we obtain the desired isomorphism $I_T\boxtimes I_T\cong I_T$.\\
\item If $T\neq T'$ then $P=\initial$ implies that the composition along the top right picks out the initial objects $\initial$ of $\calV(\calX_\1)$. Hence we get $I_T\boxtimes I_T\cong\initial$ as desired.
\end{itemize}
It is straightforward to check that the coproduct of the $I_T$'s (in the cocomplete category $\calV(\calX_\1)$) does indeed give the unit $I$; the proof is thus concluded.
\end{Prf}

\ifx\isstandaloneODLE\undefined
\else
\newpage
\bibliographystyle{amsalpha}
\bibliography{literatur}
\end{document}
\fi

\subsection{Functoriality}
\ifx\preambleloaded\undefined
   
   \def\isstandaloneBMDJSE{}
\fi

We have seen in Section~\ref{sectionHallfunctoriality} that the assignment
\[\HallV{}\colon \{\text{$2$-Segal simpl. $\Dop\to\Grpd$}\}\lra \{\text{monoidal categories}\}^\lax\]
is actually a contravariant functor. We can describe this functor a little bit more precisely:\\
Every morphism $\alpha\colon \calY\to\calX$ of $2$-Segal simplicial groupoids gives rise to a lax monoidal functor $\alpha^\star=\HallV\alpha\colon \HallV\calX\to\HallV\calY$ which is given on underlying categories by pullback $\pV{\alpha_\1}\colon \calV(\calX_\1)\to\calV(\calY_\1)$.\\
	The transformation $\alpha^\star(\blank)\boxtimes_\calY\alpha^\star(\blank)\lRa \alpha^\star(\blank \boxtimes_\calX\blank)$ is given by the pasting
		\begin{equation}		
			\begin{tikzcd}
			\calV(\calX_\1)\times\calV(\calX_\1)\rar{\fotV{2}}\ar[d,"\pV\alpha\times\pV\alpha"]&\calV(\calX_\1\times\calX_\1)\rar{\pV{c_2}}\ar[d,"\pV{\alpha\times\alpha}"]&\calV(\calX_\2)\rar{\sV{\nu_2}}\ar[d,"\pV\alpha"]&\calV(\calX_\1)\ar[d,"\pV\alpha"]\\
			\calV(\calY_\1)\times\calV(\calY_\1)\rar{\fotV{2}}&\calV(\calY_\1\times\calY_\1)\rar{\pV{c_2}}&\calV(\calY_\2)\rar{\sV{\nu_2}}\ar[ru,Rightarrow]&\calV(\calY_\1)\\
			\end{tikzcd}
		\end{equation}
	where the first two squares commute (up to canonical isomorphisms) and the third is inhabited by a canonical pull-push mate.\\

	In a very similar spirit the assignment $\calX\mapsto\HallV\calX$ can also be made into a covariant functor
	\[\HallV{}\colon \{\text{$2$-Segal simpl. $\Dop\to\Grpd$}\}\lra \{\text{monoidal categories}\}^{\op\lax}\]
	where on the right we consider oplax monoidal functors between monoidal categories:\\
	Every morphism $\alpha\colon \calX\to\calY$ of $2$-Segal simplicial groupoids gives rise to an oplax monoidal functor $\alpha_\shriek\colon \HallV\calX\to\HallV\calY$ which is given on underlying categories by left Kan extension $\sV{\alpha_\1}\colon \calV(\calX_\1)\to\calV(\calY_\1)$.\\
	The transformation $\alpha_\shriek(\blank \boxtimes_\calX\blank)\lRa\alpha_\shriek(\blank)\boxtimes_\calY\alpha_\shriek(\blank)$ is given by the pasting
		\begin{equation}		
			\begin{tikzcd}
			\calV(\calX_\1)\times\calV(\calX_\1)\rar{\fotV{2}}\ar[d,"\sV\alpha\times\sV\alpha"]&\calV(\calX_\1\times\calX_\1)\rar{\pV{c_2}}\ar[d,"\sV{\alpha\times\alpha}"]&\calV(\calX_\2)\rar{\sV{\nu_2}}\ar[d,"\sV\alpha"]&\calV(\calX_\1)\ar[d,"\sV\alpha"]\\
			\calV(\calY_\1)\times\calV(\calY_\1)\rar{\fotV{2}}&\calV(\calY_\1\times\calY_\1)\rar{\pV{c_2}}\ar[ru,Leftarrow]&\calV(\calY_\2)\rar{\sV{\nu_2}}\ar[ru,Leftrightarrow]&\calV(\calY_\1)\\
			\end{tikzcd}
		\end{equation}
	where the outer two squares commute (up to canonical isomorphisms) and the middle one is inhabited by a canonical pull-push mate.\\
	We omit a detailed proof of all the coherence properties for the oplax functor $\alpha_\shriek$.
	
\begin{War}
Both pullback $\alpha^\star$ and left Kan extension $\alpha_\shriek$ commute with coproducts, so in the finitary case they induce a $\BZ$-module homomorphism on the Grothendieck groups $\K_0\HallV{?}$. However this maps $\K_0(\alpha^\star)$ and $\K_0(\alpha_\shriek)$ are \emph{not} algebra homomorphisms in general, because we usually do not have isomorphisms of functors $\alpha^\star(\blank)\boxtimes_\calY\alpha^\star(\blank)\lLRa \alpha^\star(\blank \boxtimes_\calX\blank)$ or $\alpha_\shriek(\blank \boxtimes_\calX\blank)\lLRa\alpha_\shriek(\blank)\boxtimes_\calY\alpha_\shriek(\blank)$.
\end{War}

\ifx\isstandaloneBMDJSE\undefined
\else
\newpage
\bibliographystyle{amsalpha}
\bibliography{literatur}
\end{document}
\fi

\subsection{Hall (categorical) modules via bounded flags}\label{sectionboundedflags}
\ifx\preambleloaded\undefined
   
   \def\isstandalone{}
\fi

Let $\calA$ be a finitary proto-abelian category.\\
In Section~\ref{defgroupoidofflags} we defined the simplicial groupoid $\calX^\calA\colon \Dop\to\Grpd$ of flags in $\calA$; this allowed us to define the Hall algebra and the Hall monoidal category of $\calA$ by applying the generalized Hall construction to $\calX^\calA$.\\
Now we want to make use of our general machinery to construct some (categorical) left modules of the Hall algebra (resp.\ Hall monoidal category). To do this we must associate to $\calA$ some left relative simplicial groupoids $\Dmop\to\Grpd$. These groupoids will again consist of flags in $\calA$ but this time we impose some boundary conditions. 

\subsubsection{Quotient data}
We are looking for sensible boundary conditions to impose on flags of the shape
\[0\hra A_1\hra A_2\hra\dots\hra A_{n}\hra A_\maximal.\]
The most obvious idea is to ask for the rightmost object $A_\maximal$ to belong to some subcategory $\calQ\subseteq \calA$. This might work for a fixed $n$, but when considering transition maps $\n\la \m$ it becomes apparent that we need to put requirements on all the quotients $A_{i,\maximal}\coloneqq \Cok(A_i\hra A_\maximal)$. In order for this to work we need to impose some axioms.

\begin{Def}
We call a subcategory $0\in\calQ\subseteq\calA$ a \introduce{quotient datum} in $\calA$ if the following conditions hold
\begin{enumerate}[label=\axiomlabelstyle{Q},ref=\axiomrefstyle{Q}]
\item\label{boundaryareallepi} Every morphism in $\calQ$ is an epimorphism.
\item\label{boundarycompleteses} For every monomorphism $A\hra H$ (in $\calA$) with $H\in\calQ$ there is an epimorphism $H\thra Q$ in $\calQ$ completing a short exact sequence $0\to A\hra H\thra Q\to0$.
\item\label{boundaryuniqueses} If $H\thra Q\xra{\cong}Q'$ is a diagram in $\calA$ such that both $H\thra Q$ and the composition $H\thra Q'$ are in $\calQ$ then the isomorphism $Q\xra{\cong}Q'$ is also in $\calQ$.\qedhere
\end{enumerate}
\end{Def}

\subsubsection{The quotient datum of a group of automorphisms}\label{sectionsubgroupquotientdatum}

Let $X$ be an object in $\calA$ and let $S\subseteq \Aut(X)$ be a subgroup of automorphism of $X$.
We can define the \introduce{abstract quotient datum} associated to $S$ as the category $\Quot_S(X)$ of quotients of $X$ up to $S$:
\begin{itemize}
\item objects are epimorphisms $q\colon X\thra Y$
\item morphisms in $\Quot_S(X)$ from $X\thra Y$ to $X\thra Y'$ are  maps $Y\to Y'$ in $\calA$ such that the square
\[
\begin{tikzcd}
	X\ar[d,"s"',"\cong"]\ar[r,twoheadrightarrow]&Y\ar[d]\\
	X\ar[r,twoheadrightarrow]&Y'
\end{tikzcd}
\]
commutes for some automorphism $s\in S$ (such maps $Y\thra Y'$ are automatically epimorphisms themselves).
\end{itemize}

\begin{Expl}
If $\triv=S\subset\Aut(X)$ is the trivial group then $\Quot(X)=\Quot_\triv(X)$ is equivalent to the partially ordered set of quotient objects of $X$.
\end{Expl}

There is a canonical faithful functor $\Quot_S(X)\to \calA$ which just remembers the target $Y$ of an epimorphism $X\thra Y$. We would like to think of $\Quot_S(X)$ as a subcategory of $\calA$ via this functor. Alas this is not directly possible because $\Quot_S(X)\to \calA$ is not injective on objects; luckily there is a cheap fix:\\
We replace $\calA$ by the bloated category $\hat\calA\supset \calA$ which has additional copies $Y_q$ of each object $Y\in\calA$ indexed by epimorphisms $q\colon X\thra Y$. Clearly $\hat\calA$ is equivalent to $\calA$ via $Y_q\mapsto Y$ and we can factor the forgetful functor $\Quot_S(X)\to \calA$ as
\[\Quot_S(X)\xra{q\mapsto Y_q}\hat\calA\xra{\simeq}\calA\]
Now the functor $\Quot_S(X)\hra\hat\calA$ does in fact define a honest subcategory of $\hat\calA$ which we denote by $\calQ_S(X)$ and call the \introduce{quotient datum associated to $S$ (and $X$)}.

\begin{Lem}\label{automorphgroupgivesquotientdatum}The subcategory $\Quot_S(X)\cong \calQ_S(X)\subseteq\hat\calA$ is a quotient datum on $\hat\calA$.
\end{Lem}
\begin{Prf} Condition~\ref{boundaryareallepi} is obvious.
\begin{itemize}
\item[\ref{boundarycompleteses}] The object $H$ (of $\hat\calA$) being in $\calQ_S(X)$ means that there is an epi $q\colon X\thra Y$ in $\calA$ such that $H=Y_q$. We can complete the monomorphism $A\hra Y_q\cong Y$ to a short exact sequence with $Q\in\calA$
\[0\to A\hra Y \xthra{p} Q\to 0.\]
Then
\[0\to A\hra Y_q\xthra{p} Q_{pq}\to 0\]
is the desired short exact sequence where $Y_q\xthra{p}Q_{pq}$ comes from the morphism
\[\begin{tikzcd}[cramped]
	X\ar[rd,"pq"',twoheadrightarrow]\ar[r,"q",twoheadrightarrow]&Y\ar[d,"p",dashed,twoheadrightarrow]\\
	&Q
\end{tikzcd}\]
in $\Quot_S(X)$.
\item[\ref{boundaryuniqueses}] Consider $H\xthra{\alpha} Q\xra{\beta} Q'$ in $\hat\calA$ as in Condition~\ref{boundaryuniqueses}. The three objects $H$, $Q$ and $Q'$ are in $\calQ_S(X)$, hence we can write them as $Y_p$, $Z_q$ and $Z'_{q'}$ for some epimorphisms $p\colon X\thra Y$, $q^{(')}\colon X\thra Z^{(')}$ (in $\calA$) respectively.
Thus we have the situation in $\calA$
\[\begin{tikzcd}
	X\ar[rr,"s'", bend left]\ar[d,"p",twoheadrightarrow]\ar[r,"s"']&X\ar[r,"s's^\inv"',dashed]\ar[d,"q",twoheadrightarrow]&X\ar[d,"q'",twoheadrightarrow]\\
	Y\ar[r,"\alpha",twoheadrightarrow]&Z\ar[r,"\beta"]&Z'	
\end{tikzcd}\]
where $s,s'\in S$ are chosen such that the left square and the total rectangle commute; this is possible since $Y\thra Z$ and $Y\thra Z'$ are morphisms of $\Quot_S(X)$ by assumption. We conclude that the rightmost square commutes with $s's^\inv\in S$, hence $\beta$ is a morphism of $\Quot_S(X)$ (note that we did not use the assumption that $\beta$ is an isomorphism).\qedhere
\end{itemize}\end{Prf}

\begin{Rem}If the category $\calA$ is already bloated enough then we can choose representatives $q\colon X\thra Y_q$ for each quotient object $[q]$ of $X$ such that the $Y_q$ are pairwise non-equal objects. In this case we can replace $\Quot_S(X)$ by the dense and full (hence equivalent) subcategory $\calQ_S(X)$ spanned by the objects $q\colon X\thra Y_q$. On the other hand the restricted functor $\calQ_S(X)\hra\calA$ is now injective on objects, hence we can view it as a quotient datum on $\calA$.
\end{Rem}

\subsubsection{The $2$-Segal relative simplicial groupoid of bounded flags}
As indicated above, we can construct groupoids out of a quotient datum.
\begin{Def}
Let $\calQ$ be a quotient datum on $\calA$. We define the left relative simplicial groupoid $\XAQ\colon \Dmop\to\Grpd$ of \introduce{flags in $\calA$ bounded by $\calQ$} as follows:
\begin{itemize}
	\item Let $N^n=\{x_0<\dots< x_n\}\in \Dmop$ be an object. We define $\XAQ_N$ to be the subgroupoid of $\calX^\calA_N$ where in the case $x_n=\maximal$ we require the rightmost column of the diagram
	\begin{equation}\label{boundedflag}
	\begin{tikzcd}
				A_{x_0x_1}\ar[r, hookrightarrow]\ar[rd, "{\square}"description, phantom]\ar[d, twoheadrightarrow]&	A_{x_0x_2}\ar[r, hookrightarrow]\ar[d, twoheadrightarrow]&	\cdots\ar[r, hookrightarrow]&	A_{x_0x_{n-2}}\ar[rd, "{\square}"description, phantom]\ar[d, twoheadrightarrow]\ar[r, hookrightarrow]&	A_{x_0x_{n-1}}\ar[rd, "{\square}"description, phantom]\ar[d, twoheadrightarrow]\ar[r, hookrightarrow]	&A_{x_0x_n}\ar[d, twoheadrightarrow]\\
				0\ar[r, hookrightarrow]&	A_{x_1x_2}\ar[d, twoheadrightarrow]\ar[r, hookrightarrow]&	\cdots\ar[r, hookrightarrow]&	A_{x_1x_{n-2}}\ar[rd, "{\square}"description, phantom]\ar[d, twoheadrightarrow]\ar[r, hookrightarrow]&	A_{x_1x_{n-1}}\ar[rd, "{\square}"description, phantom]\ar[d, twoheadrightarrow]\ar[r, hookrightarrow]&	A_{x_1x_n}\ar[d, twoheadrightarrow]\ar[d, twoheadrightarrow]\\ 
			&	0\ar[r, hookrightarrow]&	\cdots\ar[r, hookrightarrow]&	A_{x_2x_{n-2}}\ar[d, twoheadrightarrow]\ar[r, hookrightarrow]&	A_{x_2x_{n-1}}\ar[d, twoheadrightarrow]\ar[r, hookrightarrow]&	A_{x_2x_n}\ar[d, twoheadrightarrow]\ar[d, twoheadrightarrow]\\ 
			&&\ddots&\vdots\ar[d, twoheadrightarrow]&\vdots\ar[d, twoheadrightarrow]&\vdots\ar[d, twoheadrightarrow]\ar[d, twoheadrightarrow]\\
			&&&0\ar[r, hookrightarrow]&A_{x_{n-2} x_{n-1}}\ar[rd, "{\square}"description, phantom]\ar[r, hookrightarrow]\ar[d, twoheadrightarrow]&A_{x_{n-2}x_n}\ar[d, twoheadrightarrow]		\\	
			&&&&0\ar[r, hookrightarrow]&A_{x_{n-1}x_n}
		\end{tikzcd}
	\end{equation}
	to be contained in $\calQ$ (both for the vertical arrows $A_{x_i,\maximal}\thra A_{x_{i+1},\maximal}$ and for isomorphisms $A_{x_i,\maximal}\xra{\cong} A'_{x_i,\maximal}$ between two such diagrams).
	\item The morphisms are induced by those of $\calX^\calA$, i.e.\ they are given by simultaneously duplicating (in case of degeneracy maps) or deleting (in the case of face maps) rows and columns.\\
	Note that the condition on the morphisms in $\Delta_{\maxdecor}$ ($x\mapsto\maximal$ if and only if $x=\maximal$) guarantees that the condition on the rightmost column being in $\calQ$ (in the case $x_n=\maximal$) is preserved.\qedhere
\end{itemize}
\end{Def}

\begin{Prop}\label{boundedflags2Segal}
For a quotient datum $\calQ\subseteq \calA$ on $\calA$ the relative simplicial groupoid $\XAQ\colon \Dmop\to\Grpd$ is $2$-Segal.
\end{Prop}
\begin{Prf}
We need to show that Diagram~\ref{sub2segalcriterium} and Diagram~\ref{sub2unitalcriterium} of Proposition~\ref{sub2segalcrit} are $2$-pullback squares for $\XAQ$; basically all we need to do is to copy the proof from the case of unbounded flags~\cite[Theorem 2.10]{DyckerhoffHigherCategoricalAspectsOfHallAlgebras}:\\
It is easy to see that the maps $\XAQb{x_i,\dots,x_j}\to \XAQb{x_i,x_j}$ and $\XAQb{x_0,\dots,x_k,x_k',\dots,n}\to\XAQb{x_k,x_k'}$ are iso-fibrations, hence it is enough to show that those diagrams are (strict) pullbacks.
\begin{itemize}
\item The map
	\begin{equation}\label{pullbackmapboundedflags}
		\Psi\colon \XAQb{x_0,\dots,x_n}\lra\XAQb{x_0,\dots,x_i,x_j,\dots,x_n}\times_{\XAQb{x_i,x_j}}\XAQb{x_i,\dots,x_j}
	\end{equation}
just forgets all the components $A_{x_ax_b}$ (in Diagram~\ref{boundedflag}) with $0\leqslant a<i<b<j\leqslant n$ or $0\leqslant i<a<j<b\leqslant n$.\\
Observe that those forgotten components can be recovered by completing the two short exact sequences
\[
\begin{tikzcd}
	0\ar[r]&A_{x_ax_b}\ar[r,hookrightarrow, dashed]&A_{x_ax_j}\ar[r,twoheadrightarrow]&A_{x_bx_j}\ar[r]&0\\
		0\ar[r]&A_{x_ix_a}\ar[r,hookrightarrow]&A_{x_ix_b}\ar[r,dashed,twoheadrightarrow]&A_{x_ax_b}\ar[r]&0
\end{tikzcd}
\]
In the case $x_b=x_n=\maximal$ we have to complete the second sequence in such a way that the dashed morphism is in $\calQ$; this is possible by Condition~\ref{boundarycompleteses}.\\
An isomorphism on the forgotten component can be recovered (uniquely!) by filling in the unique dotted vertical morphism
\[
\begin{tikzcd}
		0\ar[r]&A_{x_ix_a}\ar[d,"\cong"]\ar[r,hookrightarrow]&A_{x_ix_b}\ar[d,"\cong"]\ar[r,dashed,twoheadrightarrow]&A_{x_ax_b}\ar[r]\ar[d,"\cong",dotted]&0\\
				0\ar[r]&A'_{x_ix_a}\ar[r,hookrightarrow]&A'_{x_ix_b}\ar[r,dashed,twoheadrightarrow]&A'_{x_ax_b}\ar[r]&0
\end{tikzcd}
\]
(similarly for the other case). A priori, this dotted morphism only exists in $\calA$; in the case $x_b=x_n=\maximal$, however, it is automatically in $\calQ$ by Condition~\ref{boundaryuniqueses}.\\
Therefore the functor $\Psi$ is actually an equivalence of groupoids.
\item  The functor 
		\[\XAQb{x_0,\dots,x_n}\lra\XAQb{x_0,\dots,x_k,x_k',\dots,x_n}\]
is fully faithful with image consisting of flags
\[A_{x_0x_1}\hra\dots\hra A_{x_0x_k}\xra{=}A_{x_0x_k}\hra\dots\hra A_{x_0x_n}.\]
On the other hand $\XAQb{x_k}$ is just the trivial groupoid $\{0\}$, hence the fiber product 
\[\XAQb{x_0,\dots,x_k,x_k',\dots,x_n}\times_{\XAQb{x_k,x_k'}}\XAQb{x_k}\]
consists of those flags
		\[\begin{tikzcd}
		A_{x_0x_1}\ar[r,hookrightarrow]&\dots\ar[r,hookrightarrow]&A_{x_0x_k}\ar[d]\ar[rd, "{\square}"description, phantom]\ar[r,hookrightarrow]&A_{x_0x_k}\ar[d,twoheadrightarrow]\ar[r,hookrightarrow]&\dots\ar[r,hookrightarrow]&A_{x_0x_n}\\
		&&0\ar[r]&A_{x_kx_k'}&&
		\end{tikzcd}\]
		where $A_{x_kx_k'}\cong0$.
It follows immediately that the map
	\begin{equation}\label{pullbackmapboundedflagsunit}
		\XAQb{x_0,\dots,x_n}\lra\XAQb{x_0,\dots,x_k,x_k',\dots,x_n}\times_{\XAQb{x_k,x_k'}}\XAQb{x_k}
	\end{equation}
	is an equivalence of categories as required.
	Observe that the quotient datum $\calQ$ never enters the picture, because degeneracies cannot involve the element $\maximal$.
\qedhere
\end{itemize}
\end{Prf}

\begin{Rem}
Matthew Young pointed out that one can generalize the notion of a quotient datum by only requiring condition~\ref{boundarycompleteses} to hold for all monomorphisms $A\hra H$ with $A$ in a fixed Serre subcategory $\calB\subseteq \calA$. The corresponding left relative simplicial groupoid of flags $\calX^{\calA,\calB,\calQ}\colon\Dmop\to\Grpd$ then consists of those flags~\ref{boundedflag} in $\calA$ where the last column is in $\calQ$ and everything but the last column is in $\calB$. The same argument as for Prop~\ref{boundedflags2Segal} shows that this relative simplicial groupoid is $2$-Segal.\\
Young essentially proves~\cite[Theorem 2.7]{YoungRelative2SegalSpaces} that a stability function and a framing $\Phi$ on an abelian category $\calA$ give rise to such a generalized quotient datum on the Artin mapping cylinder $\calA_\Phi$; in this way he obtains a relative $2$-Segal simplicial groupoid which he calls the \buzzword{stable framed S-construction}.
\end{Rem}

\ifx\isstandalone\undefined
\else
\newpage
\bibliographystyle{amsalpha}
\bibliography{literatur}
\end{document}
\fi

\subsection{Duality pairing on Hall (categorical) modules}\label{sectiondualitypairing}
\ifx\preambleloaded\undefined
   
   \def\isstandaloneBJIFE{}
\fi

Let $\Phi\colon \calY\to\calX$ be a $2$-Segal morphism of simplicial groupoids such that both $\Phi^\mindecor$ and $\Phi^\maxdecor$ are regular. The way we constructed the Hall (left- or right-) module $\hall\Phi$ with underlying free $\BZ$-module $V(\calY_\0)\coloneqq \BZ^{(\pi_0\calY_\0)}$ it comes equipped with a canonical non-degenerate pairing
\[\apair{}{}'\colon \BZ^{(\pi_0\calY_\0)}\times \BZ^{(\pi_0\calY_\0)}\xra{}\BZ\]
given by 
\[\apair{\varphi_1}{\varphi_2}'\coloneqq \sum\limits_{[a]\in\pi_0\calY_\0}{\varphi_1(a)\varphi_2(a)}.\]
The questions are:
\begin{itemize}
\item How does this pairing interact with the left and right $\hall\calX$-module structures?
\item Is there an analogous statement for categorical Hall modules?
\end{itemize}

It turns out that we should renormalize the pairing $\apair{}{}'$ by keeping track of the automorphism groups of $\calY_\0$. We thus define the slightly different pairing $\apair{}{}$ by
\begin{equation}\label{Qvaluedpairingonhallmodule}
	\apair{\varphi_1}{\varphi_2}\coloneqq \sum\limits_{[a]\in\pi_0\calY_\0}\frac{\varphi_1(a)\varphi_2(a)}{\#\Aut(a)},
\end{equation}
which can be described in terms of groupoids as the pull-push composition
\begin{equation}
	\apair{}{}\colon V(\calY_\0)\times V(\calY_\0)\xra{\cdot}V(\calY_\0)\times V(\calY_\0)\xra{\Delta^\star}V(\calY_\0)\xra{\pi_\shriek}\BQ\otimes_\BZ V(\terminal)\cong\BQ.
\end{equation}

\begin{Rem}
Since the automorphism groups in $\calY_\0$ are non-trivial in most examples, we really need to define $\pair{}{}$ as a $\BQ$-valued pairing so that Equation~\ref{Qvaluedpairingonhallmodule} makes sense.
\end{Rem}

It turns out that this renormalized pairing is compatible with the left and right $\hall\calX$-module structure:

\begin{Prop}\label{Halldualitypairingbalanced}
The non-degenerate pairing
\begin{equation}
	\apair{}{}\colon \hall{\Phi^\mindecor}\times\hall{\Phi^\maxdecor}\to\BQ
\end{equation}
is $\hall\calX$-balanced, i.e.\ we have $\pair{\eta^\mindecor\actiondot\varphi}{\eta^\maxdecor}=\pair{\eta^\mindecor}{\varphi\actiondot\eta^\maxdecor}$ for all $\varphi\in\hall\calX$, $\eta^\mindecor\in\hall{\Phi^\mindecor}$ and $\eta^\maxdecor\in\hall{\Phi^\maxdecor}$
\end{Prop}

We could prove Proposition~\ref{Halldualitypairingbalanced} by hand but we will not do this. Instead we will see a more conceptual explanation and prove an analogous statement for categorical $\calV$-Hall modules.

\subsubsection{Bi-relative simplicial groupoids}
The main trick is to define a new admissible subcategory $\Delta_\mmdecor\subset \Delta$ which allows us to encode the additional duality pairing.\\
We pick distinct formal symbols $\minimal$ and $\maximal$.
\begin{Def}
We define the \introduce{bi-relative simplex category} $\Delta_\mmdecor\subset \Delta$ as the intersection $\Delta_\mmdecor\coloneqq \Delta_\mindecor\cap\Delta_\maxdecor$. Explicitly:
\begin{itemize}
\item Objects are those $N\in\Delta$ such that
	\begin{itemize}[label=\itemsecondlayer]
	\item $\{\minimal\}\neq N \neq\{\maximal\}$
	\item $\minimal\in N$ implies that $\minimal=\min N$
	\item $\maximal\in N$ implies that $\maximal=\max N$
	\end{itemize}
\item Morphisms are weakly monotone maps $f$ such that for all $x\in\Dom f$ we have:
	\begin{itemize}[label=\itemsecondlayer]
	\item $f(x)=\minimal$ if and only if $x=\minimal$
	\item $f(x)=\maximal$ if and only if $x=\maximal$
	\end{itemize}
\end{itemize}
A \hypo-simplicial groupoid $\Dnmop\to\Grpd$ is called a \introduce{bi-relative simplicial groupoid}.
\end{Def}
\begin{Rem}
The subcategory $\Delta_\mmdecor\subset\Delta$ is clearly admissible. Observe that while the singletons $\{\minimal\}$ and $\{\maximal\}$ are not objects in $\Delta_\mmdecor$, the set $\{\minimal,\maximal\}$ \emph{is}.
\end{Rem}
Some conventions:
\begin{itemize}
\item To simplify notation we will use natural numbers to denote elements of an object $N\in\Delta_\mmdecor$ which are known to be neither $\minimal$ nor $\maximal$. So for instance we write $\{\minimal,0,\dots,n\}$ instead of $\{\minimal,x_0,\dots,x_n\}$ when $x_n\neq \maximal$. We use the symbols $x_0,\dots, x_n$ when we want to leave open the possibility $x_0=\minimal$ or $x_n=\maximal$.
\item We will view $\Delta_\mindecor$ and $\Delta_\maxdecor$ as full subcategories of $\Delta_\mmdecor$ consisting of those objects which do not contain $\maximal$ and $\minimal$ respectively. Note that the objects $\{\minimal,\maximal\}$ and $\{\minimal,0,\dots,n,\maximal\}$ of $\Delta_\mmdecor$ do not belong to either $\Delta_\mindecor$ or $\Delta_\maxdecor$.
\item We denote by $\Delta_{\not\in}\subset\Delta_\mmdecor$ the full subcategory of those objects which do not contain either $\minimal$ or $\maximal$.  Observe that the inclusion $\Delta_{\not\in}\hra{\Delta}$ is an equivalence of categories.
\item For a $\calX^\mmdecor\colon \Dnmop\to\Grpd$ be a bi-relative simplicial groupoid we denote by $\calX^{\not\in}$, $\calX^\mindecor$ and $\calX^\maxdecor$ the restrictions of $\calX^\mmdecor$ to the full subcategories $\Delta_{\not\in}^\op$, $\Dnop$ and $\Dmop$.
\end{itemize}

\subsubsection{The abstract categorical pairing}
Let $\calX^\mmdecor\colon \Dnmop\to\Grpd$ be a $2$-Segal bi-relative simplicial groupoid.\\
We can apply the generalized Hall construction and obtain a pointed $1$-Segal bi-relative simplicial category $\calF^\mmdecor\colon \Dnmop\to\Grpd$ with $\calF^\mmdecor_N=\prod_{i=1}^n\calV(\calX^\mmdecor_{N_i})$ for $N^n\in \Dnmop$.\\
We obtain the categorical pairing $\pair{}{}:\calF^\mmdecor_{\{\minimal,0\}}\times \calF^\mmdecor_{\{0,\maximal\}}\to\calF^\mmdecor_{\{\minimal,\maximal\}}$ by choosing an adjoint inverse for the first map in the following span of categories and functors:
\begin{equation}
	\calF^\mmdecor_{\{\minimal,0\}}\times \calF^\mmdecor_{\{0,\maximal\}}\xla{\approxeq}\calF^\mmdecor_{\{\minimal,0,\maximal\}}\xra{}\calF^\mmdecor_{\{\minimal,\maximal\}}
\end{equation}
\begin{Claim} 
Recall that $\lcaction$ and $\rcaction$ denote the categorical left or right action of $\Hall{\calX^{\not\in}}$ on $\Hall{\calX^\maxdecor}$ and $\Hall{\calX^\mindecor}$ respectively. The pairing $\pair{}{}$ comes equipped with coherent natural isomorphisms $\pair{}{\blank \lcaction\blank }\lLRa\pair{\blank \rcaction\blank }{}$.\end{Claim}
\begin{Prf}
We define an unbiased (i.e.\ there is no preferred choice of bracketing) $n$-ary pairing by choosing an adjoint inverse to the equivalence of categories in the following span of categories and functors:
\begin{equation}
	\bpair{}{\blank \boxtimes \dots\boxtimes\blank }{}:\calF^\mmdecor_{\{\minimal,0\}}\times \calF^\mmdecor_{\{0,1\}}\times \dots\times \calF^\mmdecor_{\{n,\maximal\}}\xla{\approxeq}\calF^\mmdecor_{\{\minimal,0,\dots,n,\maximal\}}\xra{}\calF^\mmdecor_{\{\minimal,\maximal\}}
\end{equation}
Using the same techniques as in Section~\ref{sectionmontens} and Section~\ref{MacLanepent} we can define natural isomorphisms $\pair{}{\blank \lcaction\blank }\lLRa\bpair{}{}{}\lLRa\pair{\blank \rcaction\blank}{}$ and show that they are coherent.
\end{Prf}

\subsubsection{Morphisms induce bi-relative simplicial groupoids}

There is an ample supply of bi-relative simplicial groupoids for those who know where to look.

\begin{Cstr}\label{constructionforbirelative}
Every morphism $\Phi\colon \calY\to\calX$ of simplicial groupoids gives rise to a bi-relative simplicial groupoid $\Phi^\mmdecor\colon \Dnmop\to\Grpd$:
\begin{itemize}
\item Put $\Phi^\mmdecor_{\{\minimal,\maximal\}}\coloneqq \terminal$.
\item The groupoids $\Phi^\mmdecor_{\{\minimal,0,\dots,n\}}$, $\Phi^\mmdecor_{\{\minimal,0,\dots,n,\maximal\}}$ and $\Phi^\mmdecor_{\{0,\dots,n,\maximal\}}$ are all defined to be $\calY_\enu n$.
\item The inclusions $\{\minimal,0,\dots,n\}\hra\{\minimal,0,\dots,n,\maximal\}\hla\{0,\dots,n,\maximal\}$ induce the identity $\Phi^\mmdecor_{\{\minimal,0,\dots,n\}}\xla{=}\Phi^\mmdecor_{\{\minimal,0,\dots,n,\maximal\}}\xra{=}\Phi^\mmdecor_{\{0,\dots,n,\maximal\}}$.
\item The inclusions $\{\minimal,0,\dots,n,\maximal\}\hla\{\minimal,\maximal\}$ induce the unique map $\calY_\n\to\terminal$
\item The inclusions $\{\minimal,0,\dots,n\}\hla\{0,\dots,n\}\hra\{0,\dots,n,\maximal\}$ induce the maps
$\calY_\n\xra{\Phi_\n}\calX_\n\xla{\Phi_\n}\calY_\n$.
\item All maps where the presence of the symbols $\minimal$ and $\maximal$ does not change from source to target are defined via the corresponding maps in $\calX$ or $\calY$.
\end{itemize}
This construction defines a functor
\begin{equation}\label{morphismsgivebirelativesimpl}
	\{\text{maps of simpl. $\Dop\to\Grpd$}\} \xra{\Phi\mapsto\Phi^\mmdecor} \{\text{bi-rel. simpl. $\Dnmop\to\Grpd$}\}
\end{equation}
which is easily seen to be fully faithful
\end{Cstr}
\begin{War}
The functor~\ref{morphismsgivebirelativesimpl} is not dense since it only hits bi-relative simplicial groupoids $\calX^\mmdecor\colon \Dnmop\to\Grpd$ such that $\calX^\mmdecor_{\{\minimal,\maximal\}}\simeq \terminal$ and such that the inclusions $\{\minimal,0,\dots,n\}\hra\{\minimal,0,\dots,n,\maximal\}\hla\{0,\dots,n,\maximal\}$ induce equivalences of groupoids.
\end{War}

\begin{Prop}\label{2Segalmorphismverywelldef}
	For a morphism $\Phi\colon \calY\to\calX$ of simplicial groupoids the following are equivalent:
	\begin{enumerate}
	\item The morphism $\Phi\colon \calY\to\calX$ is $2$-Segal.
	\item\label{rightrelativeunitalsdf} The right-relative simplicial groupoid $\Phi^\mindecor\colon \Dnop\to\Grpd$ is $2$-Segal.
	\item\label{leftrelativeunitalsdf} The left-relative simplicial groupoid $\Phi^\maxdecor\colon \Dmop\to\Grpd$ is $2$-Segal.
	\item\label{birelativeunitalsdf} The bi-relative simplicial groupoid $\Phi^\mmdecor\colon \Dnmop\to\Grpd$ is $2$-Segal.\qedhere
	\end{enumerate}
\end{Prop}

\begin{Rem}There is an apparent clash of notation here. Starting from a morphism $\Phi\colon \calY\to\calX$ there are two ways to obtain a right-relative simplicial groupoid $\Phi^\mindecor\colon \Dnop\to\Grpd$:
\begin{itemize}
\item We can apply \autoref{constructionforbirelative} to $\Phi$ and then restrict the resulting bi-relative simplicial groupoid $\Phi^\mmdecor\colon \Dnmop\to\Grpd$ to the subcategory $\Dnop\subset\Dnmop$.
\item We can directly create a right-relative simplicial groupoid by using the equivalence of categories from \autoref{relativeisgoodnotion}.
\end{itemize}
It is straightforward to see that the two versions agree so the notation is unambiguous. Of course a similar remark applies for $\Phi^\maxdecor$.
\end{Rem}

\begin{Prf}[\Proofof{\autoref{2Segalmorphismverywelldef}}]
This is easily done by considering the various cases in Proposition~\ref{criterionrelativeunital}.
\end{Prf}

\subsubsection{The categorical pairing}
We can now explicitly describe the abstract categorical pairing in the case of ($\calV$-)Hall modules.
\begin{Prop}
Let $\Phi\colon \calY\to\calX$ be a $2$-Segal morphism of simplicial groupoids. Then the composition
\begin{equation}
	\pair{}{}\colon \calV(\calY_\0)\times\calV(\calY_\0)\xra{\fotV{2}}\calV(\calY_\0\times \calY_\0)\xra{\pV \Delta}\calV(\calY_\0)\xra{\sV\pi}\calV(\terminal)
\end{equation}
defines a $\HallV\calX$-balanced categorical pairing between the right and left $\HallV\calX$-modules $\HallV{\Phi^\mindecor}$ and $\HallV{\Phi^\maxdecor}$.
\end{Prop}
\begin{Prf}
This is basically the same as the proof of Theorem~\ref{descriptionHallmon}; we just have to keep in mind that $\Phi^\mmdecor_{\{\minimal,\maximal\}}\coloneqq \terminal$ and that the map
$\Phi^\mmdecor_{\{\minimal,0\}}\times\Phi^\mmdecor_{\{0,\maximal\}}\la\Phi^\mmdecor_{\{\minimal,0,\maximal\}}$ appearing in the pairing span is nothing but the diagonal $\calY_\0\times\calY_\0\xla{\Delta}\calY_\0$.
\end{Prf}

\begin{dCor}\label{explicitformulaforcategoricalpairing}
	In the situation $\calV=\Vect\BC$ the categorical pairing
	\[\pair{}{}\colon [\calY_\0,\Vect\BC]\times[\calY_\0,\Vect\BC]\lra\Vect\BC\]
	is given by the formula
	\begin{equation}\label{equationincorollarypairing}
		\pair {E}{F}\coloneqq \bigoplus_{i\in I}E_i\otimes_{G_i} F_i
	\end{equation}
	where we write $\coprod_{i\in I}\B G_i\simeq \calY_\0$ and $E=(E_i)_{i\in I}, F=(F_i)_{i\in I}\in\prod_{i\in I}\lRep {G_i}\BC\simeq[\calY_\0,\Vect\BC]$.
\end{dCor}

\subsubsection{The finitary case and decategorification}
Looking at Corollary~\ref{explicitformulaforcategoricalpairing} it is clear that the categorical pairing $\pair{}{}$ restricts to a $\vect\BC$-valued pairing
\[\pair{}{}\colon [\calY_\0,\vect\BC]_\finitary\times [\calY_\0,\vect\BC]_\finitary\lra\vect\BC\]
on the subcategory of finitary functors.
\begin{dCor}
By passing to Grothendieck groups the pairing $\pair{}{}$ induces a non-degenerate and $\fHall\calX$-balanced pairing
\[\K_0\pair{}{}\colon \fHall{\Phi^\mindecor}\times\fHall{\Phi^\maxdecor}\to\BZ\]
between the right and left fat Hall modules associated to $\Phi\colon \calY\to\calX$.
\end{dCor}

Next we want to study the relationship between the pairing $\K_0\pair{}{}$ and $\apair{}{}$.\\
The first idea would be to use the dimension map; alas this will almost never be useful as the next proposition shows.
\begin{Prop}
The pairing $K_0\pair{}{}\colon \fHall{\Phi^\mindecor}\times\fHall{\Phi^\maxdecor}\to\BZ$ descends to a well defined pairing $\ol{\K_0\pair{}{}}\colon \hall{\Phi^\mindecor}\times\hall{\Phi^\maxdecor}\to\BZ$ via the dimension map if and only if $\calY_\0$ is discrete.\\
In this case the pairings $\ol{\K_0\pair{}{}}$ and $\apair{}{}$ agree.
\end{Prop}
\begin{Prf}
Assume that $\calY_\0$ is not discrete. Pick any object $a\in \calY_0$ with non-trivial automorphism group $\triv\neq G\coloneqq \Aut(a)$ and consider the trivial representation $\BC_1$ and the regular representation $\lreg G\BC$ of $G\coloneqq \Aut(a)$. Setting $1\neq g\coloneqq \#G$ we have have $\cm(g\cdot[\BC_1]) =\cm([\lreg G\BC])$. On the other hand
$\K_0\pair{g\cdot[\BC_1]}{g\cdot[\BC_1]}= g^2[\BC]\xmapsto{\cm} g^2$
and $\K_0\pair{\lreg G\BC}{\lreg G\BC}=[\lreg G\BC]\xmapsto{\cm} g\neq g^2$, hence $\K_0\pair{}{}$ cannot descend to a well defined pairing via the dimension map.\\
The converse and the \qquote{In this case...}-part are obvious if we look at Corollary~\ref{explicitformulaforcategoricalpairing}.
\end{Prf}

It turns out that the regular section respects the duality pairing:

\begin{Prop}
	The regular sections $\cmreg\colon \hall{\Phi^\mindecor}\to\fHall{\Phi^\mindecor}$ and $\cmreg\colon \hall{\Phi^\maxdecor}\to\fHall{\Phi^\maxdecor}$ carry the pairing $\apair{}{}$ on the Hall modules to the pairing $\K_0\pair{}{}$ on the fat Hall modules.
\end{Prop}
\begin{Prf}
Recall that if $\calY_0\simeq\coprod_{i\in I}\B G_i$ then $\cmreg(\varphi)\coloneqq \(\frac{\left[\lreg{G_i}\BC\right]}{\#G_i}\)_{i\in I}$ in the notation of Corollary~\ref{explicitformulaforcategoricalpairing}, hence the result follows immediately from Equation~\ref{equationincorollarypairing} and Equation~\ref{Qvaluedpairingonhallmodule}.
\end{Prf}

\ifx\isstandaloneBJIFE\undefined
\else
\newpage
\bibliographystyle{amsalpha}
\bibliography{literatur}
\end{document}
\fi

\subsection{Example: vector spaces and free group representations over $\Fone$}\label{HallFGfree}
\ifx\preambleloaded\undefined
   
   \def\isstandalone{}
\fi

We want to go back to our pet example $\calA\coloneqq \FGfree$ and compute its finitary Hall monoidal category $\Hallf{\FGfree}$ as well as some categorical Hall modules corresponding to parabolic subgroups of the symmetric group.

\subsubsection{The Hall monoidal category}
Denote by $\calX$ the simplicial groupoid of flags in $\FGfree$.\\
We have already determined in Claim~\ref{multiplicationspanFG} that the multiplication span of $\calX$ looks as follows up to equivalence:
\[\coprod\limits_{n\in \BN}\B{G\wr \Sn}\times \coprod\limits_{m\in \BN}\B{G\wr \Sn}\xla{\cong}\coprod\limits_{n,m\in \BN}\B{G\wr (\Sn\times \Sm)}\xra{\nu}\coprod\limits_{l\in \BN}\B{G\wr \SG l}\]
Hence the underlying category of $\Hallf{\FGfree}$ is the category
\[\lrep{{G\wr \SG\ast}}\BC\coloneqq \bigoplus_{n\in \BN} \lrep {G\wr \Sn}\BC\]
of representations of the various wreath products $G\wr \SG\ast $.\\
Left Kan extension along the inclusion $\B {G\wr (\Sn\times \Sm)}\hra{\B{G\wr \SG{n+m}}}$ is nothing but the induction $\Ind_{G\wr(\Sn\times \Sm)}^{G\wr \SG{n+m}}$, hence the monoidal product is given by the \introduce{induction product}
\begin{equation}\label{inductionproduct}
	E_n\boxtimes E_m=\Ind_{G\wr(\Sn\times \Sm)}^{G\wr \SG{n+m}}(E_n\otimes E_m)
\end{equation}
for representations $E_n$ and $E_m$ of $G\wr \Sn$ and $G\wr \Sm$ respectively.\\

In Corollary~\ref{hallalgofFGfreedividedpowers} we identified the Hall algebra $\hall\FGfree$ with the ring $\Gamma_\BZ[x]\coloneqq \BZ\left\{\frac{x^n}{\fact n}\Bigmid n\in \BN\right\}$ of divided powers by sending the standard basis vector $\delta_\n$  to $\frac{x^n}{\fact n}$. Under this identification the dimension map $\cm\colon K_0(\lrep{G\wr \SG\ast}\BC)\to \Gamma_\BZ[x]$ is given by ${G\wr\Sn}\actson E_n\mapsto \frac{\dim E_n}{\fact n} x^n$
and the regular section is the algebra homomorphism
\[\BQ[x]=\BQ\otimes_\BZ\Gamma_\BZ[x]\xra{\cmreg}\BQ\otimes_\BZ K_0(\lrep{G\wr \SG\ast}\BC)\]
given by mapping the monomial $x^n$ to the the class $[\lreg{G\wr\Sn}\BC]$ of the regular representation of $G\wr\Sn$.
 
\begin{Rem}
If $G=\langle \tau\mid \tau^2=1\rangle$ is the cyclic group of order $2$ then $\FGfree$ has a very nice interpretation:\\
For $X\in \FGfree$ we can define a non-degenerate pairing $\pair\blank\blank$ on $X$ with $\pair xy=1$ if and only if $y=\tau x$; this pairing admits a Lagrangian i.e.\ an isotropic subspace of dimension $\frac{1}{2}\dim X$ (any set of representatives for the non-zero $G$-orbits does the job). This construction identifies the category $\FGfree$ with the category of pairs $(X, \pair\blank\blank)$ where $X$ is an $\Fone$-vector space and $\pair\blank\blank\coloneqq X\times X\to \{0,1\}$ is a non-degenerate pairing admitting a Lagrangian.\\
This point of view plays an important role in Young's work~\cite{YoungRelative2SegalSpaces}, in particular he uses it construct (categorical) modules over the Hall algebra (resp. Hall monoidal category) of $\vect\Fq$.
\end{Rem}

\subsubsection{Parabolic Hall modules}
Fix a natural number $r$ and tuple $\varphi=(\varphi_1,\dots,\varphi_r)\in\BN^r$ which we view as a composition of $N\coloneqq |\varphi|=\sum_{i=1}^r\varphi_i$. We obtain a corresponding $N$-dimensional object $[N]\in \FGfree$ and the parabolic subgroup
\[G\wr \SG\varphi\coloneqq G\wr \(\prod_{i=1}^r \SG{\varphi_i}\)\hra G\wr \SG N=\Aut_{\FGfree}([N]).\]
To this subgroup we can associate a quotient datum $\calQ(\varphi)\coloneqq \calQ(G\wr \SG\varphi)$ as in Section~\ref{sectionsubgroupquotientdatum} which defines the categorical Hall module $\calH(G,\varphi)\coloneqq \Hallf{\calX^{\calA,\calQ(\varphi)}}$ of the Hall monoidal category $\Hallf\FGfree= (\lrep{G\wr \SG\ast}\BC,\boxtimes)$.\\
To compute $\calH(G,\varphi)$ we must understand the left relative simplicial groupoid $\calX^\varphi\coloneqq \calX^{\calA,\calQ(\varphi)}\colon \Dmop\to\Grpd$ of flags in $\FGfree$ bounded by $\calQ(\varphi)$; more precisely we need to understand its action span
\begin{equation}\label{actionspanparabolicunknown}
	\xX{0,1}\times \calX^\varphi_{\{1,\maximal\}}\xla{}\calX^\varphi_{\{0,1,\maximal\}}\xra{}\calX^\varphi_{\{0,\maximal\}}.
\end{equation}
It is not hard to check that the action span~\ref{actionspanparabolicunknown} looks as follows (up to equivalence):
\begin{equation}\label{actionspanparabolicknown}
	\coprod_{n\in\BN} \B{G\wr \Sn}\times \coprod_{\mu\leqslant \varphi}\B{G\wr \SG\mu}\lla \coprod_{\substack{\gamma,\mu\\ \gamma+\mu\leqslant \varphi}}\B{G\wr{\SG\gamma\times \SG\mu}}\lra\coprod_{\mu\leqslant \varphi}\B{G\wr \SG\mu}
\end{equation}
	with components 
	\[\B{G\wr \Sn}\times \B{G\wr \SG\mu}\hla \B{G\wr{\SG\gamma\times \SG\mu}}\hra \B{G\wr \SG{\gamma+\mu}}\]
	for $r$-tuples $\gamma,\mu\in\BN^r$ with $\gamma+\mu\leqslant \varphi$ such that $n=|\gamma|$.\\
We conclude that $\calH(G,\varphi)$ is the category $\bigoplus_{\mu\leqslant \varphi}\lrep{G\wr \SG\mu}\BC$ with the categorical $\lrep{G\wr \SG\ast}\BC$-action given by
\[E_n\lcaction E_\mu= \bigoplus_{\substack{|\gamma|=n\\ \gamma+\mu\leqslant \varphi}}\Ind_{G\wr (\SG\gamma\times \SG\mu)}^{G\wr \SG{\gamma+\mu}}\(\Res_{G\wr \SG\gamma}^{G\wr \Sn} E_n\otimes E_\mu\)\]
for representations $E_n$ and $E_\mu$ of $G\wr \Sn$ and $G\wr \SG\mu$ respectively.\\

The parabolic Hall module $h(G,\varphi)\coloneqq \hall{\calX^{\calA,\calQ(\varphi)}}$ has a $\BZ$-basis given by $\{\delta_\mu\}_{\mu\leqslant \varphi}$ where $\mu$ runs through all $r$-tuples dominated by $\varphi$; in particular $h(G,\varphi)$ has finite rank.
\begin{Claim}
Mapping $\delta_\mu$ to $\frac{X^\mu}{\fact\mu}$ (where $\fact\mu\coloneqq \prod_{i=1}^r\fact\mu$) gives an isomorphism of $\hall\FGfree\cong\Gamma_\BZ[x]$-modules
\[h(G,\varphi)\xra{\cong}\rquot{\Gamma_\BZ[X_1,\dots,X_r]}{\frac{X^\varphi}{\varphi!}}\]
where $x^n$ acts on the right hand side by $(X_1+\dots+X_r)^n=\sum_{|\gamma|=n}\frac{\fact n}{\fact\gamma}X^\gamma$.
\end{Claim}

\begin{Prf}Let $R_n$ (for $n\in\BN$) and $R_\mu$ (for a sub-tuple $\mu\leqslant \varphi$) denote the regular representations $\lreg{G\wr \Sn}\BC$ and $\lreg{G\wr \SG\mu}\BC$ respectively. 
Observe that 
\[\Res_{G\wr \SG\gamma}^{G\wr \Sn} R_n \cong \bigoplus_{\multipl{\Sn}{\SG\gamma}}R_\gamma\]
and
\[\Ind_{G\wr (\SG\gamma\times \SG\mu)}^{G\wr \SG{\gamma+\mu}}R_\gamma\otimes R_\mu\cong R_{\gamma+\mu}\]
hence in the Grothendieck group $K_0\calH(G,\varphi)$ we have the identity
\[[R_n\lcaction R_\mu]=\sum\limits_{\substack{|\gamma|=n\\\gamma+\mu\leqslant \varphi}}\frac{\fact n}{\fact \gamma}[R_{\gamma+\mu}].\]
The dimension map is a module homomorphism and maps $R_n$, $R_\mu$ and $R_{\gamma+\mu}$ to $\fact n\delta_\n\triangleq x^n$, $\fact\mu\delta_\mu\triangleq X^\mu$ and $\fact{(\gamma+\mu)}\delta_{\gamma+\mu}\triangleq X^{\gamma+\mu}=X^\gamma X^\mu$ respectively. The claim follows.
\end{Prf}

\ifx\isstandalone\undefined
\else
\newpage
\bibliographystyle{amsalpha}
\bibliography{literatur}
\end{document}
\fi

\subsection{Convolution algebras}\label{convolutionalgebras}
\ifx\preambleloaded\undefined
   
   \def\isstandaloneBAFKFJE{}
\fi

Let $H\subseteq G$ be an inclusion of finite groups. Recall that the classical Hecke algebra $\Hecke(G,H)$ is the $\BZ$-algebra of $H$-$H$-biinvariant functions $\varphi\colon G\to \BZ$ (i.e.\ for all $g\in G$ and $h',h\in H$ we have $\varphi(h'gh)=\varphi(g)$ ) with the multiplication given by the $H$-convolution product:
\begin{equation}
	(\varphi_1\conv \varphi_2)(g)\coloneqq \frac{1}{\#H}\sum\limits_{x\in G}\varphi_1(x)\varphi_2(gx^\inv)
\end{equation}
Dyckerhoff and Kapranov observed~\cite[Example 8.2.11]{DyckerhoffKapranovHigherSegalSpacesI} that we can redefine this algebra as the Hall algebra of a $2$-Segal simplicial groupoid (see Section~\ref{Heckealgebraasgroupoid} below). This means that we can apply our general machinery to construct a corresponding monoidal category and some (categorical) modules.

\subsubsection{The fiber simplicial groupoid}\label{fibersimplicialgroupoid}

\begin{Def}
Let $f\colon A\to B$ be a map of groupoids. The \introduce{fiber simplicial groupoid} $\calX^f$ associated to $f$ is defined as follows:
\begin{itemize}
\item $\calX^f_{\enu n}\coloneqq A\twotimes{B}A\twotimes{B}\cdots \twotimes{B} A$.\\
More precisely $\calX^f_{\enu n}$ is the $2$-limit of the diagram
\begin{equation}\label{2pullbackforfibergroupoid}
	\begin{tikzcd}
		A\ar["f"',dr]& A\ar[d,"f"]& \cdots\ar[dl,phantom,"\cdots"]& A\ar[dll,"f"]\\
		& B
	\end{tikzcd}
\end{equation}
	with $n+1$-copies of $A$ sitting over $B$ via $f$.
	\item Face maps are given by projecting onto the corresponding components, degeneracies are given by diagonal inclusions.
\end{itemize}
If $f\colon \B H\hra \B G$ is an inclusion $H\subseteq G$ of groups then we call $\calX^f$ the \introduce{Hecke simplicial groupoid} and denote it by $\calX^{G,H}$.
\end{Def}

\begin{Lem}
The fiber simplicial groupoid $\calX^f$ satisfies the unpointed $1$-Segal condition, hence is $2$-Segal by Corollary~\ref{oneSegalgivesunital}.
\end{Lem}
\begin{Prf}
The unpointed $1$-Segal condition is verified by induction with an easy calculation
\[\calX^f_{\enu {n+1}}=\calX^f_{\enu{n}}\twotimes{B} A\simeq \calX^f_{\enu{n}}\twotimes{A}A\twotimes{B} A\simeq\calX^f_{\enu n}\twotimes A \calX^f_{\{n,n+1\}}.\qedhere\]
\end{Prf}
\begin{Rem}
The construction $f\mapsto\calX^f$ is clearly functorial, i.e.\ induces a functor
\[\Grpd^\mors\lra \{\text{unptd. 1-Segal simpl. $\Dop\to\Grpd$}\}\qedhere\]
\end{Rem}

\begin{Lem}\label{Lemmarelativefibersimplgroupoid}
Any $2$-pullback square of groupoids, seen as a morphism $\alpha\colon f\to f'$ in $\Grpd^\mors$, induces a $2$-Segal morphism $\alpha_\star\colon \calX^f\to\calX^{f'}$ via the fiber simplicial groupoid construction.
\end{Lem}
\begin{Prf}
Omitted.
\end{Prf}

\subsubsection{The case of group actions}

A right action by a group $G$ on a set $Z$ can be interpreted as the map $q\colon \qquot ZG\to \B G$ collapsing the set $Z$ to a point. Since this map $q$ is an iso-fibration, we can compute the $2$-pullback~\ref{2pullbackforfibergroupoid} as a strict pullback. Hence we can take
\[\calX^q_\n\coloneqq \qquot{Z^\n}G,\]
where $G$ acts on $Z^\n=Z^{n+1}$ diagonally. The face and degeneracy maps are given by projections and diagonal inclusions respectively.\\
Dyckerhoff and Kapranov~\cite[\S 2.6]{DyckerhoffKapranovHigherSegalSpacesI} call this $\calX^q$ the \buzzword{Hecke-Waldhausen space} associated to the action $G\actson Z$.

\begin{Lem}
If $G$ and $Z$ are finite then the simplicial groupoid $\calX^q$ is regular.
\end{Lem}
\begin{Prf}
All the finiteness conditions are automatically satisfied since $G$ and $Z$ are finite. The maps $\nu_\n$ are faithful since on morphisms they are given by $g\mapsto g$.
\end{Prf}
The multiplication span for the groupoid $\calX^q$ looks as follows
\begin{equation}\label{multiplspanforgroupaction}\begin{tikzcd}[column sep=large]
	\qquot {Z^\1} G \times \qquot {Z^\1} G&\qquot {Z^\2}G\ar[l,"{(\pr_{01},\pr_{12})}"']\ar[r,"{\pr_{02}}"]&\qquot{Z^\1}G
\end{tikzcd}\end{equation}
and the $2$-fibers of the extremal map can be computed as the strict pullback
\[\begin{tikzcd}
	\{z_0\}\times Z\times \{z_2\}\ar[d,"{\pr_{02}}"]\ar[r,hookrightarrow]&\qquot{Z^\2}G\ar[d,"{\pr_{02}}"]\\
	\{z_0\}\times\{z_2\}\ar[r,hookrightarrow]& \qquot{Z^\1}G
\end{tikzcd}\]
since $\pr_{02}$ is an iso-fibration.
Hence the product $\conv$ in $\hall{\calX^q}$ is given on functions $\hat\varphi_1,\hat\varphi_2\colon \pi_0\(\qquot{Z^\1}G\)\to\BZ$ as
\begin{equation}\label{generalconvolutioncalculation}
\begin{aligned}
(\hat\varphi_1\conv\hat\varphi_2)(z_0,z_2)&=\int_{\{z_0\}\times Z\times\{z_2\}}\hat\varphi_1(z_0,\blank)\cdot\hat\varphi_2(\blank,z_2)\\
&=\sum_{z\in Z}\hat\varphi_1(z_0,z)\hat\varphi(z,z_2).
\end{aligned}
\end{equation}

\subsubsection{The classical Hecke algebra, revisited}\label{Heckealgebraasgroupoid}
If the action of $G$ on $Z$ is transitive then picking an element $z\in Z$ identifies $Z$ with $\lquot GH$ where $H\coloneqq \Stab_G(z)$. Thus we obtain an equivalence of groupoids
\[\begin{tikzcd}
	\B H\ar[r,"\simeq"',"\star\mapsto H"]\ar[dr,"i"',hookrightarrow,]&\qquot {(\lquot GH)}G\ar[d,"q"]\ar[r,"\simeq"',"H\mapsto z"]&\qquot ZG\ar[dl,"q"]\\
	&\B G
\end{tikzcd}\]
which induces an equivalence of simplicial groupoids $\calX^{G,H}\simeq  \calX^{q}$.
\begin{dCor}\label{Heckesimplgroupoidisregular}
For an inclusion $H\subseteq G$ of finite groups the Hecke simplicial groupoid $\calX^{G,H}$ is regular.
\end{dCor}
\begin{Rem}
We have taken the easy way and required the groups $G$ and $H$ to be finite. Lemma~\ref{Heckesimplgroupoidisregular} is still true with some more precise finiteness assumptions on certain indices of subgroups $H$, $g^\inv Hg$ and their intersections.~\cite[Proposition 8.2.12]{DyckerhoffKapranovHigherSegalSpacesI}
\end{Rem}
Now we want to show that the Hecke simplicial groupoid $\calX^{G,H}$ deserves its name.

\begin{Lem}\label{Lemmaidentifydoublecosets} The assignment $(Hg,Hg')\mapsto Hg'g^\inv H$ defines a bijection
\[\pi_0\(\qquot{Z^\1}G\) \xra{\cong} \lrquot GHH\]
with inverse given by $HgH\mapsto (H,Hg)$.
\end{Lem}
\begin{Prf}
Easy calculation.
\end{Prf}
Therefore we can identify the elements $\hat\varphi\colon \pi_0\(\qquot{Z^\1}G\)\to\BZ$ of $\hall{\calX^q}$ with functions $\lrquot GHH\to\BZ$, i.e.\ $H$-$H$-biinvariant functions $\varphi\colon  G\to\BZ$. Under this identification Equation~\ref{generalconvolutioncalculation} reads
\begin{align*}
(\varphi_1\conv\varphi_2)(g)&\triangleq(\hat\varphi_1\conv\hat\varphi_2)(H,Hg)=\sum\limits_{Hx\in\lquot GH}\hat\varphi_1(H,Hx)\hat\varphi_2(Hx,Hg)\\
&\triangleq \sum\limits_{Hx\in\lquot GH}\varphi_1(x)\varphi_2(gx^\inv )\\
&=\frac{1}{\#H}\sum\limits_{x\in G}\varphi_1(x)\varphi_2(gx^\inv )
\end{align*}
hence we recover exactly the convolution product of the classical Hecke algebra.

\subsubsection{The Hall-Hecke monoidal category}
What is the (finitary) Hall monoidal category of the Hecke simplicial groupoid $\calX^{G,H}$? Let us at least compute its underlying category.
\begin{Claim}
We have an equivalence of groupoids
\[\coprod_{HgH\in \lrquot GHH}\B H_g \xra{\simeq}\({\qquot{Z^\1}G}\)\coloneqq \calX_\1^q\]
(where $H_g\coloneqq H\cap g^\inv H g$) given on objects by $\star_{H_g}\mapsto (H,Hg)$ and on morphisms by the inclusions $H_g\hra G$.
\end{Claim}
\begin{Prf}
In Lemma~\ref{Lemmaidentifydoublecosets} we have identified the isomorphism classes of $\qquot{Z^\1}G$ with double cosets of $H\subseteq G$. For each double coset $HgH\triangleq (H,Hg)$ the automorphism group in $\qquot{Z^\1}G$ consists of those $x\in G$ such that $Hx=H$ and $Hgx=Hg$, i.e.\ $x\in H_g\coloneqq H\cap g^\inv H g$.
\end{Prf}
Therefore the underlying category of the finitary Hall monoidal category $\Hallf{\calX^{G,H}}$ is the category
\[\bigoplus_{HgH\in\lrquot GHH} \lrep {H_g}\BC.\]

\subsubsection{Convolution modules}\label{sectionconvolutionmodules}
Let $H\subseteq P\subseteq G$ be an inclusion of finite groups. Inside the Hecke algebra $\Hecke(G,H)$ we can find the right-submodule (i.e.\ right-ideal) $\Hecke(G,H,P)$ of those functions $\varphi\colon G\to\BZ$ which are not only $H$-$H$-biinvariant but even $H$-$P$-biinvariant. We will now describe $\Hecke(G,H,P)$ as the (right) Hall module of a certain $2$-Segal morphism of simplicial groupoids.\\
We consider a slightly more general situation. Let $P\subseteq G$ be an inclusion of finite groups and let $Z$ be a right $G$-set viewed as the map $q\colon \qquot ZG\to \B G$. By restricting $q$ to a right action $q_P\colon \qquot ZP\to\B P$ of $P$ on $Z$ we obtain the following commutative $2$-pullback square.
\begin{equation}\label{pullbackforinducedgroupaction}
\begin{tikzcd}
\qquot ZP\rar{q_P}\ar[d,hookrightarrow]&		\B P\ar[d,hookrightarrow]\\
\qquot ZG\rar{q}&		\B G
\end{tikzcd}
\end{equation}
Using Lemma~\ref{Lemmarelativefibersimplgroupoid} we see that Diagram~\ref{pullbackforinducedgroupaction} induces a $2$-Segal morphism $\Phi_{q,P}\colon \calX^{q_P}\to\calX^{q}$ of simplicial groupoids.
\begin{Lem}
The right-relative simplicial groupoid $\Phi_{q,P}^\mindecor$ corresponding to $\Phi_{q,P}$ is regular.
\end{Lem}
\begin{Prf}
All finiteness conditions are automatic. Faithfulness conditions are satisfied because on morphisms all structure maps are of the form $P\ni p\mapsto p\in G$ or $G\ni g\mapsto g$.
\end{Prf}
Hence $\hall{\Phi_{q,P}^\mindecor}$ is a right-$\hall{\calX^q}$-module structure on the abelian group $V(\pi_0\calX^{q_P}_\0)$ of functions $\rquot ZP=\pi_0(\qquot{Z^\0}P)\to\BZ$.
\begin{Claim}
Explicitly the $\hall{\calX^q}$-action is computed by the formula
\begin{equation}\label{formulaconvolutionaction}
(\hat\eta\actiondot\hat\varphi)(z)=\suml_{x\in Z} \hat\eta(x)\hat\varphi(x,z)
\end{equation}
for $\hat\varphi\colon \pi_0(\qquot{Z^\1}G)\to\BZ$ and $\hat\eta\colon \rquot ZP\to\BZ$.
\end{Claim}
\begin{Prf}
Straightforward computation.
\end{Prf}

Let us return to the case $Z\cong\lquot GH$ which defines the Hecke algebra $\Hecke(G,H)$. The quotient $\rquot ZP$ is nothing but $\lrquot GHP$, hence we can identify $V(\pi_0\calX^{q_P}_\0)$ with the $\BZ$-module $\Hecke(G,H,P)$ of $H$-$P$-biinvariant functions.
Under this identifications Equation~\ref{formulaconvolutionaction} reads:
\begin{equation}\label{convolutionactiononHeckeGHP}
(\eta\actiondot\varphi)(g)=\suml_{Hx\in \lquot GH}\eta(x)\varphi(gx^\inv)=\frac{1}{\#H}\suml_{x\in G}\eta(x)\varphi(gx^\inv)=(\eta\conv\varphi)(g)
\end{equation}
for every $H$-$H$-biinvariant function $\varphi\colon G\to \BZ$ and every $H$-$P$-biinvariant function $\eta\colon G\to\BZ$.\\
We have thus realized $\Hecke(G,H,P)$ as the right Hall module of the $2$-Segal morphism $\Phi_{q,P}$.
\begin{Rem}
The only reason why we might need the subgroup $P$ to contain $H$ is if we want $\Hecke(G,H,P)$ to be included in $\Hecke(G,H)$. Of course we could just abstractly define $\Hecke(G,H,P)$ to consist of all $H$-$P$-biinvariant functions and define a $\Hecke(G,H)$-module structure on $\Hecke(G,H,P)$ by Equation~\ref{convolutionactiononHeckeGHP}.
\end{Rem}

\ifx\isstandaloneBAFKFJE\undefined
\else
\newpage
\bibliographystyle{amsalpha}
\bibliography{literatur}
\end{document}
\fi


\newpage\section{Equivariant polynomial functors}\label{sectionpolyfunctors}
\ifx\preambleloaded\undefined
   
   \def\isstandalone{}
\fi

\subsubsection{Motivation}
A classical tool for studying the representations of the symmetric group is the characteristic map which assigns to each representation of $\Sn$ a symmetric function and in fact gives a ring isomorphism~\cite[I.7.3]{MacDonaldSymmetricFunctionsHallPolynomials}
\[\ch\colon \K_0(\lrep{S_\ast}\BC)\xra{\cong} \Lambda,\]
where $\Lambda$ is the ususal ring of symmetric polynomial functions $\BC^\ast \to\BC$.\\
For a finite group $G$, this isomorphism can be generalized to a similar isomorphism~\cite[\S I.B.6]{MacDonaldSymmetricFunctionsHallPolynomials}
\begin{equation}\label{unliftedringiso}
	\ch\colon R(G)\coloneqq \K_0(\lrep{G\wr S_\ast}\BC)\xra{\cong} \Lambda(G),
\end{equation}
where $\Lambda(G)$ is the ring of symmetric polynomial class functions $\BC[G]^\ast \to\BC$ (and will be introduced properly in Section~\ref{subsectionsymmpolclassfunctions}).\\
It is a pity, though, that we have to forget so much structure when passing from the monoidal category $\lrep{G\wr S_\ast}\BC$ to its underlying Grothendieck-ring. So we categorify the notion of $G$-equivariant (with respect to conjugation) symmetric polynomial \emph{functions} $\BC[G]^\ast\to\BC$ to get the notion of  polynomial \emph{functors} $\lrep G\BC\to \vect\BC$ which we interpret as \roughly{$G$-equivariant} analogues of the classical polynomial functors $\vect\BC\to\vect\BC$. Further, we replace the \emph{ring} $\Lambda(G)$ by the \emph{monoidal category} $\Pol_\BC^\finitary(\lrep G\BC,\vect\BC)$ (of polynomial functors of bounded total degree). \\
Our goal will then be to construct an equivalence of monoidal categories
\begin{equation}\label{polynomialequivalenceintro}
	\Xi\colon \lrep{G\wr S_\ast }\BC\xra{\simeq}\Pol_\BC^\finitary(\lrep G\BC,\vect\BC)
\end{equation}
and show that it categorifies the ring isomorphism $\ch$. This equivalence of categories has the remarkably easy formula
\[G\wr \Sn\ni E_n\mapsto \Hom_{G\wr \Sn}(T^n(\blank^\dual),E_n)\cong E_n\otimes_{G\wr \Sn} T^n\]
where $\dual$ denotes the dual representation and $T^n(\blank)$ is the $n$-th tensor power. This means that every polynomial functor $F$ can be written as an expression of the form
\[F\cong\bigoplus_{n\in \BN} E_n\otimes_{G\wr \Sn} T^n(\blank)\]
which we can interpret as some sort of power series with coefficients $E_n$.\\
Hence the equivalence of categories~\ref{polynomialequivalenceintro} can also be seen as a categorification of the isomorphism
\[\BC[X]\cong\bigoplus_{n\in\BN} \BC \xlra{\cong} \{\text{polynomial functions }f\colon \BC\to\BC\}\]
which maps a polynomial $\sum_{n\in\BN}^{\finite}e_n X^n$ (viewed as its sequence $(e_0,e_1,\dots)\in\bigoplus_{n\in\BN}\BC$ of coefficients) to the polynomial function $\sum_{n\in \BN}e_n (\blank)^n\colon \BC\to\BC$.
\begin{Rem}
The notion of polynomial functors $\vect k\to\vect k$ has been studied extensively, see Kock's notes~\cite{KockNotesOnPolynomialFunctors} for a historical overview.\\
All the ingredients that go into the monoidal equivalence~\ref{polynomialequivalenceintro} are present in Macdonald's paper~\cite{MacDonaldPolynomialFunctorsAndWreathProducts}, our proof can be seen as a slightly more systematic repackaging of Macdonald's ideas.
\end{Rem}

\subsubsection{Overview}
To prove that the map \ref{polynomialequivalenceintro} is an equivalence we follow the same general strategy as the classical one for $G=\triv$~\cite[\S I.A]{MacDonaldSymmetricFunctionsHallPolynomials} and in many places the generalization works almost word-by-word. We divide the proof into a few steps:
\begin{itemize}
\item We show that each polynomial functor decomposes into a direct sums of homogeneous ones so that we only have to worry about homogeneous functors of some degree $n$.
\item We use the auxiliary category ${({\lrep G\BC})}\wr \Sn$ (see Appendix~\ref{Appendixsemidirectproduct}) and construct an equivalence of categories (Proposition~\ref{linearizationequivalence})
\[\Pol_\BC^{(1)}({({\lrep G\BC})}\wr \Sn, \vect\BC)\xra{\simeq}\Pol_\BC^n(\lrep G\BC,\vect\BC),\]
where the left side consists of multilinear functors ${(\lrep G\BC)}^n\to \vect\BC$ with some extra structure which encodes the action of the symmetric group.
\end{itemize}
The arguments in the first two steps are not specific to group representations and we will carry them out by considering general $k$-linear categories $\calG$ instead of $\lrep{G}\BC$ where the field $k$ has to have characteristic $0$.
\begin{itemize}
\item Embedded in the category $\Polone$ of componentwise linear functors we find the nicer category $\Rep^\on$ of componentwise representable functors. We will study such multi-representable functors abstractly and derive an equivalence of categories (\ref{GmodtonRep})
\[[\calD\wr \Sn, \calV]\xra{\simeq} \Rep^\on([\calD,\calV]_\op\wr \Sn,\calV)\]
in a completely formal way for any closed symmetric monoidal ground category $\calV$ (in our case $\calV\coloneqq \Vect\BC$ is the category of \emph{all} $\BC$-vector spaces) and any $\calV$-category $\calD$ (in our case $\calD\coloneqq \B \BC[G]$; hence $[\calD\wr \Sn,\calV]$ and $[\calD,\calV]$ are nothing but $\lRep{G\wr \Sn}\BC$ and ${(\lRep G\BC)} \wr \Sn)$ respectively.
\item If $\calD=\B A$ is a finite dimensional algebra (e.g.\ $A\coloneqq \BC[G]$) then we can pass from the case of all vector spaces and all modules to the finite dimensional case and obtain an equivalence of categories (Proposition~\ref{tensordecompositionfinite})
\[\lmod {A\wr \Sn} \xra{\simeq} \Rep^\on(\lmod A^\op\wr \Sn,\lmod R)\]
\item Combining the previous steps we obtain a fully faithful functor
\[\Xi_n:\lmod{A\wr \Sn}\hra \Pol^n_\BC(\lmod A^\op, \vect\BC )\]
If we are in the semisimple case (e.g.\ for $G$-representations) then all linear functors are representable, so we don't actually lose anything by passing from $\Polone$ to $\Rep^\on$; hence $\Xi_n$ is an equivalence in this case.\\
In the case of finite groups we can drop the $(\blank)^\op$, since the category of finite dimensional $G$-representations is self-dual.
\item Finally we need to check that the the resulting equivalences $\Xi_n$ assemble to the desired equivalence 
\[\Xi\colon \lrep{G\wr S_\ast }\BC\xra{\simeq}\Pol_\BC^\finitary(\lrep G\BC,\vect\BC)\]
of monoidal categories (Proposition~\ref{wreathmainequivalenceismonoidal}).
\end{itemize}

We will also see in which sense the constructed equivalence $\Xi$ of monoidal categories categorifies the ring isomorphism~\ref{unliftedringiso} in the case where $G$ is abelian:\\
When we want to decategorify a polynomial functor $F\colon \vect\BC\to\vect \BC$ into a polynomial function $f\colon \BC\to\BC$ the obvious way is to define $f(x)=\tr F(\cdot x\colon k\to k)$. The obvious generalization to polynomial functors $F\colon \lrep G\BC\to \vect\BC$ is the formula $f(x)=\tr F(\cdot x\colon \lreg G\BC\to \lreg G\BC)$ to obtain a polynomial function $f\colon \BC[G]\to\BC$.\\
In this spirit we can build a canonical decategorification map
\[\chi\colon \K_0(\Pol^\finitary_\BC(\lrep G\BC))\xra{\cong} \Lambda(G)\]
and show that it is an isomorphism and that it turns the equivalence $\Xi$ (\ref{polynomialequivalenceintro}) into the ring isomorphism $\ch$ (\ref{unliftedringiso}).\\
Actually, we will proceed the other way around: We show that map $\chi$ corresponds to $\ch$ under the identification coming from \ref{polynomialequivalenceintro} and conclude that it is a ring isomorphism.

\ifx\isstandalone\undefined
\else
\newpage
\bibliographystyle{amsalpha}
\bibliography{literatur}
\end{document}
\fi

\subsection{Polynomial functors and linearization}\label{Part1ofbigequiv}
\ifx\preambleloaded\undefined
   
   \def\isstandaloneBFHER{}
\fi


Let $k$ be a field of characteristic $0$ and denote by $\vect k$ the category of finite dimensional vector spaces over $k$. Let $\calG$ and $\calV$ be $k$-linear categories and $\cal C$ be an arbitrary category. Let $r$ and $n$ be natural numbers.\\
Let $\varphi=(\varphi_1,\dots,\varphi_r)$ denote $r$-tuples of natural numbers with $n=|\varphi|\coloneqq \varphi_1+\dots+\varphi_r$. If $r=n$ then we write $(1)$ for the tuple $(1,\dots,1)$. For a tuple $\lambda=(\lambda_1,\dots,\lambda_r)\in k^r$ (which we view as a variable) we abbreviate $\lambda^\varphi\coloneqq \lambda_1^{\varphi_1}\cdots\lambda_r^{\varphi_r}$.\\
For a tuple $X_1,\dots,X_r\in\calG$ we denote by $\pi_i=\pi^X_i\colon X_1\oplus\dots\oplus X_r\thra X_i$ and by $\iota_i=\iota^X_i\colon X_i\hra X_1\oplus\dots\oplus X_r$ the $i$-th projection and the $i$-th inclusion respectively.\\
We denote by $(\lambda)=(\lambda)_X$ both the automorphism $(\lambda_1,\dots,\lambda_r)=\sum\limits_{i=1}^r\lambda_i\iota_i\pi_i$ of $X_1\oplus\dots\oplus X_r$ in $\calG$  and  the automorphism $(\lambda_1,\dots,\lambda_r)$ of $(X_1,\dots,X_n)$ in $\calG^r$.\\

\begin{Def}
\begin{itemize}
	\item We call a functor $F\colon\calG\to\calV$ \introduce{polynomial} if for all objects $X,Y\in\calG$ the mapping $F\colon \Hom(X,Y)\to\Hom(F(X), F(Y))$ of $k$-vector spaces is polynomial, i.e.\ for every $f_1,\dots,f_s\colon X\to Y$ the expression $F(\lambda_1f_1+\dots+\lambda_rf_s)$ is a polynomial in the variables $\lambda_1,\dots,\lambda_s\in k$ with coefficients in $\Hom(F(X),F(Y))$.
	\item If all this polynomials are homogeneous of degree $n$ then we say that $F$ is \introduce{homogeneous of degree $n$}. We denote by $\Pol_k^n(\calG,\calV)$ the category of polynomial functors that are homogenous of degree $n$.
	\item A functor $\calG^r\to \calV$ is called \introduce{homogeneous of degree $\varphi$} if it is homogenous of degree $\varphi_i$ in the $i$-th variable for $i=1,\dots,r$. We denote by $\Pol_k^\varphi(\calG^r,\calV)$ the category of such functors. If $n=r$ and $\varphi=(1)$ then we call such functors \introduce{$n$-multilinear}.\qedhere
	\end{itemize}
\end{Def}

If $\calC$ is a category then we can form the wreath product $\calC\wr \Sn$ by taking the category $\calC^n$ and adding some morphisms (but no additional objects) coming from the permutation action of $\Sn$ on $\calC^n$. In the case where $\calC=\B G$ is a group we recover the usual notion of wreath product, i.e.\ we have a canonical isomorphism $(\B G)\wr \Sn\cong \B(G\wr \Sn)$. See Appendix~\ref{Appendixsemidirectproduct} for a systematic description of wreath products (and more generally of semidirect products) of categories.

\begin{Def}
 We call a functor $\calG\wr \Sn\to\calV$ \introduce{homogeneous of degree $n$} or \introduce{$n$-multilinear} if it has the corresponding property after restricting to the subcategory $\calG^n\subset\calG\wr \Sn$. Denote by $\Pol_k^n(\calG\wr \Sn, \calV)$ and $\Pol_k^{(1)}(\calG\wr \Sn,\calV)$ the corresponding categories.
\end{Def}
We will always abbreviate $\Pol_k^\? (\?)$ for $\Pol_k^\? (\?,\vect k)$.

\subsubsection{Decomposition into homogeneous pieces}
Whenever we have a polynomial, we can decompose it into a sum of monomials, i.e.\ homogeneous pieces. Now we want to do the same thing for polynomial functors.

\begin{Lem}\label{LemmaFoflambdaX}
	If $F\colon \calG^r\to\vect k$ is homogeneous of degree $\varphi$ then $F((\lambda)_X)=\lambda^\varphi\Id_X$
\end{Lem}
\begin{Prf}
	Writing $(\lambda_i)_X$ for the map $(1,\dots,\lambda_i,\dots,1)\colon X\to X$ we can decompose $(\lambda)_X$ as $(\lambda_1)_X\cdots(\lambda_r)_X$. Since $F$ is homogeneous of degree $\varphi_i$ in the $i$-th coordinate we obtain $F((\lambda_i)_X)=\lambda_i^{\varphi_i}$; the claim follows since $F$ respects compositions.
\end{Prf}
\begin{Prop}\label{decompositionintohomogeneouspieces}
	The functor $\bigoplus\limits_{|\varphi|=n}\colon \prod\limits_{|\varphi|=n}\Pol_k^\varphi(\calG^r)\lra\Pol_k^n(\calG^r)$ is an equivalence of categories.
\end{Prop}
\begin{Prf}
	The functor $\bigoplus\limits_{|\varphi|=n}$ is clearly faithful.\\
	Every morphism $\bigoplus\limits_{|\varphi|=n}F_\varphi\to\bigoplus\limits_{|\psi|=n}G_\psi$ is a matrix of transformations $\alpha_{\varphi,\psi}\colon F_\varphi\to G_\psi$. To see fullness we need to prove that all the entries off the diagonal are zero, or in other words that every natural transformation $\alpha\colon F\to G$ is zero if $F\in\Pol_k^\varphi$, $G\in\Pol_k^{\psi}$ and $\varphi\neq\psi$.
	To prove this let $X\in\calG^r$ be any object. By naturality of $\alpha$ we have
	\[\alpha_X\circ F(\lambda\Id_X)=G(\lambda\Id_X)\circ\alpha_X\]
	Using Lemma~\ref{LemmaFoflambdaX} we conclude $\lambda^\varphi\alpha_X=\lambda^\psi\alpha_X$ for all $\lambda\in k^r$, hence $\alpha_X=0$ if $\varphi\neq\psi$.\\
	Lastly we want to prove essential surjectivity; to do this let $F\colon \calG^r\to\calV$ be a polynomial functor of homogeneous degree $n$. 
	By the conditions on $F$, we can write 
	\begin{equation}\label{decompositionofFprimeoflambda}
		F((\lambda)_X)=\sum\limits_{|\varphi|=n}u_{\varphi}(X)\lambda^\varphi
	\end{equation}
	for some coefficients $u_\varphi(X)\in\Hom(F(X), F(X))$.\\
	Now take any map $f\colon X\to Y$ between objects $X, Y\in\calG^n$. Then we have
	\[\sum\limits_{|\varphi|=n}F(f)u_\varphi(X)\lambda^\varphi=F(f)F((\lambda)_X)=F((\lambda)_Y)F(f)=\sum\limits_{|\varphi|=n}u_\varphi(Y)F(f)\lambda^\varphi.\]
	Comparing coefficients gives $F(f)u_\varphi(X)=u_\varphi(Y)F(f)$, i.e.\ $u_\varphi\colon F\to F$ is a natural transformation.\\
	For two tuples of scalars $(\lambda)_X, (\mu)_X$ we have
	\[\sum\limits_{|\varphi|=n}u_\varphi(X)(\lambda\mu)^\varphi=F((\lambda\mu)_X)=F((\lambda)_X)F((\mu)_X)=\left(\sum\limits_{|\varphi|=n}u_\varphi(X)\lambda^\varphi\right)\left(\sum\limits_{|\varphi|=n}u_\varphi(X)\mu^\varphi\right)\]
	Comparing coefficients we see that we must have $u_\varphi^2=u_\varphi$ and $u_\varphi u_{\psi} =0$ for $\varphi\neq \psi$. Moreover by putting $\lambda_1=\dots=\lambda_n=1$ in Equation~\ref{decompositionofFprimeoflambda} we obtain $\sum\limits_{|\varphi|=n}u_\varphi=1_F$.\\
	This means that the $u_\varphi$ form a complete set of pairwise orthogonal idempotents, hence we have a direct sum decomposition
	\begin{equation}\label{decompositionofFprime}
		F=\bigoplus\limits_{|\varphi|=n}F_\varphi
	\end{equation}
where $F_\varphi=\Im(u_\varphi)$. This shows that $F$ lies in the essential image of $\bigoplus\limits_{|\varphi|=n}$.
\end{Prf}

If $F\colon \calG\to\calV$ is a polynomial functor (not necessarily homogeneous) then in a similar spirit we can define idempotent natural endo-transformations $u_n$ for $n\in \BN$ by taking $u_n(X)$ to be the coefficient of $\lambda^n$ in the polynomial expression $F(\lambda\colon X\to X)$. 
One can use exactly the same methods as in the proof of Proposition~\ref{decompositionintohomogeneouspieces} to show that we can decompose $\Pol_k(\calG)$ into the various $\Pol_k^n(\calG)$ and write every polynomial $F$ as direct sum $\bigoplus\limits_{n\in \BN}F_n$, where $F_n$ is the image of the idempotent endomorphism $u_n(X)$ and is homogeneous of degree $n$. In other words the functor
\[{\bigoplus\limits_{n\in \BN}}\colon \prod\limits_{n\in \BN}\Pol_k^n(\calG)\xra{\simeq} \Pol_k(\calG)\]
is an equivalence of categories.

\begin{Rem}
	Note that the direct sum of functors $F=\bigoplus\limits_{n\in \BN}F_n$ might be infinite even though $\vect k$ does not admit infinite direct sums. For example consider the exterior power functor
	$\Lambda^\ast=\bigoplus\limits_{n\in \BN}\Lambda^n\colon \vect k\to\vect k$ sending a finite dimensional vector space $X$ to the finite dimensional (!) space $\Lambda^\ast X=\bigoplus\limits_{n\in\BN}\Lambda^nX$.\\
	Evaluating a polynomial functor $F$ at a fixed object $X\in \calG$ will always result in a finite sum decomposition $F(X)=\bigoplus\limits_{n\in\BN}^{\finite}F_n(X)$ because only finitely many $u_n(X)$ can appear in the \emph{polynomial} expression $F(\lambda\colon X\to X)$.
\end{Rem}

Sometimes we are interested in those polynomial functors which can be written as a finite sum of homogenous ones. We call such a polynomial functor \introduce{finitary} or \introduce{of finite total degree} and we denote by $\Pol_k^\finitary(\calG)\subset\Pol_k(\calG)$ the full subcategory of all finitary polynomial functors. Clearly we get an equivalence of categories 
\[{\bigoplus\limits_{n\in \BN}}\colon \bigoplus\limits_{n\in \BN}\Pol_k^n(\calG)\xra{\simeq} \Pol_k^\finitary(\calG).\]

\subsubsection{Lifting the multilinear component to the wreath product}
From now on we assume $n=r$. \\
For each functor $F\in\Pol_k^n(\calG\wr \Sn)$, we have the restricted functor $F'=\restr F{\calG^n}$ which we decompose as $F'=\bigoplus\limits_{|\varphi|=n}F'_\varphi\colon \calG^n\to\calV$.\\
Let $\alpha=(\alpha,\sigma)\colon X\to Y$ be a morphism in $\calG\wr \Sn$. Then $F(\alpha)$ is a morphism $\bigoplus\limits_{|\varphi|=n}F'_\varphi(X)=F'(X)\to F'(Y)=\bigoplus\limits_{|\varphi|=n}F'_\varphi(Y)$ and we can ask ourselves whether this morphism descends to the summands $F'_{(1)}(X)\to F'_{(1)}(Y)$.\\
For a tuple $(\lambda)$ of scalars we have the equation $\alpha\circ(\lambda)_X=(\sigma\lambda)_Y\circ\alpha$
of morphisms $X\to Y$.
where $(\sigma\lambda)\coloneqq (\lambda_{\sigma^{-1}(1)},\dots,\lambda_{\sigma^{-1}(n)})$. Applying $F$ and using the decomposition $F((\lambda))=F'((\lambda))=\sum\limits_{|\varphi|=n}u_\varphi\lambda^\varphi$ we obtain
\[\sum\limits_{|\varphi|=n}F(\alpha)u_\varphi(X)\lambda^\varphi=\sum\limits_{|\varphi|=n} u_\varphi (Y)F(\alpha) (\sigma\lambda)^\varphi.\]
Notice that $\lambda^{(1)}=(\sigma\lambda)^{(1)}$, hence comparing the coefficients of $\lambda^{(1)}$ gives $F(\alpha) u_{(1)}(X)=u_{(1)}(Y)F(\alpha)$.
But $u_{(1)}$ is the idempotent corresponding to the summand $F'_{(1)}$ of $F'$ and so we see that $F(\alpha)$ descends to a map $F'_{(1)}(X)\to F'_{(1)}(Y)$ which we shall call $F_{(1)}(\alpha)$.\\
In other words we have constructed an extension $F_{(1)}\colon \calG\wr \Sn\to \calV$ of the functor $F'_{(1)}\colon \calG^n\to \calV$.
\begin{Claim}The construction $F\mapsto F_{(1)}$ is functorial in $F$; hence it provides a lift of functors in the following commutative diagram:
\begin{equation*}
	\begin{tikzcd}
		\Pol_k^n(\calG^n)\rar{\simeq}&\prod\limits_{|\varphi|=n}\Pol_k^\varphi(\calG^n)\rar{\pr_{(1)}}&\Polone(\calG^n)\\
		\Pol_k^n(\calG\wr \Sn)\uar{\Res}\ar[dashed]{rr}&&\Polone(\calG\wr \Sn)\uar{\Res}\\
	\end{tikzcd}
\end{equation*}\end{Claim}
\begin{Prf}
Consider a natural transformation $\eta\colon F\to G$ between functors $F,G\in \Pol_k^n(\calG\wr \Sn)$. Then $\eta$ is a fortiori a natural transformation $F'\to G'$ and thus induces a natural transformation $\eta_{(1)}\colon F'_{(1)}\to G'_{(1)}$. Now it is straightforward to show that $\eta_{(1)}$ is natural also with respect to the additional arrows $\alpha\colon X\to Y$ in $\calG\wr \Sn$, since $\eta_{(1)}(X)$, $\eta_{(1)}(Y)$, $F_{(1)}(\alpha)$ and $G_{(1)}(\alpha)$ are just restrictions of  $\eta(X)$, $\eta(Y)$, $F(\alpha)$ and $G(\alpha)$ respectively; hence the naturality equation is inherited.
\end{Prf}

\subsubsection{Linearization gives an equivalence of categories}
Recall that we have the functor $\Delta\colon \calG\times \B \Sn \to \calG\wr \Sn$ induced by the invariant diagonal embedding $\Delta\colon \calG\to \calG^n$ (see Appendix~\ref{appendixwreathproducts}). Pulling back along this functor gives a functor 
\[[\calG\wr \Sn,\vect k]\xra{\Delta^\star}[\calG\times\B \Sn,\vect k]\cong [\calG,\lrep \Sn k],\]
where $\lrep \Sn k[\B \Sn,\calV]$ is the category of $\Sn$-representations over $k$.\\
It is clear that if $L\colon \calG\wr \Sn\to \vect k$ is polynomial of homogeneous degree $(1)$ then $\Delta^\star(L)\colon \calG\to\lrep \Sn k$ is polynomial of homogenous degree $n$. Hence we get an induced functor
\[\Delta^\star\colon \Polone(\calG\wr \Sn)\to\Pol_k^n(\calG,\lrep \Sn k).\]
The category $\lrep \Sn k$ comes equipped with the functor
\[(-)^{\Sn}=\lim_{\B \Sn}\colon \lrep \Sn k\to\vect k\]
of $\Sn$-invariants which assigns to each diagram $\B \Sn\to \vect k$ (a.k.a. $\Sn$-representation) its limit.
\begin{Prop}\label{linearizationequivalence}
	The functors
	\[\calL\colon \Pol_k^n(\calG)\xra{(\bigoplus\colon \calG\wr \Sn\to\calG)^\star}\Pol_k^n(\calG\wr \Sn)\xra{(-)_{(1)}}\Polone(\calG\wr \Sn)\] 
	and \[D\colon \Polone(\calG\wr \Sn)\xra{\Delta^\star}\Pol_k^n(\calG,\lrep\Sn k)\xra{\((-)^{\Sn}\)_\star}\Pol_k^n(\calG)\]
	form mutually inverse equivalence of categories.
\end{Prop}
\begin{Rem}
The functor $\calL\colon  \Pol_k^n(\calG)\to\Polone(\calG\wr \Sn)$ is called \introduce{linearization} because it takes polynomial functors and turns them into multi-linear functors. We can view multi-linear functors $\calG\wr {\Sn}\to\vect k$ as \emph{linear} functors $\calG^\on \rtimes \Sn\to\vect k$ where now $\rtimes$ is the $k$-enriched version of the semi-direct product of categories.
\end{Rem}
\begin{Prf}
	We show that the composition $\calL\circ D$ is naturally isomorphic to the identity.\\
	Let $L\in \Polone(\calG\wr \Sn)$ be a functor which is linear in each variable. Let $F=D(L)$, i.e.\ $F\colon X\mapsto L(X,\dots,X)^{\Sn}$.
	Then the image $F$ under $\bigoplus^\star$ is given by
	\[F'\colon(X_1,\dots,X_n)\mapsto L(X_1\oplus\dots\oplus X_n,\dots,X_1\oplus\dots\oplus X_n)^{\Sn}\] Since $L$ is linear in each argument, we can write 
	\[L(X_1\oplus\dots\oplus X_n,\dots,X_1\oplus\dots\oplus X_n)=\bigoplus\limits_{0\leqslant\gamma_1,\dots\gamma_n\leqslant n}L(X_{\gamma_1},\dots, X_{\gamma_n})\]
	The question now is: what is the $(1)$-homogeneous part $F'_{(1)}$ of $F'$?\\
	The endomorphism $(\lambda)=(\lambda_1,\dots,\lambda_n)$ of $(X_1,\dots,X_n)$ induces the map
	\[\bigoplus\limits_{0\leqslant\gamma_1,\dots,\gamma_n\leqslant n}L(X_{\gamma_1},\dots, X_{\gamma_n})\xra{\bigoplus L(\lambda_{\gamma_1},\dots,\lambda_{\gamma_n})}\bigoplus\limits_{0\leqslant\gamma_1,\dots,\gamma_n\leqslant n}L(X_{\gamma_1},\dots, X_{\gamma_n}),\]
	where $L(\lambda_{\gamma_1},\dots,\lambda_{\gamma_n})$ is nothing but $\lambda_{\gamma_1}\cdots\lambda_{\gamma_n}$ by $(1)$-linearity of $L$. We deduce that the $(1)$-homogenous part of $\bigoplus\limits_{0\leqslant\gamma_1,\dots,\gamma_n\leqslant n}L(X_{\gamma_1},\dots, X_{\gamma_n})$ is given by the summands corresponding to those $\gamma$, where $\lambda_{\gamma_1}\cdots\lambda_{\gamma_n}=\lambda_1\cdots\lambda_n$; this means precisely that $\gamma\colon \{1,\dots,n\}\to\{1,\dots,n\}$ is a bijection, i.e.\ $\gamma\in\Sn$.\\
	We conclude that $F'_{(1)}$ is given by
	\[F'_{(1)}\colon (X_1,\dots,X_n)\mapsto\(\bigoplus\limits_{\sigma\in \Sn}L(X_{\sigma(1)},\dots,X_{\sigma(n)})\)^{\Sn}\]
	since the $\Sn$-action is compatible with the restriction to the $(1)$-homogenous component.\\
	The $\Sn$-action simply identifies the various summands $L(X_{\sigma(1)},\dots,X_{\sigma(n)})$, hence we have a natural isomorphism
	\[L(X_1,\dots,X_n)\xra[\cong]{(L(P_\sigma))_{\sigma\in \Sn}}\(\bigoplus\limits_{\sigma\in \Sn}L(X_{\sigma(1)},\dots,X_{\sigma(n)})\)^{\Sn}=F'_{(1)}(X_1,\dots, X_n).\]
	It is straightforward to check that the resulting isomorphism $L\xra{\cong} F'_{(1)}=(\cal L\circ D)(L)$ is natural in $L$, thus the proof of $\cal L\circ D\cong \Id$ is completed.\\
	
	Next we show $D\circ \cal L\cong \Id$.\\
	Let $F\in \Pol_k^n(\calG)$ be a polynomial functor of homogenous degree $n$. Denote by $F'$ the functor $\bigoplus\nolimits^\star(F)\colon (X_1,\dots,X_n)\mapsto F(X_1\oplus\dots\oplus X_n)$ and abbreviate $L_F\coloneqq\Delta^\star(F'_{(1)})$. Then we need to show that $L_F(X)^{\Sn}\cong F(X)$ naturally in $F$ and $X$.\\
	Recall that $L_F(X)$ was defined to be the direct summand of $F(X^n)$ induced by the idempotent endomorphism $u\coloneqq u_{(1)}$ which was defined to be the coefficient of $\lambda_1\cdots\lambda_n$ in the polynomial $F((\lambda_1,\dots,\lambda_n)\colon X^n\to X^n)$. Recall the notation $\iota_j\colon X\to X^n$ and $\pi_j\colon X^n\to X$ for the $j$-th inclusion and projection respectively. Denote by $\Delta=\sum_j\iota_j\colon X\to X^n$ and $\nabla=\sum_j\pi_j\colon X^n\to X$ the diagonal and codiagonal on $X$ respectively and let $P_\sigma=\sum_j\iota_{\sigma(j)}\pi_i\colon X^n\to X^n$ be the permutation map associated to $\sigma\in \Sn$.\\
	\begin{Claim}\label{equationclaimsinlinearization} We have the following equations:
	\begin{align}
		u F(\Delta)F(\nabla) u = u\sum\limits_{\sigma\in \Sn}F(P_\sigma) u.\label{nablaDeltaissumsigmasonL}\\
		F(\nabla)uF(\Delta) = \fact n\label{NablaUDeltaisNfactorial}
	\end{align}\end{Claim}
	\begin{Prf}[\Proofof{Claim~\ref{equationclaimsinlinearization}}]
	Consider the expression $F(\sum_{ij}\xi_{ij}\iota_j\pi_i)$ which is a polynomial in the variables $\xi_{ij}\in k$ and denote by $w_\sigma$ the coefficient in front of $\xi_{\sigma(1),1}\cdots\xi_{\sigma(n),n}$.
	The left side of Equation~\ref{nablaDeltaissumsigmasonL} is the coefficient of $\lambda_1\cdots\lambda_n\cdot\mu_1\cdots\mu_n$ in the expression $F((\lambda)\Delta\nabla(\mu))=F(\sum_{ij}\lambda_j\mu_i\iota_j\pi_i)$. By putting $\xi_{ij}=\lambda_j\mu_i$ we see that this coefficient is precisely $\sum_\sigma w_\sigma$.\\
	On the other hand $u F(P_\sigma) u$ is the coefficient of $\lambda_1\cdots\lambda_n\cdot\mu_1\cdots\mu_n$ in the expression $F((\lambda)P_\sigma(\mu))=\sum_i\lambda_{\sigma(i)}\mu_i\iota_{\sigma(i)}\pi_i$. By putting $\xi_{ij}=\lambda_{\sigma(i)}\mu_i$ if $j=\sigma(i)$ and $\xi_{ji}=0$ otherwise we see that this coefficient is nothing but $w_\sigma$.
	By summing over all $\sigma\in \Sn$ we get Equation~\ref{nablaDeltaissumsigmasonL}.\\
	To get Equation~\ref{NablaUDeltaisNfactorial} observe that the left side is the coefficient of $\lambda_1\cdots\lambda_n$ in $F(\nabla(\lambda)\Delta)=F(\lambda_1+\dots+\lambda_n)=(\lambda_1+\dots+\lambda_n)^n$; it is elementary to see that this coefficient is indeed $\fact n$.
	\end{Prf}
	Denote by $b$ and $q$ the inclusion of and the projection onto the direct summand $L_F(X)$ of $F(X\oplus\dots\oplus X)$ respectively, hence $bq=u$  and $qb=\Id_{L_F(X)}$. We can define morphisms
	\begin{align*}
	\varepsilon&\colon F(X)\xra{F(\Delta)} F(X\oplus\dots\oplus X)\xthra{q} L_F(X)\\
	\eta&\colon L_F(X)\xhra{b} F(X\oplus\dots\oplus X)\xra{F(\nabla)} F(X).
	\end{align*}
	Equation~\ref{NablaUDeltaisNfactorial} can then be rewritten as $\eta\varepsilon = \fact n$; restricting Equation~\ref{nablaDeltaissumsigmasonL} to $L_F(X)$ gives $\varepsilon\eta = \sum_\sigma L_F(P_\sigma)$.
	We know that $\frac{1}{\fact n}\sum_\sigma L_F(P_\sigma)\colon L_F(X)\to L_F(X)^{\Sn}$ is a retraction of the inclusion $L_F(X)^{\Sn}\to L_F(X)$; hence defining
	\begin{align*}
		\varepsilon'&\colon F(X)\xra{\varepsilon}L_F(X)\xthra{\frac{1}{\fact n}\sum_\sigma L_F(P_\sigma)}L_F(X)^{\Sn}\\
		\eta'&\colon L_F(X)^{\Sn}\hra L_F(X)\xra{\eta} F(X)
	\end{align*}
	gives $\varepsilon'\eta'= \fact n$ and $\eta'\varepsilon'=\eta\(\frac{1}{\fact n}\varepsilon\eta\)\varepsilon=\fact n$. Therefore $\varepsilon'$ and $\eta'$ are mutually inverse up to the invertible scalar $\fact n\in k^\units$.\\
	It is clear that $\varepsilon$ and $\eta$ (hence $\varepsilon'$ and $\eta$) are natural in $X$ and $F$ since they are constructed out of the maps $F(\Delta)$, $F(\nabla)$, $b$, $q$ and $L_F(P_\sigma)$; all of which are natural in $X$ and $F$.
	\end{Prf}

\ifx\isstandaloneBFHER\undefined
\else
\newpage
\bibliographystyle{amsalpha}
\bibliography{literatur}
\end{document}
\fi

\subsection{Multi-representable functors}\label{Part2ofbigequiv}	
	\ifx\preambleloaded\undefined
   
   \def\isstandaloneANDJG{}
   \else
\fi

Let $\calV$ be a complete and cocomplete closed symmetric monoidal ground category, e.g.\ $\calV=\lMod R$ for any commutative ring $R$. Everything (categories, functors, natural transformation) that appears in Section~\ref{Part2ofbigequiv} shall be understood to be enriched over $\calV$ even if it is not explicitly mentioned. See Kelly's book\cite{KellyBasicConceptsOfEnrichedCategoryTheory} for an introduction to enriched category theory.\\ The symbol $\calD$ will denote a $\calV$-category.

\begin{Def}
A functor $\calD^\on\to\calV$ is said to be \textbf{$n$-representable} if it is representable when fixing all but one of the variables.\\
A functor $\calD\wr \Sn\to\calV$ is said to be \textbf{$n$-representable} if the restriction to $\calD^\on$ is $n$-representable.\\
The $\calV$-category of all $n$-representable functors is denoted by $\Rep^\on(\blank,\calV)$ with $\Rep=\Rep^{\otimes 1}$.
\end{Def}

We denote by $\bigotimes=\bigotimes_n\colon [\calD,\calV]^\on\to [\calD^\on,\calV]$ the functor given by
\begin{align*}(\rho_1,\dots,\rho_n)&\mapsto& (x_1,\dots,x_n)\mapsto \rho_1(x_1)\otimes\dots\otimes \rho_n(x_n).\end{align*}

\subsubsection{Decomposing multi-representable functors}
We use the notations $Y_\calD\colon \calD\hra[\calD_\op,\calV]$ and $Y^\calD\colon \calD\hra[\calD,\calV]_\op$ for the co- and contravariant (enriched) Yoneda-embedding respectively. See also Appendix~\ref{appendixenrichedcattheory}

\begin{Prop}\label{decomposingnrep}
The functor $\bigotimes\colon [\calD,\calV]^\on\to[\calD^\on,\calV]$ induces an equivalence of categories
\[\Rep([\calD^\on,\calV]_\op,\calV)\xra{\bigotimes_\op^\star}\Rep^\on([\calD,\calV]^\on_\op,\calV)\] with an inverse given by the restriction of
\[\left[[\calD,\calV]^\on_\op,\calV\right]\xra{\left[\(Y^\calD\)^\on,\calV\right]}[\calD^\on,\calV]\xra[Y_{[\calD^\on,\calV]}]{\simeq}\Rep([\calD^\on,\calV]_\op,\calV)\]
to the full subcategory of $n$-representable functors.
\end{Prop}
\begin{Prf}
Let $L\colon [\calD,\calV]_\op^\on\to \calV$ be an $n$-representable functor.
	Then we need to show that $L$ can be written as the composition
	\[[\cal D,\calV]^\on_\op\xra{\bigotimes_\op}[\cal D^\on,\calV]_\op\xra{\hat L}\calV,\]
	where $\hat L=\hom{[\cal D^\on,\calV]}(\blank,\restr L{\calD^\on})$ is the functor represented by the object $\restr L{\calD^\on}\colon \cal D^\on\xhra{\(Y^\calD\)^\on}[\cal D,\calV]^\on_\op\xra{L}\calD$ of $[\calD^\on,\calV]$.
	For functors $\rho_1,\dots,\rho_n\in [\calD,\calV]$ we calculate:
	\begin{align*}
		L(\rho_1,\dots,\rho_n) &\cong \hom{[\calD,\calV]}(\rho_1,\restr{L(\blank,\rho_2,\dots,\rho_n)}\calD)\\
		&=  \hom{[\calD,\calV]}(\rho_1, d_1\mapsto L(\hom\calD(d_1,\blank),\rho_2,\dots,\rho_n))\\
		&\cong  \hom{[\calD,\calV]}(\rho_1, d_1\mapsto\hom{[\calD,\calV]}(\rho_2, d_2\mapsto L(\hom\calD(d_1,\blank),\hom\calD(d_2,\blank),\rho_3,\dots,\rho_n)))\\
		&\cong  \hom{[\calD\otimes\calD,\calV]}(\rho_1\otimes\rho_2, (d_1,d_2)\mapsto L(\hom\calD(d_1,\blank),\hom\calD(d_2,\blank),\rho_3,\dots,\rho_n))\\
		&\dots\\
		&\cong \hom{[\calD^\on,\calV]}(\rho_1\otimes\dots\otimes\rho_n, L(\hom\calD(d_1,\blank),\dots,\hom\calD(d_n,\blank)))\\
		&=\(\hat L\circ \bigotimes\)(\rho_1,\dots,\rho_n)
	\end{align*}
	naturally in $\rho_1,\dots,\rho_n$ and $L$. The first and second $\cong$-steps are given by applying the  Yoneda Lemma~\ref{YonedaLemma},~\ref{variantYoneda} to the representable functors $L(\blank,\rho_2,\dots,\rho_n)$ and $L(\hom\calD(d_1,\blank),\blank,\dots,\rho_n)$ respectively. The third $\cong$ follows from applying Lemma~\ref{partialhomtensoradjunction} to $\rho_1$, $\rho_2$ and $F\colon (d_1,d_2)\mapsto L(\hom\calD(d_1,\blank),\hom\calD(d_2,\blank),\rho_3,\dots,\rho_n)$ and then the rest is repeating the same steps again and again.\\
	
	On the other hand, consider a functor $\hat L\colon [\calD^\on,\calV]_\op\to\calV$, represented as $\hom{[\calD^\on,\calV]}(\blank,E)$ by an object $E\in [\calD^\on,\calV]$. We have to show that the composition
	\[\calD^\on\xra{\(Y^\calD\)^\on}[\calD,\calV]_\op^\on\xra{\bigotimes_\op}[\calD^\on,\calV]_\op\xra{\hat L}\calV\]
	is naturally isomorphic to $E$. So we calculate, for $x_1,\dots,x_n\in \calD$:
	\begin{align*}
		E(x_1,\dots,x_n) &\cong\hom{[\calD,\calV]}(\hom\calD(x_1,\blank),E(\blank,x_2,\dots,x_n))\\
		&\cong\hom{[\calD,\calV]}(\hom\calD(x_1,\blank),d_1\mapsto\hom{[\calD,\calV]}(\hom\calD(x_2,\blank),E(d_1,\blank,x_3,\dots,x_n))\\
		&\cong \hom{[\calD\otimes\calD,\calV](\hom\calD(x_1,\blank)\otimes\hom\calD(x_2,\blank),E(\blank,\blank,x_3,\dots,x_n))}\\
		&\dots\\
		&\cong \hom{[\calD^\on,\calV]}(\hom\calD(x_1,\blank)\otimes\dots\otimes\hom\calD(x_n,\blank),E).
	\end{align*}
	Here the first and the second isomorphisms are the Yoneda Lemma~\ref{YonedaLemma},~\ref{classicalYoneda} and the third is Lemma~\ref{partialhomtensoradjunction}; then we keep repeating the same steps. Since all the provided isomorphisms are natural in $x_1,\dots,x_n\in \calD$ and $E\in [\calD^\on,\calV]$ the proof is completed.
\end{Prf}

\subsubsection{Extending to the semidirect product}

Now we want to apply Lemma~\ref{restricttosemidirectprod} to the following situation:\\
Let $S\coloneqq \Sn$. We set $\calC\coloneqq \calD^\on$ and $\calC'\coloneqq [\calD,\calV]^\on_\op$ with the obvious permutation actions $\Gamma\colon \Sn\actson \calD^\on$ and $\Gamma'\colon \Sn\actson[\calD,\calV]^\on_\op$. Further declare $\calF$ to be the whole category $[\calC,\calV]$ and $\calF'$ to be the full subcategory $\Rep^\on([\calD,\calV]^\on_\op,\calV)$ of $[\calC',\calV]$.
Then the composition
\[\Phi\colon \calF\coloneqq[\calD^\on,\calV]\xra[Y_{[\calD^\on,\calV]}]{\simeq}\Rep([\calD^\on,\calV]_\op,\calV)\xra[\bigotimes^\star]{\simeq}\Rep^\on([\calD,\calV]^\on_\op,\calV)\eqqcolon \calF'\]
is an equivalence of categories by Proposition~\ref{decomposingnrep}. It is an easy calculation to see that $\Phi$ is compatible with the $\Sn$ actions on both sides.\\
Hence $Lemma~\ref{restricttosemidirectprod}$ provides:

\begin{dCor}\label{GmodtonRep}
	The functor $\Phi_{\Sn}\colon [\calD\wr \Sn, \calV]\xra{\simeq} \Rep^\on([\calD,\calV]_\op\wr \Sn,\calV)$, which explicitly is given by
	\begin{align*}E&\mapsto& (\rho_1,\dots,\rho_n)\mapsto \Mor{[\calD^\on,\calV]}{\rho_1\otimes \dots\otimes\rho_n}{E},\end{align*}
	is an equivalence of categories with inverse given by the formula
	\begin{align*}L&\mapsto& \calD\wr \Sn\xra{Y^\calD\wr \Sn}[\calD,\calV]_\op\wr \Sn\xra{L}\calV\end{align*}
\end{dCor}

\ifx\isstandaloneANDJG\undefined
\else
\newpage
\bibliographystyle{amsalpha}
\bibliography{literatur}
\end{document}
\fi

\subsection{In representation theory...}
\ifx\preambleloaded\undefined
   
   \def\isstandalone{}
\fi

Set $\calV\coloneqq \lmod R$ for some commutative ring $R$. If $\calD$ is a $\calV$- (i.e.\ $R$-)linear category, then $\calV$- (i.e.\ $R$-)linear functors $\calD\to\calV$ are what we would call representations of $\calD$ over $R$.

\subsubsection{...of finite algebras...}
 Let $A$ be an $R$-algebra, we can view it as the $R$-linear category $\calD\coloneqq \B A$.\\
Corollary~\ref{GmodtonRep} applies with $\calV\coloneqq \lMod R$, so we get an equivalence of categories
\begin{equation}\label{tensordecompositionforA}
	\lMod {A\wr \Sn} \xra{\simeq} \Rep^\on(\lMod A^\op\wr \Sn,\lMod R)
\end{equation}
given by the formula
\begin{equation}\label{forumlaforwreathprodofalgs}
E\mapsto \hat E\coloneqq \Hom_{A^\on}(\blank \otimes \dots\otimes \blank, E)
\end{equation}

\begin{Rem}
Formula~\ref{forumlaforwreathprodofalgs} is slightly ambiguous since it only describes $\hat E$ as a functor
\[\hat E\colon \lMod A_\op^\on\to\lMod R.\]
By the universal property of the wreath product, we still need to specify compatible natural transformations $Q_\sigma\colon \hat E\Ra \hat E\circ \sigma$ for all $\sigma\in\Sn$.\\
Consider the $R$-linear permutation map
\[P_\sigma(X_1,\dots, X_n)\colon X_{\sigma^\inv(1)}\otimes\dots\otimes X_{\sigma^\inv(n)}\lra X_1\otimes\dots\otimes X_n\]
which is natural in $X_1,\dots,X_n\in \lMod A$.
Define the natural transformation $Q_\sigma\coloneqq E_\sigma\circ \blank \circ P_\sigma\colon \hat E\Ra\hat E\circ \sigma$ which is given on the $(X_1,\dots,X_n)$-th component by mapping a morphism $f\colon X_1\otimes \dots\otimes X_n\to E$ to the composition
\[X_{\sigma^\inv(1)}\otimes\dots\otimes X_{\sigma^\inv(n)}\xra{P_\sigma} X_1\otimes\dots\otimes X_n\xra{f} E\xra{E_\sigma}E\]
where $E_\sigma$ is the action of $\sigma$ on $E$. A direct calculation shows that $E_\sigma\circ f\circ P_\sigma$ is an $A^\on$-module homomorphism even though $P_\sigma$ and $E_\sigma$ are only $R$-linear, therefore $Q_\sigma\coloneqq E_\sigma\circ \blank \circ P_\sigma$ is well defined.\\
It is easy to see that the various $Q_\sigma$ are compatible, hence give rise to the desired functor
\[\hat E\colon \lmod A^\op\wr \Sn\to\lMod R\qedhere\]
\end{Rem}

\begin{Prop}\label{tensordecompositionfinite}
If the algebra $A$ is finite over $R$ (i.e.\ finitely generated as an $R$-module) then Equivalence~\ref{tensordecompositionforA} restricts to an equivalence of categories
\begin{equation}\label{tensordecompositionforAfinite}
	\lmod {A\wr \Sn} \xra{\simeq} \Rep^\on(\lmod A^\op\wr \Sn,\lmod R)
\end{equation}
where we consider \emph{finitely generated} modules.
\end{Prop}
\begin{Prf}
By looking at the explicit formulas from Corollary~\ref{GmodtonRep} it is clear that under Equivalence~\ref{tensordecompositionforA} the full subcategory $\lmod{A\wr \Sn}\subset\lMod{A\wr \Sn}$ (on the left) corresponds to the full subcategory of those $n$-representable functors
\[L\colon \lMod{A}^\op\wr \Sn\lra \lMod{R}\]
which map finite modules to finite modules (on the right).\\
Since every $A$-module is a direct limit of finite $A$-modules and $n$-representable functors preserve direct limits in each variable, we see that every such functor $L$ is uniquely determined by its restriction to finite modules. The result follows.
\end{Prf}

\subsubsection{...over a field of characteristic zero}
We go back to the situation where our ground ring $R$ is a field $k$ of characteristic $0$; we still consider a finite (dimensional) $k$-algebra $A$. We abbreviate $\calG_n\coloneqq \lmod{A^\on}$ with $\calG\coloneqq \calG_1=\lmod A$.\\

Consider the following composition $\Xi_n$ of functors
\begin{equation}\label{defineXiasbigcomp}
		\begin{tikzcd}
			\lmod{A\wr \Sn}\ar[d, dashed, "\Xi_n"]\ar[r, "\simeq"] & \Rep^\on({\calG^\op}\wr \Sn,\vect k)\ar[r, hookrightarrow]& {[{\calG^\op}\wr \Sn, \vect k]}_k\ar[d, "\cong"]\\
			\Pol_k^n(\calG^\op,\vect k)&\Pol_k^n(\calG^\op,\lrep{\Sn}k)\ar[l,"\((\blank)^{\Sn}\)_\star"']&\Polone(\calG^\op\wr \Sn, \vect k)\ar[l,"\Delta^\star"']\ar[ll, "\simeq", bend left=9]
		\end{tikzcd}
\end{equation}

where the first equivalence is the one from Proposition~\ref{tensordecompositionfinite} and the equivalence in the bottom row is the one constructed in Section~\ref{Part2ofbigequiv} (Proposition~\ref{linearizationequivalence}). The vertical isomorphism arises from comparing the two different meanings of $\blank \wr \Sn$; in the first case we use the $\vect k$ enriched version $\calG_\op^\on\rtimes \Sn$ (and talk about $k$-linear functors), in the second case we use the non-enriched version $\calG_\op^n\wr \Sn$ (and talk about multilinear, i.e.\ $(1)$-homogeneous functors).\\
Tracking the explicit formulas for the various functors appearing in Diagram~\ref{defineXiasbigcomp} we see that the dashed composition is given by
\begin{equation}\label{formulaXi}
	\Xi_n\colon E\mapsto \(\Hom_{A^\on}\(T^n(\blank),E\)\)^{\Sn},
\end{equation}
where $T^n\colon \lmod A\to \lmod{A\wr \Sn}$ is the $n$-th tensor product functor.

\begin{dCor}
The functor $\Xi_n\colon \lrep{A \wr \Sn}k\lra \Pol_k^n(\lrep{A}k^\op,\vect k)$ given by Formula~\ref{formulaXi} is fully faithful.
\end{dCor}

If the algebra $A$ is semi-simple then every linear functor ${\lrep A k}^\op\to\vect k$ is representable, hence the inclusion $\Rep^\on({\calG^\op}\wr \Sn,\vect k)\hra {[{\calG^\op}\wr \Sn, \vect k]}_k$ appearing in the definition of $\Xi_n$ is actually the identity.

\begin{dCor}
If the algebra $A$ is semi-simple then the functor $\Xi_n$ is an equivalence of categories.
\end{dCor}

For a finite group $G$ we can identify $\lrep G\BC^\op$ with $\lrep G\BC$ via the duality functor $E\mapsto E^\dual$. Moreover there is a canonical isomorphism between invariants and coinvariants, hence we don't have to worry about the distinction.

\begin{dCor}\label{equivalenceforgroupsmain}
	For a finite group $G$ the functor
	\[\Xi_n\colon \lrep{G \wr \Sn}k\lra \Pol_k^n(\lrep Gk,\vect k)\]
	given by the formula
	\[E_n\mapsto \Hom_{G\wr \Sn}(T^n(\blank^\dual),E_n)\cong E_n\otimes_{G\wr \Sn} T^n\cong \Hom_{G\wr\Sn}(E_n^\dual,T^n(\blank))\]
	is an equivalence of categories.
\end{dCor}

\ifx\isstandalone\undefined
\else
\newpage
\bibliographystyle{amsalpha}
\bibliography{literatur}
\end{document}
\fi

\subsection{The Hall monoidal category of $\FGfree$}\label{polyfunctbacktoFGfree}
\ifx\preambleloaded\undefined
   
   \def\isstandalone{}
\fi

Recall from Section~\ref{HallFGfree} that $\lrep {G\wr S_\ast}\BC\coloneqq \bigoplus\limits_{n\in \BN}\lrep{G\wr S_n}\BC$ is the underlying category of the finitary Hall monoidal category $\Hallf{\FGfree}$.
The monoidal product is given by the induction product
\[E_n\boxtimes E_m\coloneqq \Ind_{(G\wr S_n)\times (G\wr S_m)}^{G\wr{S_{n+m}}}(E_n\otimes E_m)\]
for representations $E_n$ and $E_m$ of $G\wr S_n$ and $G\wr S_m$ respectively.\\
We can assemble the $\Xi_n$'s for all $n\in\BN$ to an equivalence
\[\Xi\coloneqq \bigoplus_{n\in\BN}\Xi_n\colon \lrep {G\wr S_\ast}\BC\xra{\simeq}\bigoplus_{n\in\BN}\Pol_\BC^n(\lrep G\BC)\simeq\Pol_\BC^\finitary(\lrep G\BC).\]

\begin{Prop}\label{wreathmainequivalenceismonoidal}
The functor $\Xi\colon E_\ast \mapsto E_\ast \otimes_{G\wr S_\ast} T^\ast$ is an equivalence of monoidal categories
\[\Xi\colon \Hallf{\FGfree}\lra \(\Pol_\BC^\finitary(\lrep G\BC), \otimes\),\]
where $\otimes$ denotes the pointwise tensor product of polynomial functors.
\end{Prop}

\begin{Prf} An easy calculation shows that $\boxtimes$ corresponds to $\otimes$ under $\Xi$:
\begin{align*}
	\Xi(E_n\boxtimes E_m)&\cong\(\Ind_{(G\wr S_n)\times (G\wr S_m)}^{G\wr{S_{n+m}}}(E_n\otimes E_m)\)\otimes_{G\wr S_{n+m}}T^{n+m}\\
	&\cong (E_n\otimes E_m)\otimes_{(G\wr S_n)\times(G\wr S_m)}\Res_{(G\wr S_n)\times (G\wr S_m)}^{G\wr{S_{n+m}}}T^{n+m}\\
	&\cong (E_n\otimes E_m)\otimes_{(G\wr S_n)\times(G\wr S_m)}(T^n\otimes T^m)\\
	&\cong (E_n\otimes_{G\wr S_n} T^n) \otimes (E_m\otimes_{G\wr S_m} T^m)\\
	&\cong  \Xi(E_n)\otimes \Xi(E_m)
\end{align*}
naturally in $E_n$ and $E_m$.
\end{Prf}

\ifx\isstandalone\undefined
\else
\newpage
\bibliographystyle{amsalpha}
\bibliography{literatur}
\end{document}
\fi

\subsection{The characteristic map}\label{sectioncharacteristicmap}
\ifx\preambleloaded\undefined
   
   \def\isstandalone{}
\fi

In this section we will write $\K_0$ for the $\BC$-valued Grothendieck group/ring which we would otherwise denote by $\BC\otimes_\BZ\K_0$

\subsubsection{Symmetric polynomial class functions}\label{subsectionsymmpolclassfunctions}
Let $R(G\wr \Sn)$ be the $\BC$-vector space spanned by the irreducible characters of $G\wr \Sn$; this is the same as the vector space of class functions $G\wr \Sn\to \BC$. Out of this vector spaces we can construct the graded ring $R(G)\coloneqq\bigoplus_{n\in \BN}R(G\wr \Sn)$ with the multiplication given by the induction product.\\ Note that $R(G)$ is just the $K_0$-ring of the monoidal category $\lrep{G\wr S_\ast}\BC$, i.e.\ the fat Hall algebra $\fHall[\BC]{\FGfree}$. Hence the equivalence of (graded) monoidal categories $\lrep {G\wr S_\ast} \BC \simeq \Pol^\finitary_k(\lrep G\BC,\vect\BC)$ provides an isomorphism of graded rings $R(G)\simeq K_0(\Pol^\finitary_k(\lrep G\BC,\vect\BC))\eqqcolon \frakF(G)$.\\
For each integer $d\geqslant 1$ we consider the ring $\calC_d(G)$ of polynomial class functions on $\BC[G]^d$, i.e.\ functions $f\colon \BC[G]^d\to \BC$ which can be written as $\(\sum_ga_{ig}g\)_{i=1}^d\mapsto \hat f((a_{ic})_{c\in G_\star, i})$ where $\hat f\in \BC[x_{ic}\mid 1\leqslant i \leqslant d; c\in G_\star]$ is a polynomial and $a_{ic}=\sum_{g\in c}a_i$. Since the irreducible characters $\gamma\in G^\star$ form a basis of the linear class functions $\BC[G]\to \BC$, we immediately obtain that $\calC_d(G)$ is the polynomial ring $\BC[y_{i\gamma}\mid 1\leqslant i \leqslant d; \gamma\in G^\star]$ in the polynomial functions $y_{i\gamma}\colon \BC[G]^d\ni (a_i)_{i=1}^d\mapsto \gamma(a_i)$.\\
Now we can build the graded ring $\Lambda(G)$ of \emph{symmetric} polynomial class functions as the inverse limit $\lim_{d\la\infty}\Lambda_d(G)$ (in the category of graded rings) of the system defined by $\Lambda_d(G)\coloneqq\calC_d(G)^{S_d}$ (where $S_d$ acts by permuting the variables) with the maps $\Lambda_{d+1}(G)\to\Lambda_{d}(G)$ given by restriction of functions.\\

\subsubsection{An intrinsic description of the characteristic map}\label{subsectioncharacteristicmap}
A main tool for understanding characters of the wreath product $G\wr \Sn$ is the characteristic map $\ch\colon R(G)\to\Lambda(G)$. We don't need the original definition since we are going to give a new one soon anyway, but we need the following facts:

\begin{enumerate}[label=(ch\roman*), ref=(ch\roman*)]
	\item\label{chisiso} The map $\ch$ is an isomorphism of graded rings.~\cite[\S I.B.6]{MacDonaldSymmetricFunctionsHallPolynomials}
	\item\label{chmapsEgammatoh} For any representation $E_\gamma$ of $G$ with character $\gamma\in G^\star $, the map $\ch$ sends the $G\wr \Sn$-representation $E_\gamma^{\otimes n}$ to the $n$-th complete symmetric polynomial $h_n(y_{\ast,\gamma})$ in the variables $y_{1,\gamma},y_{2,\gamma},\dots$.~\cite[\S I.B.8]{MacDonaldSymmetricFunctionsHallPolynomials}
	\item\label{chmapsXtoSchur} A complete set of pairwise non-isomorphic representations of $G\wr \Sn$ is given by $\{X_\lambda\mid \lambda\in \calP_n(G^\star)\}$, where $X_\lambda$ corresponds under $\ch$ to the Schur function $S_\lambda(y_{\ast,\ast})=\prod_{\gamma\in G^\star}s_{\lambda(\gamma)}(y_{\ast,\gamma})$.~\cite[\S I.B.9]{MacDonaldSymmetricFunctionsHallPolynomials}
\end{enumerate}

The goal of Section~\ref{subsectioncharacteristicmap} is to give an intrinsic definition of $\ch$ using the language of polynomial functors.\\

Let $F\colon \lrep G\BC\to\vect\BC$ be a polynomial functor and let $d\in \BN$ be the number of variables. The functor $F$ induces a polynomial function
\[\chi_d(F)\colon \BC[G]^d\lra\End_G({\lreg G\BC}^d)\xra{F}\End_k(F(\lreg G\BC))\xra{\tr}\BC,\]
where the first map is given by componentwise right-multiplication with the inverse, i.e.\ $(g_i)_i\mapsto (x_i)_i\mapsto(x_ig_i^\inv)_i$.\\
For $g_i,h_i\in G$ we have
\[\chi_d(F)((h_i^\inv g_ih_i)_i)=\tr \(F((\cdot h_i)_i)F((\cdot g_i^\inv)_i)F((\cdot h_i)_i)^\inv\)=\tr F((\cdot g_i^\inv)_i)=\chi_d(F)((g_i)_i),\]
hence $\chi_d(F)$ is a class function. Similarly, if $\sigma\in S_d$ is a permutation and $P_\sigma$ is the corresponding permutation matrix then
\[\chi_d(F)((g_{\sigma(i)})_i)=\tr F(P_{\sigma}^\inv\circ (\cdot g^\inv_i)_i\circ P_\sigma)=\tr F((\cdot g^\inv_i)_i)=\chi_d(F)((g_i)_i);\]
therefore $\chi_d(F)$ is also symmetric.\\
Since $\tr$ maps $\oplus$ to $+$ and $\otimes$ to $\cdot$ we obtain a ring homomorphism $\chi_d\colon \frakF(G)\to\Lambda_d(G)$. Clearly $\chi_d$ is graded, i.e.\ sends homogenous polynomial functors of degree $n$ to polynomial functions of degree $n$.\\
Using Lemma \ref{triviallemmaabouttracespluszero} below it is clear that $\chi_d(F)(g_1,\dots,g_d)=\chi_{d+1}(F)(g_1,\dots,g_d,0)$, hence the various $\chi\colon \frakF(G)\to\Lambda_d(G)$ glue to a homomorphism $\chi\colon \frakF(G)\to\Lambda(G)$ of graded rings.

\begin{Lem}\label{triviallemmaabouttracespluszero}
	Let $F\colon \calG\to \vect\BC$ be a functor and let $g\colon A\to A$ and $B$ be morphisms and objects in $\calG$. Then $\tr F(g\oplus 0\colon A\oplus B\to A\oplus B)=\tr F(g)$
\end{Lem}
\begin{Prf}
	We write $g\oplus 0$ as the composition $A\oplus B\xra{\pi_A} A\xra{g}A\xra{\iota_A} A\oplus B$ and obtain $\tr F(g\oplus 0) = \tr \(F(\iota_A)\circ F(g)\circ F(\pi_A)\)=\tr \(F(g)\circ F(\pi_A)\circ F(\iota_A)\)=\tr F(g)$
\end{Prf}

\begin{Prop}
	Let $G$ be a finite abelian group. Then the ring homomorphisms $\chi\colon \frakF(G)\to\Lambda(G)$ and $\ch\colon R(G)\xra{\cong}\Lambda(G)$ agree under the identification  $R(G)\xra\cong \frakF(G)$ induced by the monoidal equivalence $\Xi\colon (G\wr \Sn\actson E)\mapsto (E\otimes T^n(\blank))^{\Sn}$.\\
	In particular $\chi$ is an isomorphism of graded rings.
\end{Prop}

\begin{Prf}
	The complete symmetric polynomials generate the ring of symmetric polynomials, hence by \ref{chisiso}\ the $E_\gamma^{\otimes n}$ generate the ring $R(G)$. Therefore we see that \ref{chmapsEgammatoh}\ determines the ring homomorphism $\ch$ uniquely which means we only need to show that $\chi$ maps the functor $F_{\gamma^{\otimes n}}\coloneqq \(E_\gamma^{\otimes n}\otimes T^n(\blank)\)^{G\wr \Sn}$ to $h_n(y_{\ast,\gamma} )$.\\
	We compute:
	\begin{align*}
		F_{\gamma^{\otimes n}}(\lreg G\BC^d) &\cong \(\Hom_{G^n}\(\lreg G\BC^d\otimes\dots\otimes\lreg G\BC^d, E^\dual_\gamma\otimes\dots\otimes E^\dual_\gamma\)^\dual\)^{\Sn}\\
		&\cong \(\Hom_G\(\lreg G\BC^d, E_\gamma^\dual\)^\dual\otimes\dots\otimes \Hom_G\(\lreg G\BC^d, E_\gamma^\dual\)^\dual\)^{\Sn}\\
	&\cong \(E_\gamma^d\otimes\dots\otimes E_\gamma^d\)^{\Sn} = \Sym^n(E_\gamma^d)
	\end{align*}
	It is straightforward to track the endomorphism $(\cdot g_i^\inv)_i\colon \lreg G\BC^d\to\lreg G\BC^d$ through the isomorphisms and see that it is the usual left-multiplication on ${E_\gamma^d}^{\otimes n}$ given by $(g_i)_i^{\otimes n}: {(e^{(1)}_i)}_i\otimes\dots\otimes {(e^{(n)}_i)}_i\mapsto {(g_ie^{(1)}_i)}_i\otimes\dots\otimes {(g_ie^{(n)}_i)}_i$.
	The proof is completed by the next Lemma.
\end{Prf}

\begin{Lem}
	Let $G$ be a finite abelian group and $\gamma_1,\dots,\gamma_d\in G^\star $ some (not necessarily distinct) irreducible characters of $G$ with corresponding representations $E_1,\dots,E_d$. Then the character of the $G$-representation $\Sym^n(\bigoplus_{i=1}^dE_i)$ is given by the complete symmetric polynomial $h_n(\gamma_1,\dots,\gamma_d)$.
\end{Lem}

\begin{Prf}
	Since $G$ is abelian, all the $E_i$ are one-dimensional; choose a basis element $e_i$ of $E_i$. Then for each $g\in G$, the element $e_i$ is an eigenvector of $g\colon E_i\to E_i$ with eigenvalue $\gamma_i(g)$. A basis of $\Sym^n(\bigoplus_iE_i)$ is given by monomials $e_{i_1}\cdot\dots\cdot e_{i_n}$ with $1\leqslant i_1\leqslant \dots\leqslant i_n\leqslant d$. This monomials are eigenvectors of $g$ with eigenvalue $\gamma_{i_1}(g)\cdot\dots\cdot\gamma_{i_n}(g)$, hence the trace of $g$ is given by $\sum\limits_{1\leqslant i_1\leqslant \dots\leqslant i_n\leqslant d}\gamma_{i_1}\cdot\dots\cdot\gamma_{i_n}(g)\eqqcolon h_n(\gamma_i,\dots,\gamma_d)(g)$.
	\end{Prf}

\ifx\isstandalone\undefined
\else
\newpage
\bibliographystyle{amsalpha}
\bibliography{literatur}
\end{document}
\fi


\newpage\section{Schur-Weyl duality for wreath products}\label{sectionSchurWeylduality}
	\ifx\preambleloaded\undefined
   
   \def\isstandalone{}
\fi

\subsection{Overview}
We work over the complex numbers $\BC$. Let $G$ be a finite abelian group and let $V\cong\lreg G\BC^d$ be a free $G$-representation of rank $d$.\\
The $n$-fold tensor product $V^\on$ carries canonical actions by the wreath product $G\wr \Sn$ and by the automorphism group $\Aut_G(V)$; clearly those two actions commute with each other.\\
A careful study of this actions allows us to find a reciprocity between irreducible representations of $G\wr \Sn$ and irreducible polynomial representations of $\Aut_G(V)$. More specifically, we will construct for each $n\in \BN$ an assignment (with $0\mapsto 0$)
\begin{equation}\label{SchurWeylassignment}
	\begin{tikzcd}
	{\{\text{irred.\ $G\wr \Sn$-reps (up to iso)}\}}\cup{\{0\}}\dar{\Phi}\\
	{\{\text{irred.\ polyn.\ $\Aut_G(V)$-reps of homog.\ deg.\  $n$ (up to iso)}\}}\cup {\{0\}}
	\end{tikzcd}
\end{equation}
such that the $n$-fold tensor product decomposes as 
\[V^\on\cong \bigoplus\limits_\lambda X_\lambda\otimes \Phi(X_\lambda),\]
where $\lambda$ labels pairwise non-isomorphic irreducible representations $X_\lambda$ of the wreath product $G\wr \Sn$.\\
The representations $\Phi(X_\lambda)$ of $\Aut_G(V)$ are constructed by sending (the dual of) the representation $X_\lambda$ through the equivalence of categories
\[\Xi_n\colon \lrep{G\wr \Sn}\BC\xra\simeq\Pol^n_\BC(\lrep G\BC)\]
from Corollary~\ref{equivalenceforgroupsmain} and then evaluating the resulting functor $\Xi_n(X_\lambda)$ at $V$.

\begin{Goal}[Theorem~\ref{SchurWeylDualityforG}] The assignment $\Phi$ is always injective outside of its kernel (i.e.\ non-isomorphic representations which do not get sent to zero stay non-isomorphic) and $\Phi$ is always surjective. For $n\leqslant d$ the kernel vanishes, hence in this case the assignment is a bijection.
\end{Goal}

The cornerstone of our approach is the fact (established in Section~\ref{sectioncharacteristicmap}) that the monoidal equivalence
\[\Xi\colon \lrep{G\wr \SG\ast}\BC\xra{\simeq} \Pol_\BC^\finitary(\lrep G\BC)\]
categorifies the classical ring isomorphism
\[\ch\colon \K_0(\lrep{G\wr S_\ast}\BC)\xra{\cong} \Lambda(G).\]

\begin{Rem}
The structure and content of this section are a direct adaptation of the corresponding chapter in Macdonald's book~\cite[\S I.A]{MacDonaldSymmetricFunctionsHallPolynomials} where the same methods are used in the case $G=\triv$.\\
Various other roads to a Schur-Weyl duality for wreath products $G\wr \Sn$ are known; see~\cite{MazorchukStropppelGlkdModulesViaGroupoids} and references therein.\end{Rem}

\subsection{Matrix coefficients}
Let $A$ be any (not necessarily finite) group and let $\rho\colon A\to \Aut(W)$ be a representation of $A$ on the $m$-dimensional vector space $W$. If we choose a basis of $W$ then $\rho$ gives a map $A\to \Aut(W)\hra\End(W)\cong \Mat_{m\times m}(\BC)$ or equivalently a collection of $m\times m$ many maps $\rho_{ij}\colon A\to \BC$. We call the map $\rho_{ij}$ the \textbf{matrix coefficient} of $\rho$ for the pair $(i,j)$ (with respect to the given basis).

\begin{Prop}\label{matrixcoeffstotestforisomirred}~\cite[I.A.8.2]{MacDonaldSymmetricFunctionsHallPolynomials}~\cite[(27.8)]{CurtisReinerRepresentationTheoryofFiniteGroupsandAssAlgebras}
	Let $\rho, \rho', \rho'',\dots$ be a bunch of non-zero $A$-representations of dimensions $m, m', m'',\dots$ respectively. Then the following are equivalent:
	\begin{enumerate}
		\item\label{representationspairwisenonisomorphic} The representations $\rho,\rho',\dots$ are irreducible and pairwise not isomorphic.
		\item\label{matrixcoefficientslinearindependent} The set $\{\rho_{ij},\rho'_{i'j'},\ldots\mid 1\leqslant i,j\leqslant m; 1\leqslant i',j'\leqslant m'; \dots\}$ of matrix coefficients is linearly independent as a subset of the vector space $\BC^A$ of functions $A\to \BC$.\qedhere
	\end{enumerate}
\end{Prop}
\begin{Rem}
	Choosing a different basis of a representation will of course change the $\rho_{ij}$'s but it will not change the subspace $\langle\rho\rangle\coloneqq \langle\rho_{ij}\mid 1\leqslant i,j\leqslant m\rangle_\BC\subseteq \BC^A$ generated by the $\rho_{ij}$'s inside the space of all functions $A\to \BC$. Hence Condition~\ref{matrixcoefficientslinearindependent} in Proposition~\ref{matrixcoeffstotestforisomirred} is independent of any choice of basis. This of course also follows a posteriori, since Condition~\ref{representationspairwisenonisomorphic} is basis-free.
\end{Rem}

\begin{Rem}If we have representations $\rho_1,\dots,\rho_n$ of $A$ then
\begin{equation}\label{trivialequationaboutspans}
	\langle\rho_1^{d_1}\oplus\dots\oplus\rho_n^{d_n}\rangle=\langle\rho_1\rangle+\dots+\langle\rho_n\rangle
\end{equation}
as subspaces of the vector space $\BC^A$ because in the correct basis the matrices describing the $A$-action on $\rho_1^{d_1}\oplus\dots\oplus\rho_n^{d_n}$ are in block-form.\end{Rem}

\subsection{Polynomial representations of matrix groups}\label{matrixgroups}

If $A\subseteq \GL_I(\BC)\coloneqq \Aut(\oplus_{i\in I}\BC\cdot e_i)\subset \Mat_{I\times I}(\BC)$ is an affine algebraic  matrix group (where $I$ is a finite set) then we have the ring $\calO(\ol A)\coloneqq \BC[x_{ij}\mid i,j\in I]\subset \BC^A$ of polynomial functions on $A$, where $x_{ij}\colon A\to \BC$ is the $(i,j)$-th coordinate function $(a_{st})_{s,t\in I}\mapsto a_{ij}$.\\
Note that the $x_{ij}$ are usually \emph{not} algebraically independent. In fact we have 
\[\calO(\ol A)\xla[x_{ij}\mapsfrom X_{ij}]{\cong}\rquot{\BC[X_{ij}\mid i,j\in I]}{\frakA},\]
where $\frakA$ is the ideal of polynomials vanishing on $A$ and the $X_{ij}\colon \Mat_{I\times I}(\BC)\to \BC$ are the algebraically independent coordinate functions on $\Mat_{I\times I}(\BC)$.\\
If $A$ is \textbf{homogeneous}, i.e.\ the ideal $\frakA$ can be generated by homogenous polynomials, then we have a well defined decomposition $\calO(\ol A)=\bigoplus_n\calO(\ol A)$ into the subspaces $\calO^n(\ol A)\subset \calO(\ol A)$ consisting of all polynomial functions of homogeneous degree $n$. 
\begin{Def}
We call a representation $\rho\colon A\to\Aut(W)$ of the affine algebraic matrix group $A\subseteq\GL_I(\BC)$ a \textbf{polynomial representation} if all the matrix coefficients $\rho_{ij}\colon A\to \BC$ of $\rho$ are polynomial functions (for one, hence all bases of $W$). In other words we require the subspace $\langle\rho\rangle$ of $\BC^A$ to be already contained in $\calO(\ol A)$.\\
In case $A$ is homogeneous: we say that $\rho$ is \textbf{polynomial of homogeneous degree $n$} if $\langle\rho\rangle$ is already contained in $\calO^n(\ol A)$.
\end{Def}
With the same proof as for polynomial functors one can prove:
\begin{Lem}\label{Lemmairredpolyrepsarehomogeneous}
Every polynomial representation can be decomposed into a sum of homogeneous ones. In particular every irreducible polynomial representation is homogeneous.
\end{Lem}

\subsection{The $n$-standard representation}
Let $A\subseteq \GL_I(\BC)$ be a homogeneous affine algebraic matrix group and fix a natural number $n$. We define the \textbf{$n$-standard representation} of $A$ to be the $n$-fold tensor-power $T^n(\BC^{(I)})=\BC^{(I)}\otimes\dots\otimes \BC^{(I)}$ with $A$ acting diagonally on all components. We want to compute the matrix coefficients of this action.\\
The vector space $T^n(\BC^{(I)})$ has a canonical basis given by monomials $e_{j_1}\otimes\dots\otimes e_{j_n}$ indexed by tuples $(j_1,\dots,j_n)\in I^n$. Let $a=(a_{ij})_{i,j\in I}$ be an element of $A$. We compute
\begin{align*}
	a.(e_{j_1}\otimes\dots\otimes e_{j_n})&=\(\sum_{i_1\in I}a_{i_1j_1}e_{i_1}\)\otimes\dots\otimes \(\sum_{i_n\in I}a_{i_nj_n}e_{i_n}\)\\
	&=\sum_{(i_1,\dots,i_n)\in I^n}\prod_{l=1}^na_{i_lj_l}(e_{i_1}\otimes\dots\otimes e_{i_n})
\end{align*}
which means that the matrix coefficient for the pair $((i_1,\dots,i_n),(j_1,\dots,j_n))$ is given by the monomial $\prod_{l=1}^nx_{i_lj_l}\colon a\mapsto\prod_{l=1}^na_{i_lj_l}$.
Since by definition the monomials $\prod_{l=1}^nx_{i_lj_l}$ span $\calO^n(\ol A)$ we obtain
\begin{equation}\label{monomialsspaneverything}
	\langle T^n(\BC^{(I)})\rangle=\calO^n(\ol A).
\end{equation}

\subsection{The ring of polynomial functions on $A_d(G)$}
We denote by $A_d(G)$ the group $\Aut_G(\lreg G\BC^d)$ of $G$-automorphims of the free $G$-module $\lreg G\BC^d$. Note that $A_d(G)$ is an affine algebraic matrix group in the sense of Section~\ref{matrixgroups} with respect to the standard basis $dG\coloneqq \{1,\dots,d\}\times G$ of $\lreg G\BC^d$.\\
The subvariety $\ol{A_d(G)}=\End_G(\lreg G\BC^d)$ of $\End_\BC(\lreg G\BC^d)=\Mat_{dG\times dG}(\BC)$ is described by the equations $a\circ g=g\circ a$ for $g\in G$. Hence the ideal $\frakA$ of $\BC[T_{if,jh}\mid 1\leqslant i,j\leqslant d; g,h\in G]$ which describes this subvariety is generated by the (homogeneous) polynomials $T_{if,jh}-T_{i(fg),j(gh)}$ (for $f,g,h\in G$) as can be seen by an easy computation.\\
To proceed we need an easy lemma about the abelianization $G^\ab\coloneqq \rquot{G}{[G,G]}$ of the group $G$.

\begin{Lem}\label{abelianizationrelation}
	Let $\sim$ be the equivalence relation on the set $G\times G$ generated by $(f,h)\sim (fg,gh)$ for all $f,g,h\in G$. Then the map
	\[q\colon G\times G\xra{(f,h)\mapsto f^{-1}h} G\thra G^\ab\]
	induces a bijection $\ol q\colon \rquot{G\times G}{{\sim}}\xra{\cong} G^\ab$.
\end{Lem}

\begin{Prf}
	For any $f,g,h\in G$ we have $q(fg,gh)=g^{-1}f^{-1}gh=f^{-1}h=q(f,h)$ since $g^{-1}f^{-1}g=f^{-1}$ in $G^\ab$, hence $\ol q$ is well defined. For elements $g_1,g_2,g\in G$ we have $(1_G,g_1g_2g)\sim (g_1^{-1},g_2g)\sim(g_1^{-1}g_2^{-1},g)=((g_2g_1)^{-1},g)\sim (1_G,g_2g_1g)$, hence the map
	\[G^\ab\xra{g\mapsto(1,g)} \rquot{G\times G}{{\sim}}\]
	is a well defined inverse to $\ol q$.
\end{Prf}

Using Lemma~\ref{abelianizationrelation} we deduce that $\calO(\ol{A_d(G)})$ is the polynomial ring $\BC[x_{ijg}\mid 1\leqslant i,j\leqslant d; g\in G^\ab]$ in the $d^2|G^\ab|$-many algebraically independent variables $x_{ijg}$ which we view as the coordinate functions $x_{if,jh}$ for any pair $(f,h)\in G\times G$ such that $f^{-1}h\equiv g \mod [G,G]$.\\
This implies in particular that the space $\calO^n(\ol{A_d(G)})$ of homogeneous polynomial functions of degree $n$ on $A_d(G)$ has dimension
\begin{equation}\label{dimensionofhomopolyring}
	\dim \calO^n(\ol{A_d(G)}) = \multiset{d^2|G^\ab|}{n},
\end{equation}
where $\multiset mn$ is the number of multisets of cardinality $n$ contained in a given set of size $m$ and counts the number of monomials of degree $n$ in $m$ variables.

\subsection{Decomposing the $n$-standard representation}
Denote by $T(n,d)$ the $n$-fold tensor product $\(\lreg G\BC^d\)^\on$ which is a representation both of $G\wr \Sn$ and of $A_d(G)$ (and those actions commute).\\
Recall that the partition valued maps $\lambda\in P_n(G^\star)$ of total size $n$ index the irreducible representations $X_\lambda$ of $G\wr \Sn$. Hence by Corollary~\ref{equivalenceforgroupsmain} they also index irreducible polynomial functors of homogenous degree $n$ via the formula $F_\lambda\coloneqq \Xi_n(X_\lambda)\cong\Hom_{G\wr \Sn}(X^\dual_\lambda,T^n(\blank))$.\\
We obtain representations $R_\lambda^d$ of the group $A_d(G)$ by evaluating the functor $F_\lambda$ at $V=\lreg G\BC ^d$; more precisely, $R_\lambda^d$ is given as the composition 
\[R^d_\lambda\colon \B A_d(G)=\B\Aut_G(V)\hra\lrep G\BC\xra{F_\lambda}\vect\BC.\]
This representations are clearly polynomial of homogeneous degree $n$.\\
Note that we have an isomorphism of $G\wr \Sn$-representations
\[T^n(X)\cong \bigoplus\nolimits _\lambda X^\dual_\lambda\otimes \Hom_{G\wr Sn}(X^\dual_\lambda,T^n(X))=\bigoplus\nolimits_\lambda X^\dual_\lambda\otimes F_\lambda(X),\]
naturally in the $G$-representation $X$. Hence by specializing to $X=\lreg G\BC^d$ we obtain the decomposition 
\begin{equation}\label{SchurWeyldecomposition}
T(n,d)\cong \bigoplus_{\lambda\in \calP_n(G^\star )}X^\dual_\lambda\otimes R^d_\lambda
\end{equation}
as $\Sn$-modules and as $A_d(G)$-modules.

\begin{Expl}
	For $G=\triv$ and $n=2$ we obtain the decomosition
	\[T^2(\BC^d)\cong {\BC^\dual_{-1}\otimes\Lambda^2(\BC^d)}\oplus {\BC_{+1}^\dual\otimes\Sym^2(\BC^d)},\]
	where $\BC_{-1}$ resp. $\BC_{+1}$ are the sign representation resp. the trivial representation of $S_2$.\\
	This example also shows what can go wrong if we don't have $n\leqslant d$: if specialize to $d=1$ then $\Lambda^2(\BC)$ will vanish in the above identity, hence the decomposition will degenerate. In other words this is the case where some of the $R_\lambda^d$'s might be zero.
\end{Expl}

We will see soon that this degeneration phenomenon cannot appear for $n\leqslant d$.

\subsection{Schur-Weyl duality for wreath products}

\begin{Thm}[Schur-Weyl duality for $G\wr S_\ast$]\label{SchurWeylDualityforG}
Assume $G$ is an abelian group and let $n,d$ be natural numbers.
\begin{enumerate}[label=(\arabic*), ref=(\arabic*)]
	\item\label{SWinjwherenotzero} Some of the $A_d(G)$-representations $R_\lambda^d$ (for $\lambda\in \calP_n(G^\star)$) might be zero. All the others are irreducible and pairwise non-isomorphic polynomial representations of homogeneous degree $n$.
	\item\label{SWinjfornleqd} For $n\leqslant d$ all the $R_\lambda^d$ (for $\lambda\in \calP_n(G^\star)$) are irreducible and pairwise non-isomorphic polynomial representations of $A_d(G)$ of homogeneous degree $n$.
	\item\label{SWalwayssurj} Every irreducible polynomial representation of $A_d(G)$ is homogeneous of some degree $n$ and appears as $R_\lambda^d$ for some $\lambda\in \calP_n(G^\star)$.\qedhere
\end{enumerate}
\end{Thm}

\begin{Rem}
Note that Theorem~\ref{SchurWeylDualityforG} is basically saying that the assignment $X_\lambda\mapsto R_\lambda^d$ (and $0\mapsto 0$) gives a well defined mapping
\[	\begin{tikzcd}
	{\{\text{irred.\ } G\wr \Sn\text{-reps (up to iso)}\}}\cup{\{0\}}\dar\\
	{\{\text{irred.\ polyn. } A_d(G)\text{-reps of hom.\ degree } n \text{ (up to iso)}\}}\cup {\{0\}}
	\end{tikzcd}
\]
which is always surjective, is always injective outside of its kernel and is a bijection if $n\leqslant d$.
\end{Rem}

Before we can start with the proof of the theorem we need a couple of computations.

\begin{Lem}\label{technicallemmaaboutSchur}
For a finite set $X$ we have the identity
\[\sum_{\lambda\in \calP_n(X)}\(\prod_{x\in X}s_{\lambda(\gamma)}(1^d)\)^2=\multiset{d^2|X|}n\]
\end{Lem}
\begin{Prf}
We start with the following identity of formal power series in the variables $x_1,x_2,\dots$ and $y_1,y_2,\dots$~\cite[I.4.3]{MacDonaldSymmetricFunctionsHallPolynomials}:
\[\prod_{i,j}(1-x_iy_j)^{-1}=\sum_{\lambda\in \calP}s_\lambda(x_1,\dots)s_\lambda(y_1,\dots)\]
By setting $x_i=y_i=t$ for $1\leqslant i\leqslant d$ and $x_i=y_i=0$ for $i>d$ we obtain the identity
\[(1-t^2)^{-d^2}=\sum_{\lambda\in \calP}s_\lambda(t,\dots,t)^2\]
which we can take to the $|X|$-th power to get
\[(1-t^2)^{-d^2|X|}=\sum_{\lambda\colon X\to \calP}\prod_{x\in X}s_{\lambda(x)}(t,\dots,t)^2\]
By restricting to the monomial $t^{2n}$ on both sides we get
\[ct^{2n}=\sum_{\lambda\in \calP_n(X)}\prod_{x\in X}s_{\lambda(x)}(t,\dots,t)^2,\]
where $c$ is the coefficient of $t^{2n}$ in $(1-t^2)^{-d^2|X|}$, which is the same as the coefficient of $t^n$ in $(1-t)^{-d^2|X|}=(1+t+t^2+\dots)^{d^2|X|}$.
It is an elementary calculation that this coefficient is precisely $\multiset{d^2|X|}{n}$, hence specializing to $t=1$ finishes the proof.
\end{Prf}

\begin{Lem}\label{dimofRlambdad}
Let $G$ be an abelian group, $n,d\in \BN$ and $\lambda\in \calP_n(G^\ast)$. Then the degree of the representation $R_\lambda^d$ is
	\[\dim R_\lambda^d = \prod_{\gamma\in G^\ast}s_{\lambda(\gamma)}(1^d).\qedhere\]
\end{Lem}
\begin{Prf}
	The degree of the representation $R^d_\lambda$ can be computed as
	\[\dim F_\lambda(\lreg G\BC^d) = \tr F_\lambda(\Id_{\lreg G\BC^d})=\tr F_\lambda((\cdot {1_G}^{-1})_i)=\chi(F_\lambda)(1_G,\dots,1_G),\] where $1_G$ is the unit of the group $G$.
	By definition $F_\lambda$ corresponds to  $X_\lambda$ under the isomorphism $\frakF(G)\cong R(G)$ and we saw in Section~\ref{subsectioncharacteristicmap} that $\chi$ corresponds to $\ch$, hence we can replace $\chi(F_\lambda)$ by $\ch(X_\lambda)$. Moreover, $\ch(X_\lambda)(1_G,\dots,1_G)$ is nothing but the Schur function $S_\lambda(y_{\ast,\ast })(1_G,\dots,1_G)\coloneqq \prod_{\gamma\in G^\star}s_{\lambda(\gamma)}(\gamma(1_G),\dots,\gamma(1_G))$ by the Property~\ref{chmapsXtoSchur}. Since the group $G$ is abelian we have $\gamma(1_G)=1$ for all irreducible characters $\gamma\in G^\star$ and the result follows.
\end{Prf}

\begin{Rem}\label{nleqdrangeremark}
The Schur function $s_{\lambda(\gamma)}$ can be defined by the sum $s_{\lambda(\gamma)}(x_1,\dots)\coloneqq \sum_T x^T$ where $T$ runs over the semistandard tableaux of shape $\lambda(\gamma)$~\cite[Definition 4.4.1]{SaganTheSymmetricGroup}. Hence is clear that the expression $s_{\lambda(\gamma)}(1^d)$ counts semistandard Young tableaux with entries contained in $\{1,\dots, d\}$. This number is always positive for $|\lambda(\gamma)|\leqslant d$, hence $R_\lambda^d$ is non-zero in the range $n\leqslant d$.
\end{Rem}

After these preparations we can finally prove the Schur-Weyl duality for wreath products.

\begin{Prf}[\Proofof{Theorem~\ref{SchurWeylDualityforG}}]
	 Putting together Lemma~\ref{dimofRlambdad}, Lemma~\ref{technicallemmaaboutSchur} (applied to $X=G^\star$) and Equation~\ref{dimensionofhomopolyring} we obtain (observe $|G^\star|=|G|=|G^\ab|$ since $G$ is abelian) 
	\begin{equation}\label{numberofmatrixcoeffs}
		\sum_{\lambda\in \calP_n(G^\star)}\dim^2 R_\lambda^d=\sum_{\lambda\in \calP_n(G^\ast )}\(\prod_{x\in G^\ast }s_{\lambda(\gamma)}(1^d)\)^2=\multiset{d^2|G^\star|}n=\dim \calO^n(\ol {A_d(G)}).
	\end{equation}
	Equation~\ref{SchurWeyldecomposition} is a decomposition of $T(n,d)$ into copies of $R_\lambda^d$, hence by Equation~\ref{monomialsspaneverything} and Equation~\ref{trivialequationaboutspans} we have 
	\begin{equation}\label{Rlambdasspaneverything}
		\calO^n(\ol {A_d(G)})=\langle T(n,d)\rangle=\sum_{\lambda\in\calP_n(G^\star)}\langle R_\lambda^d\rangle,
	\end{equation} i.e.\ the matrix coefficients of the $R_\lambda^d$ jointly span $\calO^n(\ol {A_d(G)})$. The total number of this matrix coefficients is described in Equation~\ref{numberofmatrixcoeffs}, hence the matrix coefficients of all the $R_\lambda^d$ are actually a basis of $\calO^n(\ol {A_d(G)})$ by reasons of dimensions.\\
	It follows from Proposition~\ref{matrixcoeffstotestforisomirred} that all the non-zero $R_\lambda^d$'s are irreducible and non-isomorphic representations of $A_d(G)$, hence part~\ref{SWinjwherenotzero} is proved.\\
	Conversely, let $\rho$ be some irreducible polynomial representation of $A_d$. The representation $\rho$ must be homogeneous of some degree $n$ by Lemma~\ref{Lemmairredpolyrepsarehomogeneous}, hence we obtain $\langle\rho\rangle\subseteq \calO^n(\ol {A_d(G)})=\sum_{\lambda\in\calP_n(G^\star)}\langle R_\lambda^d\rangle$ by Equation~\ref{Rlambdasspaneverything}. Therefore by Proposition~\ref{matrixcoeffstotestforisomirred} the representation $\rho$ must be isomorphic to one of the $R_\lambda^d$'s; this proves part~\ref{SWalwayssurj}.\\
	Part~\ref{SWinjfornleqd} now follows directly from Remark~\ref{nleqdrangeremark} which guarantees that $R_\lambda^d$ is non-zero for $n\leqslant d$.
\end{Prf}

\ifx\isstandalone\undefined
\else
\newpage
\bibliographystyle{amsalpha}
\bibliography{literatur}
\end{document}
\fi


\newpage\begin{appendices}
\section{Abstract Nonsense}

\subsection{The Grothendieck construction}\label{AppendixGrothendieckConstruction}
\ifx\preambleloaded\undefined
   
   \def\isstandalone{}
\fi

When trying to define a category valued functor $\calB\to\CAT$ one often faces the problem of having to make some arbitrary choices.\\ The idea of the Grothendieck construction is that one can view a (pseudo-)functor $\calB\to\CAT$ as a \roughly{fibration} over $\calB$ and that this \roughly{fibration} $?\to\calB$ can often be defined by making no choice at all, or rather, by taking \emph{all} possible choices.
\begin{Expl}
Consider the {(pseudo-)functor}
\[\lMod?\colon  \category{Ring}\lra \CAT\]
which assigns to each ring $A$ the category of $A$-modules and to each ring homomorphism $f\colon A\to B$ the functor $B\otimes_A\blank$. If one is really pedantic then this means that we \emph{choose} a specific model of the tensor product $B\otimes_A M$ for each $A$-module $M$ amongst all the possible models which satisfy the universal property.\\
What we could do instead is to define the category $\calM$ of pairs $(A,M)$ where $A$ is a ring and $M$ is an $A$-module. A morphism $(A,M)\to(B,N)$ would then be a ring homomorphism $f\colon A\to B$ together with a $A$-module homomorphism $M\to N$ (where $A$ acts on $N$ via $f$). The category $\calM$ admits an obvious functor $\calM\to\category{Ring}$ which just forgets the second component.\\
Note that this whole construction is completely choice-free.\\
We can recover the category $\lMod A$ as the fiber $\calM_A$ over $A$. Moreover if $(A,M)$ is in the fiber over $A$ and $f\colon A\to B$ is a ring homomorphism we can say that a model for $B\otimes_A M$ is any $T\in \calM_B$ together with a natural isomorphism $\Hom_\calM(T,\blank )\cong \Hom_\calM((A,M),\blank )$ of functors $\calM_B\to\Set$.\\
We can thus recover the original functor $\lMod ?$ by making lots of choices of such models $T=(B,B\otimes _A M)$.
\end{Expl}
Keeping this example in mind we can give a general construction.

\begin{Cstr}Let $\Gamma\colon \calB\to\CAT$ be a pseudo-functor. The \introduce{Grothendieck construction} $\Groth\Gamma\calB$ is the following category:
\begin{enumerate}
	\item Objects are pairs $(c,b)$, where $b$ is an object of $\calB$ and $c$ is an object of the category $\Gamma_b$.
	\item A morphism $(c,b)\to(c',b')$ is a pair $(\alpha\colon \Gamma_f(c)\to c',f\colon b\to b')$, where $f$ is a morphism in $\calB$ and $\alpha$ is a morphism in $\Gamma_{b'}$.
\end{enumerate}
The identity on the object $(c,b)$ is $(\Gamma_{\Id_b}(c)\xra{\cong}c,\Id_b)$, the composition of composable arrows $(c,b)\xra{(\alpha,f)}(c',b')\xra{(\beta,g)}(c'',b'')$ is given by
\[(\Gamma_{gf}(c)\cong \Gamma_g(\Gamma_f(c))\xra{\Gamma_g(\alpha)}\Gamma_g(c')\xra{\beta}c'', g\circ f)\qedhere\]
\end{Cstr}
The category $\Groth\Gamma\calB$ comes equipped with the obvious projection functor $p\colon \Groth\Gamma\calB\to\calB$.\\
We want to classify the functors $p\colon \calD\to\calB$ that arise this way. For any object $b\in \calB$ we denote by $\calD_b$ the fiber of $p$ over $\Id_b$, i.e.\ the subcategory $\calD_b\subset\calD$ of morphisms $\psi$ such that $p(\psi)=\Id_b$.

\begin{Def}\label{Defopfibration}
	A functor $p\colon \calD\to\calB$ is called an \textbf{op-fibration of categories (over $\calB$)}, if the following condition is satisfied:
	\begin{enumerate}[label=\axiomlabelstyle{M}, ref=\axiomrefstyle{M}]
		\item \label{opfibration} 
Let $f\colon b\to b'$ be a map in $\calB$ and $F\in \calD_b$ an object in the fiber. Then we require the existence of an object $G\in\calD_{b'}$ and a morphism $\varphi\colon F\to G$ lifting $f$ (i.e.\ $p(\varphi)=f$) such that for any object $H\in \calD$ the square
		\begin{equation}\label{inducedfunctoronfibers}
		\begin{tikzcd}[column sep=50]
			\Mor\calD{G}{H}\rar{\Mor\calD{\varphi}{\Id}}\dar{p}&	\Mor\calD{F}{H}\dar{p}\\
			\Mor{\calB}{b'}{p(H)}\rar{\Mor{\calD}{f}{p(H)}}&	\Mor{\calB}{b}{p(H)}
		\end{tikzcd}
		\end{equation}
is a pullback of sets.\qedhere
	\end{enumerate}
\end{Def}

\begin{Prop}\label{Grothendieckisopfibration}
	The projection functor $p\colon \Groth\Gamma\calB\to\calB$ is an op-fibration of categories.
\end{Prop}
\begin{Prf}Abbreviate $\calD\coloneqq \Groth\Gamma\calB$.
	Let a morphism $f\colon b\to b'$ in $\calB$ and an object $F=(c,b)\in \calD_b$ be given. I claim that the morphism $\varphi\coloneqq(\Id_{\Gamma_f(c)},f)\colon (c,b)\to(\Gamma_f(c),b')\eqqcolon G$ is the desired lift of $f$.\\
	Indeed if for any other object $H\coloneqq(\hat c,\hat b)$ we are given a morphism $g\colon b'\to\hat b=p(H)$ and a morphism $\psi=(\alpha\colon \Gamma_h(c)\to\hat c,h\colon b\to \hat b)\colon F\to H$ such that $h=p(\psi)=g\circ f$ then $\hat\psi\coloneqq(\alpha\colon \Gamma_g(\Gamma_f(c))\to\hat c, g)$ is the unique lift of $g$ such that $\psi=\hat\psi\circ \varphi$. In other words, elements $g\in \Mor\calB{b'}{p(H)}$ and $\psi\in \Mor\calD F H$ with $g\circ f=p(\psi)$ glue uniquely to the element $\hat\psi\in \Mor\calD G H$, i.e.\ the square \ref{inducedfunctoronfibers} is a pullback.
\end{Prf}

Now we want to study the converse. Fix an op-fibration $p\colon \calD\to\calB$ of categories. Observe that \ref{opfibration} is saying in particular (by restricting the pullback to the fiber of the element $\Id_{b'}\in \Mor\calB {b'}{b'}$)  that for all $F\in\calD_{b}$ the functor 
\begin{equation}
f_\star(F)=\Mor\calD{F}{\blank}_f=\{\varphi\in\Mor\calD{F}{\blank}\mid p(\varphi)=f\}\colon\calD_{b'}\lra \Set
\end{equation}
is represented by the object $G\in\calD_{b'}$ (via the universal object $\varphi\in f_\star(F)(G)$). The assignment $F\mapsto f_\star(F)$ is contravariantly functorial in $F\in\calD_{b}$, i.e.\ we have a functor $f_\star\colon \calD_{b}\to[\calD_{b'}, Set]^\op$. Saying that all the $f_\star(F)$ are representable means that $f_\star$ factors through the Yoneda embedding $\calD_{b'}\hra[\calD_{b'}, Set]^\op$, i.e.\ we get an induced functor $f_\star\colon\calD_{b}\ra\calD_{b'}$.\\
Moreover, for composable morphisms ${b}\xra{f}{b'}\xra{g}{b''}$ in $\calB$ we can restrict the pullback in \ref{opfibration} to the fiber of the element $g\in \Mor\calB{b'}{b''}$ and obtain that
\[\Mor\calD\varphi\Id\colon \Mor\calD{f_\star(F)}\blank_g\xra{\cong} \Mor\calD F\blank_{g\circ f}\]
is an isomorphism of functors $\calD_{b''}\to\Set$ (which is also natural in $F$). Thus we get a canonical isomorphism of the corresponding representing objects $g_\star(f_\star(F))\cong (g\circ f)_\star(F)$ which is also natural in $F$.\\
Therefore we obtain a pseudofunctor $\calD\colon \calB\to\CAT$ given by ${b}\mapsto\calD_{b}$ on objects and by $(f\colon b\ra b')\mapsto (f_\star\colon \calD_{b}\to\calD_{b'})$ on morphisms.\\
From the proof of Proposition~\ref{Grothendieckisopfibration} we know that if $p$ was of the shape $\calD=\Groth\Gamma\calB\to \B$ then for $F=(c,b)$ and $f\colon b\to b'$ as above a representing object of the functor $\Mor \calD F\blank_f$ can be chosen to be $(\Gamma_f(c),b')$. Hence if we identify $\Gamma(b)$ with $\calD_b$ by mapping $\alpha\colon c\to c'$ to $(\Gamma_{\Id_b}(c)\cong c\xra{\alpha}c',\Id_b)\colon (c,b)\to(c',b)$ then the functor $\calD\colon \calB\to\CAT$ becomes nothing but $\Gamma$ itself.
\begin{dCor}
We have a one-to-one correspondence
\[\{\text{pseudo-functors } \calB\to\CAT\}\xra{\blank\rtimes\calB}\{\text{op-fibrations over }\calB\}\qedhere\]
\end{dCor}
\begin{Rem}
Both sides of this correspondence are in fact $2$-categories and the Grothendieck construction can be upgraded to an actual equivalence of $2$-categories. See the literature~\cite[Theorem 1.3.6]{JohnstoneSketchesofanElephant} for more details.
\end{Rem}

\ifx\isstandalone\undefined
\else
\newpage
\bibliographystyle{amsalpha}
\bibliography{literatur}
\end{document}
\fi

\subsection{Semidirect products of categories}\label{Appendixsemidirectproduct}
\ifx\preambleloaded\undefined
   
   \def\isstandalone{}
\fi

Let $S$ be a monoid with unit $\unit\in S$.
\begin{Def} An \introduce{action} of S on a category $\calC$ (written $\Gamma\colon S\actson \cal C$) is a functor $\Gamma\colon \B S\to\CAT$ with $\cal C=\Gamma(\star)$.
\end{Def}
Fix such an action $\Gamma\colon S\actson \cal C$. We define the \textbf{semidirect product} $C\rtimes_\Gamma S$ of $\cal C$ with $S$ (along $\Gamma$) by applying the Grothendieck construction (see Appendix~\ref{AppendixGrothendieckConstruction}) to the functor $\Gamma\colon \B S\to\CAT$. Explicitly, $\cal C\rtimes_\Gamma S$ has an object $(c,\star)$ for each object $c\in \cal C$ and a morphism between two objects $(c,\star),(c',\star)\in \cal C\rtimes_\Gamma S$ is a pair $(f,s)$, where $s\in S$ and $f\colon \Gamma_s(c)\to c'$ is a morphism in $\cal C$. Composition is given by 
\[(g\colon \Gamma_t(c')\to c'',t)\circ (f\colon \Gamma_s(c)\to c',s)\coloneqq (g\circ \Gamma_t(f)\colon \Gamma_{ts}(c)\to c'',ts)\]
Note that this generalizes the usual semidirect product of groups in the sense that $\B(G\rtimes_\Gamma S)=(\B G)\rtimes_\Gamma S$ when $S$ is a group and acts on a group $G$ (or, equivalently, on the category $\B G$) via $\Gamma$.

\begin{Rem}
	It is immediate from the definition that if the action $\Gamma\colon S\actson \cal C$ is trivial then we recover the usual product of categories in the sense that $\cal C\rtimes_\Gamma S=\cal C\times \B S$.
\end{Rem}

\subsubsection{Universal property}
The category $\cal C\rtimes_\Gamma S$ comes equipped with some extra structure:
\begin{itemize}
	\item The functor $(-,\star)\colon \cal C\hra{\cal C\rtimes_\Gamma S}$ given by $c\mapsto (c,\star)$ on objects and by $(f\colon c\to c')\mapsto(f\colon \Gamma_\unit (c)=c\to c',\unit)\colon (c,\star)\to(c',\star)$.
	\item For each $s\in S$ the natural transformation $P_s\coloneqq(\Id_{\Gamma_s}, s)\colon (-,\star)\Rightarrow (-,\star)\circ \Gamma_s$ given on the object $c\in \cal C$ by the morphism $(\Id_{\Gamma_s(c)},s)\colon (c,\star)\to (\Gamma_s(c),\star)$.
\end{itemize}
satisfying $P_t\Gamma_s\circ P_s=P_{ts}$ for all $s,t\in S$.

\begin{Prop}\label{univpropofsemidirectproduct}
	The semidirect product $\cal C\rtimes_\Gamma S$ together with the inclusion $(-,\star)\colon \cal C\hra{\cal C\rtimes_\Gamma S}$ and the natural transformations $\{P_s\}_{s\in S}$ satisfy the following universal property:
	\begin{enumerate}
	\item Let $\cal D$ be a category equipped with a functor $F\colon \cal C\to\cal D$ and a set $\{Q_s\colon F\Rightarrow F\circ \Gamma_s\}_{s\in S}$ of natural transformations such that $Q_t\Gamma_s\circ Q_s=Q_{ts}$ for all $s,t\in S$.
	Then there is a unique functor $\hat F=F\rtimes_\Gamma Q\colon \cal C\rtimes_\Gamma S\to\cal D$ such that $\hat F\circ (-,\star)=F$ and $\hat FP_s=Q_s$.
	\item Let $\cal D$ be equipped with a second such structure $(F'\colon \cal C\to\cal D,\{Q'_s\colon F'\Rightarrow F'\circ \Gamma_s\}_{s\in S})$ and let $\alpha\colon F\Rightarrow F'$ be a natural transformation. Then $\alpha$ is a natural transformation $F\rtimes_\Gamma Q\Rightarrow F'\rtimes_\Gamma Q'\colon \cal C\rtimes_\Gamma S\to \cal D$ (i.e.\ $\alpha$ is natural also with respect to the extra morphisms of the category $\cal C\rtimes_\Gamma S$) if and only if $\alpha$ is compatible with the $Q_s$ and $Q'_s$ in the sense that $Q'_s\circ \alpha=\alpha\Gamma_s\circ Q_s$ for all $s\in S$.\qedhere
	\end{enumerate}
\end{Prop}

\begin{Rem}
This proposition generalizes the well known statement for representations of semidirect products that a representation $\rho\colon \B G\rtimes S\to \vect k$ is the same thing as a representation $\rho\colon \B G\to \vect k$ together with $G$-maps $\rho_s\colon \rho\to\Gamma_s^\star(\rho)$ (i.e.\ natural transformations $\rho\to\rho\circ \Gamma_s$) satisfying $\Gamma_s^\star\rho_t\circ \rho_s=\rho_{ts}$.
\end{Rem}

\begin{Prf}[\Proofof{Proposition~\ref{univpropofsemidirectproduct}}]\begin{enumerate}
		\item We unravel the defining equations of $\hat F$: The first says that $\hat F(c,\star)=F(c)$ and $\hat F(f,\unit)=F(f)$ for objects $c$ and morphisms $f$ in $\cal C$. The second says that for each object $c\in C$ we have $\hat F(\Id_{\Gamma_s(c)},s)\eqqcolon\hat F(P_s(c))=Q_s(c)$. Each object of $\cal C\rtimes_\Gamma S$ is of the shape $(c,\star)$ and each morphism $(f\colon \Gamma_s(c)\to c',s)\colon (c,\star)\to(c',\star)$ of $\cal C\rtimes_\Gamma S$ can be written as the composition $(c,\star)\xra{(\Id_{\Gamma_s(c)},s)}(\Gamma_s(c),\star)\xra{(f,\unit)}(c',\star)$. Hence $\hat F$ is uniquely determined.\\
		To show existence we define $\hat F$ by the above formula, i.e.\ we define $\hat F(c,\star)=F(c)$ for $c\in \cal C$ and $\hat F(f\colon \Gamma_s(c)\to c',s)\coloneqq F(f)Q_s(c)$. If $(g\colon \Gamma_t(c')\to c'',t)$ and $(f\colon \Gamma_s(c)\to c',s)$ are composable  with composite $(g\circ \Gamma_t(f)\colon \Gamma_{ts}(c)\to c'',ts)$ we have to see
		\[F(g)Q_t(c')F(f)Q_s(c)=\hat F(g,t)\hat F(f,s)=\hat F((g,t)(f,s))=F(g)F(\Gamma_t(f))Q_{ts}(c)\] which is immediate from the properties of $\{Q_s\}_{s\in S}$.
		\item Because of the factorization $(f,s)=(f,1)(\Id_{\Gamma_s},s)$ we see that for $\alpha$ to be a natural transformation of functors $\cal C\rtimes_\Gamma S\to D$ is the same as being compatible with the morphisms of the shape $(\Id_{\Gamma_s(c)},s)\colon (c,\star)\to(\Gamma_s(c),\star)$ for $c\in \cal C$. This compatibility is expressed by the equation $\hat{F'}(\Id_{\Gamma_s(c)},s)\alpha(c)=\alpha(\Gamma_s(c))\hat F(\Id_{\Gamma_s(c)},s)$ which translates to $Q'_s\circ \alpha=\alpha\Gamma_s\circ Q_s$ using the definition of $\hat F$ and $\hat{F'}$.\qedhere
	\end{enumerate}
\end{Prf}
We identify the object $c\in \cal C$ with the object $(c,\star)\in \cal C\rtimes_\Gamma S$ and view the functor $(-,\star)\colon \cal C\hra\cal C\rtimes_\Gamma S$ as an inclusion of categories.

\subsubsection{Functoriality}
Let $\theta\colon S\to S'$ be a monoid homomorphism, let $F\colon \cal C\to\cal C'$ be a functor and let $S$ and $S'$ act on $\cal C$ and $\cal C'$ via $\Gamma\colon S\actson \cal C$ and $\Gamma'\colon S'\actson\cal C'$ respectively. Assume further that $F$ has an $\Gamma$-$\Gamma'$-equivariant structure, i.e. isomorphisms $F\circ \Gamma_s\cong \Gamma'_{\theta(s)}\circ F$ for all $s\in S$ (satisfying some coherence).\\
Then we obtain a functor $F\rtimes_{\Gamma,\Gamma'}\theta\colon \cal C\rtimes_\Gamma S\to \cal C'\rtimes_{\Gamma'} S'$ given by $F$ on objects and $(f\colon \Gamma_s(c)\to c',s)\mapsto \(\Gamma'_{\theta(s)}(F(c))\cong F(\Gamma_s(c))\xra{F(f)} F(c'),\theta(s)\)$ on morphisms. It is straightforward to check that this construction is functorial in the pair $(F,\theta)$ up to pseudo-natural isomorphisms.

\subsubsection{Enriched semidirect product}
Let $(\calV,\otimes)$ be a complete and cocomplete closed symmetric monoidal category and let $\CAT_\calV$ denote the (not enriched) $2$-category of $\calV$-categories, $\calV$-functors and $\calV$-transformations.
\begin{Def}
A functor $\Gamma\colon \B S \to \CAT_\calV$ is called a $\calV$-action of $S$ on the $\calV$-category $\calC\coloneqq \Gamma(\star)$.
\end{Def}
We can define an enriched version of the semidirect product
\begin{Cstr}
The $\calV$-category $\calC\rtimes_\Gamma S$ has the same objects as $\calC$ and we define
\[\Mor{(\calC\rtimes_\Gamma S)}{c}{c'}\coloneqq \coprod_{s\in S}\Mor\calC{\Gamma_s(c)}{c'}\]
with composiiton given by
\[\Mor\calC{\Gamma_s(c)}{c'}\otimes\Mor\calC{\Gamma_t(c')}{c'}\xra{\Gamma\otimes\Id}\Mor\calC{\Gamma_{ts}(c)}{\Gamma_t(c')}\otimes\Mor\calC{\Gamma_t(c')}{c''}\xra{\circ}\Mor\calC{\Gamma_{ts}}{c''}\qedhere\]
\end{Cstr}
\begin{Rem}The enriched semidirect product satisfies an universal property which is the $\calV$-enriched analogue of Proposition~\ref{univpropofsemidirectproduct} and characterizes $\calV$-functors $\calC\rtimes_\Gamma S\to\calE$ out of the semidirect product (and $\calV$-natural transformations between them) as $\calV$-functors $F\colon \calC\to\calE$ together with a set $\{Q_s\colon F\Ra F\circ \Gamma_s\}$ of $\calV$-natural transformations satisfying the compatibility conditions $Q_t\Gamma_s\circ Q_s=Q_ts$ for all $s,t\in S$.

\end{Rem}

We can deduce the following useful Lemma.
\begin{Lem}\label{restricttosemidirectprod}
	Let $\calE$ be a $\calV$-category. Let $S$ be a group and let $\calC$ and $\calC'$ be two $\calV$-categories with $\calV$-linear actions $\Gamma\colon S\actson \calC$ and $\Gamma'\colon S\actson\calC'$. Let $\calF$ and $\calF'$ be full $\calV$-subcategories of $[\calC,\calE]$ and $[\calC',\calE]$ respectively and let $\Phi\colon \calF\to\calF'$ be a $\calV$-equivalence.\\
	Assume that $\Phi$ is $S$-equivariant in the sense that for all $s\in S$ the functors $\Gamma_s^\star\colon [\calC,\calE]\to[\calC,\calE]$ and $\Gamma_s'^\star\colon [\calC',\calE]\to[\calC',\calE]$ restrict to the subcategories $\calF$ and $\calF'$ respectively and correspond to each other under the equivalence $\Phi$.\\
	Denote by $\calF_S$ the full subcategoriy of $[\calC\rtimes_\Gamma S,\calE]$ consisting of those $\calV$-functors $F\colon \calC\rtimes_\Gamma S\to \calE$ whose restriction $\calC\hra\calC\rtimes_\Gamma S\xra{F}\calE$ lies in $\calF$; similarly define $\calF'_S$.\\
	Then we can construct a certain $\calV$-equivalence of categories $\Phi_S$ (see the proof) making the following diagram commute:
	\begin{equation}\begin{tikzcd}
		\calF\ar[r,"\Phi"',"\simeq"]&	\calF'\\
		\calF_S\ar[r,"\simeq"',"\Phi_S", dashed]\uar{\Res^{\calC\rtimes S}_\calC}&	\calF'_S\ar[u,"\Res^{\calC\rtimes S}_\calC"']
	\end{tikzcd}\end{equation}	
	Moreover, this construction is functorial in $\Phi$ in the sense that $(\Phi^1\circ \Phi^2)_S\cong\Phi^1_S\circ \Phi^2_S$ for two composable equivalences $\Phi^1$ and $\Phi^2$ as above.
\end{Lem}
\begin{Prf}
By the universal property of the semidirect product, a $\calV$-functor in $\calF_S$ is the same thing as a $\calV$-functor $F\in\calF$ together with a coherent $S$-indexed set $\{Q_s\colon F\Ra F\circ \Gamma_s\}_{s\in S}$ of $\calV$-natural transformations. Under the equivalence $\Phi$ (which is $S$-equivariant) this datum corresponds to a $\calV$-functor $F'\in\calF'$ together with another such coherent set $\{\Phi(Q_s)\}_{s\in S}$of $\calV$-natural transformations which then assemble to a unique $\calV$-functor in $\calF_S'$.\\
We have thus described the functor $\Phi_S$ on objects; the description on morphisms is analogous. It is clear that the assignment $\Phi\mapsto\Phi_S$ is compatible with composition.
\end{Prf}

\subsubsection{Wreath products}\label{appendixwreathproducts}
We are interested in the case where $S\leqslant  S_n$ is a subgroup of the symmetric group with the canonical action $\Gamma$ on the category $\cal C^n=\cal C\times \dots\times \cal C$ by permuting the factors, i.e.\ $\Gamma_\sigma\colon (f_1,\dots,f_n)\mapsto (f_{\sigma^{-1}(1)},\dots,f_{\sigma^{-1}(n)})$ for $\sigma\in S$. In this case we denote the category $\cal C^n\rtimes_\Gamma S$ by $\cal C\wr S$ and call it the \textbf{wreath product} of $\cal C$ by $S$. If $S\leqslant  T\leqslant S_n$ are nested subgroups then we obtain an inclusion of categories $\cal C^n=\cal C\wr\triv\hra\cal C\wr S\hra\cal C\wr T\hra{\calC\wr S_n}$.\\
Of course we can generalize to the $\calV$-enriched setting and define the wreath product of a $\calV$-category $\calC$ with $S$ as the semidirect product $\calC^\on\wr S$ where $S$ acts on the $n$-fold tensor product of $\calV$-categories by permuting the factors.\\

The diagonal embedding $\Delta\colon \cal C\hra\cal C^n$ is invariant under the action of $S$, i.e.\ $\Gamma_\sigma\circ \Delta=\Delta$ for all $\sigma\in S$. This means that $\Delta$ is equivariant if we equip $\cal C$ with the trivial action of $S$, hence we obtain a functor $\cal C\times \B S\to \cal C\wr S$ which we also call $\Delta$.\\

\ifx\isstandalone\undefined
\else
\newpage
\bibliographystyle{amsalpha}
\bibliography{literatur}
\end{document}
\fi

\subsection{Calculus of canonical mates}\label{calculusofmates}
\ifx\preambleloaded\undefined
   
   \def\isstandalone{}
\fi

	We work in a $2$-category $\calC$. We will mostly be interested in the case $\calC=\CAT$ and we will use the terminology from $\CAT$ (e.g.\ map or morphism for a $1$-cell, natural transformation (isomorphism) for an (invertible) $2$-cell, etc.) where it makes sense.\\
	We collect some abstract statements. All the proofs are purely formal manipulations and are omitted.
See Groth's book~\cite[\S 8.1]{GrothBookProjectOnDerivatorsVolumeI} for more details.
	\begin{Def} Let $f\colon x\to y$ be a map. A map $f_\shriek\colon y\to x$ in the other direction together with natural transformations $f_\shriek\circ f\Ra\Id$ (called \textbf{counit}) and $\Id\Ra f\circ f_\shriek$ (called \textbf{unit}) is called a \textbf{left adjoint} of $f$ if the \textbf{triangle equations} are satisfied, i.e.\ the composite transformations $f = f\circ \Id \La f\circ f_\shriek\circ  f\La \Id \circ f$ and $f_\shriek = f_\shriek\circ \Id \Ra f_\shriek\circ f\circ  f_\shriek\Ra \Id \circ f_\shriek$ are the identity.\\
	The map $f_\shriek$ is called an adjoint inverse, if the unit and counit are both natural isomorphisms. In this case we denote it by $f^\inv$.
	\end{Def}
	\begin{Rem}
		We will abuse notation and call $f_\shriek$ the left adjoint; the unit and counit will be carried along silently and sometimes decorated in a way to make clear that they belong to $f$.\\
		If we say that a map $f$ admits a left adjoint, it shall always be understood that we are thereby also fixing such a left adjoint and denoting it with the standard notation $f_\shriek$ (with the unit and counit implicit).\\
		When drawing diagrams we decorate a map with $\hasadj$ if it admits a left adjoint and by $\approxeq$ if it admits an adjoint inverse.
	\end{Rem}

Consider the following $2$-cells where all arrows admit left adjoints.
	\begin{equation}\label{commtriangleformates}
	\begin{tikzcd}
		\bullet\rar{g}\ar[rr, bend right=50, "h"',""{name=A}]&	\bullet\rar{f}\ar[from=A,Rightarrow, "\alpha"]&	\bullet&&		\bullet\rar{g}\ar[rr, bend right=50, "h"',""{name=A}]&	\bullet\rar{f}\ar[from=A,Leftarrow, "\beta"]&	\bullet
	\end{tikzcd}
	\end{equation}
we obtain natural transformations
	\begin{equation}\label{compositionofinverses}
	\begin{aligned}
	&\alpha_\shriek\colon & g_\shriek\circ f_\shriek\xRa{h} g_\shriek\circ f_\shriek\circ h\circ h_\shriek \xRa{\alpha}g_\shriek\circ f_\shriek\circ f\circ g\circ h_\shriek\xRa{f}g_\shriek\circ g\circ h_\shriek\xRa{g} h_\shriek\\
	&\beta_\shriek\colon & g_\shriek\circ f_\shriek\xLa{h} h_\shriek\circ h\circ g_\shriek\circ f_\shriek\xLa{\beta}h_\shriek\circ f\circ g\circ g_\shriek\circ f_\shriek \xLa{g}h_\shriek\circ  f\circ f_\shriek\xLa{f}h_\shriek
	\end{aligned}
	\end{equation}
	Next consider the $2$-cell
\begin{equation}\label{commutativediagramformate}
\begin{tikzcd}
	\bullet\rar{a^1}\ar[dr,"\gamma",Rightarrow]&		\bullet\\
	\bullet\ar[r,"a^2"']\ar[u,"f^1"', "\hasadj"]&		\bullet\ar[u,"f^2"', "\hasadj"]\\
\end{tikzcd},
\end{equation}
where we assume $f^1$ and $f^2$ to admit left adjoints. We obtain a new square
\begin{equation}\label{mateofsquare}
\begin{tikzcd}
	\bullet\rar{a^1}\ar[d,"f^1_\shriek"']&		\bullet\ar[d,"f^2_\shriek"]\\
	\bullet\ar[r,"a^2"']\ar[ur, Leftarrow, "\gamma_\shriek"]&		\bullet
\end{tikzcd}
\end{equation}
inhabited by the following natural isomorphism:
\begin{equation}
\begin{aligned}
&\gamma_\shriek\colon& f^2_\shriek\circ a^1\xRa{f^1}f^2_\shriek\circ a^1\circ f^1\circ f^1_\shriek\xRa{\gamma}f^2_\shriek\circ  f^2\circ a^2\circ f^1_\shriek\xRa{f^2}a^2\circ f^1_\shriek
\end{aligned}
\end{equation}

We call the natural transformations $\alpha_\shriek$, $\beta_\shriek$ and $\gamma_\shriek$ the \textbf{canonical mates} of the $2$-cells $\alpha$, $\beta$ and $\gamma$ respectively.

\begin{Lem}
	Canonical mates are compatible with pasting of $2$-cells.
	\if{false}\begin{enumerate}
		\item Let $\bullet \xra{h}\bullet\xra{g}\bullet \xra{f}\bullet$ be a composable chain of arrows such that all possible partial composites admit an adjoint inverse. Then the following two pastings of natural isomorphisms give the same result:\\
		\[\begin{tikzcd}[column sep=large]
		|[alias=TL]|\bullet&&	\bullet\ar[ll, "h_\shriek"']&	&	\bullet\ar[from=drr, "(gh)_\shriek" description, ""{name=C, below, near end}, ""{name=D, above, near start}]&&	\bullet\ar[ll, "h_\shriek"']\ar[from=D,Leftrightarrow]\\
		\bullet\ar[u, "(hh)_\shriek"]\ar[rr, "f_\shriek"']\ar[urr, "h_\shriek" description, ""{name=A,  below, near end}, ""{name=B,  above, near start}]&&	\bullet\ar[u, "g_\shriek"']\ar[from=A, Leftrightarrow]&&	\bullet\ar[u, "(f(gh))_\shriek"]\ar[rr, "f_\shriek"']\ar[to=C,Leftrightarrow]&&	\bullet\ar[u, "g_\shriek"']
		\ar[from=TL, to=B, Leftrightarrow]
		\end{tikzcd}\]
		In other words, the natural isomorphisms
		\begin{align*}
		h_\shriek\circ (g_\shriek\circ f_\shriek) &\LRa h_\shriek\circ h_\shriek \LRa (hh)_\shriek \text{ and}\\
		(h_\shriek\circ g_\shriek)\circ f_\shriek &\LRa (gh)_\shriek\circ f_\shriek \LRa (f(gh))_\shriek
		\end{align*} agree.
		\item If we have a pasting of $2$-cells
			\[\begin{tikzcd}
				\bullet\rar{a^1}\ar[rd, Rightarrow]&	\bullet\rar{b^1}\ar[rd, Rightarrow]&	\bullet\\
				\bullet\rar{a^2}\ar[u,"f^1"', "\hasadj"]&	\bullet\ar[u,"f^2"', "\hasadj"]\rar{b^2}&	\bullet\ar[u,"f^3"', "\hasadj"]
			\end{tikzcd}\]
			where the vertical arrows admit an adjoint inverse then the canonical mates of the two individual cells can be pasted to the obtain the mate of the pasted cell.
		\item Consider the commutative diagram
			\[\begin{tikzcd}
				\bullet\rar{a^1}&	\bullet\\
				\bullet\rar{a^2}\ar[u,"f^1"', "\hasadj"]&	\bullet\ar[u,"f^2"', "\hasadj"]\\
				\bullet\rar{a^3}\ar[u,"g^1"', "\hasadj"]&	\bullet\ar[u,"g^2"', "\hasadj"]
			\end{tikzcd}\]
		where the vertical arrows admit an adjoint inverse. The mate of the outer rectangle can be computed by the following pasting:
			\[\begin{tikzcd}
				\bullet\rar{a^1}\ar[d,"f^1_\shriek"']\ar[dd, ""{name=L, right},bend right=90, "(f^1g^1)_\shriek"']&	\bullet\dar{f^2_\shriek}\ar[dd, ""{name=R, left},bend left=90, "(f^2g^2)_\shriek"]\\
				\bullet\rar{a^2}\ar[ur, Leftrightarrow]\ar[d, "g^1_\shriek"']\ar[Leftrightarrow, from=L]&	\bullet\dar{g^2_\shriek}\ar[Leftrightarrow, from=R]\\
				\bullet\rar{a^3}\ar[ur, Leftrightarrow]&	\bullet
			\end{tikzcd}\]
	\end{enumerate}\fi
\end{Lem}

\begin{dCor}
	If the natural transformations $\alpha$ and $\beta$ in diagram \ref{commtriangleformates} are inverse to each other then so are their canonical mates $\alpha_\shriek$ and $\beta_\shriek$.
\end{dCor}

\begin{War}Even if the $2$-cell $\gamma$ in diagram \ref{commutativediagramformate} is invertible, it is usually not true that the canonical mate $\gamma_\shriek$ is an isomorphism. For example, the mate of the identity cells
	\begin{equation}
		\begin{tikzcd}
			\bullet\rar{\Id}&\bullet&&			\bullet\rar{f}&\bullet\\
			\bullet\ar[ur, phantom, bend right, ""{name=B}]\ar[ur, phantom, bend left, ""{name=A}]\ar[u, "f"', "\hasadj"]\rar{f}&\bullet\uar{\Id}&&		\bullet\ar[ur, phantom, bend right, ""{name=D}]\ar[ur, phantom, bend left, ""{name=C}]\uar{\Id}\rar{\Id}&\bullet\ar[u, "f"', "\hasadj"]
			\ar[from=A,to=B, equal]
			\ar[from=C,to=D, equal]
		\end{tikzcd}
	\end{equation}
	are nothing but the unit $\Id\Ra f\circ f_\shriek$ and the counit $f_\shriek\circ f\Ra\Id$ of the adjunction.
\end{War}

As the next Lemma shows, the above example is the only obstruction.

\begin{Lem}
	Assume the maps $f^1$ and $f^2$ in the diagram \ref{commutativediagramformate} admit not only a left adjoint but an adjoint inverse. Then the $2$-cell $\gamma$ is invertible if and only if its canonical mate $\gamma_\shriek$ is an isomorphism.
\end{Lem}

\ifx\isstandalone\undefined
\else
\newpage
\bibliographystyle{amsalpha}
\bibliography{literatur}
\end{document}
\fi

\subsection{Some enriched category theory}\label{appendixenrichedcattheory}
\ifx\preambleloaded\undefined
   
   \def\isstandalone{}
\fi

We will briefly list some lemmata which will be used in Section~\ref{Part2ofbigequiv}. For a systematic theory of enriched category theory, see Kelly's book~\cite{KellyBasicConceptsOfEnrichedCategoryTheory}.\\
We work over a complete and cocomplete closed symmetric monoidal category $\calV$.\\

The Yoneda principle is: if one can build two (potentially different) functors out of instances of the Yoneda-maps $Y^D\colon \cal D\hra[\calD,\calV]_\op$ and $Y_\calD\colon \calD\hra[\calD_\op,\calV]$, evaluation functor $[\calD,\calV]\otimes\calD\to\calV$ and the Hom-functor $\Hom_\calD\colon \calD_\op\otimes\calD\to\calV$ then those two functors will be $\calV$-naturally isomorphic.
\begin{Lem}[Enriched Yoneda]\label{YonedaLemma}~\cite[5.1]{DayKellyEnrichedFunctorCategories}
	\begin{enumerate}[label=\roman*.]
	\item\label{classicalYoneda} There is an isomorphism $\rho(d)\cong \Mor{[\calD,\calV]}{\Mor\calD d\blank}\rho$, natural in $d\in \calD$ and $\rho\in [\calD,\calV]$. In other words, the functor
	\[\calD\otimes[\calD,\calV]\xra{Y^\calD}[\calD,\calV]_\op\otimes[\calD,\calV]\xra{\Hom_{[\calD,\calV]}}\calV\]
	is naturally isomorphic to the evalutation map.
	\item\label{variantYoneda} If $L\colon [\calD,\calV]_\op\to V$ is representable, then a representing object can be chosen to be $\restr L\calD\colon \cal D\xhra{Y^\calD}[\calD,\calV]_\op\xra{L}\calV$. In other words the functor
	\[[\calD,\calV] \xhra{Y_{[\calD,\calV]}}[[\calD,\calV]_\op,\calV]\xra{[Y^\calD,\calV]} [\calD,\calV]\]
	is naturally isomorphic to the identity functor.\qedhere
	\end{enumerate}
\end{Lem}

\begin{Lem}\label{partialhomtensoradjunction}
	Let $\rho_1\colon \calD_1\to\calV$, $\rho_2\colon \calD_2\to\calV$ and $F\colon \calD_1\otimes\calD_2\to\calV$ be $\calV$-functors. Then we have a canonical isomorphism of objects in $\calV$
	\[\Mor{[\calD_1,\calV]}{\rho_1}{d_1\mapsto\Mor{[\calD_2,\calV]}{\rho_2}{F(d_1,\blank)}}\xra{\cong}\Mor{[\calD_1\otimes\calD_2,\calV]}{\rho_1\otimes\rho_2}F,\]
	which is natural in $\rho_1\in [\calD_1,\calV]$, $\rho_2\in [\calD_2,\calV]$ and $F\in [\calD_1\otimes\calD_2,\calV]$.
\end{Lem}
\begin{Prf}
	We calculate:
	\begin{align*}
		\LHS &=\int_{d_1\in \calD_1}\left[\rho_1(d_1),\hom{[\calD_2,\calV]}(\rho_2,F(d_1,\blank))\right]\\
		&=\int_{d_1\in D_1}\left[\rho_1(d_1),\int_{d_2\in \calD_2}[\rho_2(d_2),F(d_1,d_2)]\right]\\
		&\cong \int_{d_1\in \calD_1}\int_{d_2\in \calD_2}\left[\rho_1(d_1),[\rho_2,F(d_1,d_2)]\right]\\
		&\cong \int_{(d_1,d_2)\in\calD_1\otimes\calD_2}\left[\rho_1(d_1)\otimes\rho_2(d_2),F(d_1,d_2)\right] = \RHS
	\end{align*}
	where the first $\cong$ comes from the fact that ends commute with $[\blank,\blank]$ in the second variable  (by definition); the second is the \qquote{Fubini theorem for ends}~\cite[2.1]{KellyBasicConceptsOfEnrichedCategoryTheory} to join the two ends and the Hom-Tensor adjunction in $\calV$ to modify the argument. 
\end{Prf}

\ifx\isstandalone\undefined
\else
\newpage
\bibliographystyle{amsalpha}
\bibliography{literatur}
\end{document}
\fi

\end{appendices}


\newpage
\bibliographystyle{amsalpha}
\bibliography{literatur}


\end{document}